\documentclass[12pt]{article}
\usepackage{amsmath}
\usepackage{amssymb}
\usepackage{amsthm}
\usepackage[totalwidth=17cm, totalheight=25cm]{geometry}
\usepackage{graphicx}
\usepackage{mathabx}
\usepackage{tabularx}
	\newcolumntype{C}[1]{>{\centering\arraybackslash}m{#1}} 
	\newcolumntype{R}[1]{>{\raggedleft\arraybackslash}m{#1}} 
\usepackage{varwidth}

\usepackage{color}
\usepackage{diagrams}

\newtheoremstyle{boldplain}% name
{9pt}%      Space above
{9pt}%      Space below
{\itshape}%         Body font
{}%         Indent amount (empty = no indent, \parindent = para indent)
{\bfseries}% Thm head font
{.}%        Punctuation after thm head
{.5em}%     Space after thm head: " " = normal interword space;
{\thmname{#1}\thmnumber{ #2}\thmnote{ (#3)}}%

\newtheoremstyle{bolddefinition}% name
{9pt}%      Space above
{9pt}%      Space below
{}%         Body font
{}%         Indent amount (empty = no indent, \parindent = para indent)
{\bfseries}% Thm head font
{.}%        Punctuation after thm head
{.5em}%     Space after thm head: " " = normal interword space;
{\thmname{#1}\thmnumber{ #2}\thmnote{ (#3)}}%

\theoremstyle{boldplain}
\newtheorem{add}[equation]{Addendum}

\newtheorem{ass}[equation]{Assumption}
\newtheorem{claim}[equation]{Claim}

\newtheorem{cor}[equation]{Corollary}

\newtheorem{lem}[equation]{Lemma}
\newtheorem{lemma}[equation]{Lemma}
\newtheorem{prop}[equation]{Proposition}

\newtheorem{thm}[equation]{Theorem}
\newtheorem{theorem}[equation]{Theorem}

\theoremstyle{bolddefinition}
\newtheorem{dfn}[equation]{Definition}
\newtheorem{definition}[equation]{Definition}

\newtheorem{example}[equation]{Example}
\newtheorem{exa}[equation]{Example}

\newtheorem{rem}[equation]{Remark}

\newfont{\bigbf}{cmbx10 scaled\magstep1}
\normalbaselines
\numberwithin{equation}{section}

\def\B{{\mathcal B}}

\def\R{{\mathbb R}}
\def\H{{\mathbb H}}

\def\N{{\mathbb N}}
\def\P{\mathcal P}
\def\Z{{\mathbb Z}}
\def\Q{{\mathbb Q}}

\def\al{\alpha}
\def\ga{\gamma}
\def\Ga{\Gamma}
\def\de{\delta}
\def\De{\Delta}
\def\eps{\epsilon}
\def\la{\lambda}
\def\La{\Lambda}
\def\si{\sigma}
\def\Si{\Sigma}

\def\3{\ss}

\def\acts{\curvearrowright}
\def\amod{a_{mod}}

\def\CH{\operatorname{CH}}
\def\QCH{\operatorname{QCH}}

\def\D{\partial}

\def\diamot{\diamondsuit_{\tau_{mod}}}
\def\diamoTh{\diamondsuit_{\Theta}}

\def\Dt{\partial_{\tau_{mod}}}
\def\embed{\hookrightarrow}
\def\Fix{\operatorname{Fix}}
\def\Flag{\operatorname{Flag}}

\def\Flags{\operatorname{Flag_{\si_{mod}}}}

\def\Flagt{\operatorname{Flag_{\tau_{mod}}}}

\def\Fmod{F_{mod}}

\def\geo{\partial_{\infty}}
\def\geoc{\partial_{\infty}^{con}}

\def\geot{\partial_{\infty}^{\tau_{mod}}}
\def\geotc{\partial_{\infty}^{\tau_{mod},con}}

\def\id{\operatorname{id}}
\def\im{\operatorname{im}}

\def\inte{\operatorname{int}}

\def\Isom{\mathop{\hbox{Isom}}}
\def\Lac{\Lambda^{con}}
\def\LaGa{\Lambda(\Gamma)}
\def\Las{\Lambda_{\sigma_{mod}}}
\def\LasGa{\Lambda_{\sigma_{mod}}(\Gamma)}
\def\LatGa{\Lambda_{\tau_{mod}}(\Gamma)}
\def\Lat{\La_{\tau_{mod}}}
\def\LaXt{\Lambda_{X,\tau_{mod}}}
\def\LaXtc{\Lambda_{X,\tau_{mod}}^{con}}

\def\LaY{\Lambda_Y}
\def\LaYc{\Lambda_Y^{con}}

\def\lra{\longrightarrow}

\def\Min{\operatorname{Min}}

\def\ol{\overline}

\def\pihalf{\frac{\pi}{2}}

\def\2pithird{\frac{2\pi}{3}}

\def\rank{\mathop{\hbox{rank}}}

\def\Ra{\Rightarrow}

\def\simod{\si_{mod}}

\def\st{\operatorname{st}}

\def\ost{\operatorname{ost}}

\def\Stab{\operatorname{Stab}}

\def\taumod{\tau_{mod}}

\def\8{\infty}
\def\<{\langle}
\def\>{\rangle}
\def\ov{\overrightarrow}

\def\mini{\scriptsize}

\long\def\comment#1\endcomment{}

\begin{document}

\title{Relativizing characterizations of Anosov subgroups, I} 
\author{Michael Kapovich, Bernhard Leeb}
\date{\today}
\maketitle

\begin{abstract}
We propose several common extensions of the classes of Anosov subgroups and geometrically finite Kleinian groups
among discrete subgroups of semisimple Lie groups. 
We relativize various dynamical and coarse geometric characterizations of Anosov subgroups
given in our earlier work,
extending the class from intrinsically hyperbolic to relatively hyperbolic subgroups. 
We prove implications and equivalences between the various relativizations.
\end{abstract}

\tableofcontents

\section{Introduction}

The notion of geometric finiteness was first introduced by Ahlfors \cite{Ahlfors} 
in the context of Kleinian group actions on hyperbolic 3-space $\H^3$.
It was originally defined via the existence of finite-sided convex fundamental polyhedra. 
A few years later, Beardon and Maskit \cite{BM} 
gave a dynamical characterization of geometrically finite groups in terms the action on the limit set,
now called the {\em Beardon-Maskit condition}. 
Subsequently, alternative characterizations were given by Marden \cite{Marden}, Thurston \cite{Thurston} and many others,
see e.g.  \cite{Bowditch93} and \cite{Rat}.

While Ahlfors' original definition turned out to be unsuitable for hyperbolic space of dimension $\ge 4$, 
the Beardon-Maskit condition worked well in the context of discrete subgroups of rank one Lie groups and, 
more generally, of discrete groups of isometries acting on negatively pinched Hadamard manifolds \cite{Bowditch},
and was shown to be equivalent to a variety of other properties.
The Beardon-Maskit condition remains meaningful even in the purely dynamical setting of 
convergence actions on topological spaces,
something which we are exploiting in our work.

A particularly nice subclass of geometrically finite Kleinian groups is formed by convex cocompact subgroups
which are distinguished by the absence of parabolic elements. 
They are intrinsically word hyperbolic,
whereas a general geometrically finite Kleinian group inherits a natural structure as a are relatively hyperbolic group,
the peripheral structure given by the collection of maximal parabolic subgroups.

The notion of convex cocompact Kleinian groups was extended to discrete subgroups of 
higher rank Lie groups, starting with the notion of {\em Anosov subgroups} \cite{Labourie}, see also \cite{GW}. 
These were originally defined in terms of their dynamics on flag manifolds.
We subsequently gave various characterizations of Anosov subgroups 
in terms of their coarse geometry, dynamics and topology
along with a simplification of their original definition
\cite{morse,mlem,bordif, anolec}, see also \cite{anosov,manicures}.

As convex cocompact subgroups, Anosov subgroups are intrinsically word hyperbolic,
and as the former contain no parabolics, the latter contain no strictly parabolic elements,\footnote{That is, 
non-elliptic elements with zero infimal displacement.} e.g. no unipotents.
Our goal is to find a common extension of the classes of Anosov and geometrically finite subgroups,
that is, to complete the diagram:

$$
 \begin{diagram}
\framebox{convex cocompact} &   &	&	     	\rTo^{\hbox{allow parabolics}}			&       &             &\framebox{geometrically finite}  \\
 	\dTo^{\hbox{higher rank}}&	& 	&		& &	 &\dTo_{\hbox{higher rank}}		  \\
\framebox{~~~~~~~Anosov~~~~~~~} 	& &	 &		\rTo_{\hbox{}}	& &	&  \framebox{~~~~~~~~~~~~?~~~~~~~~~~~~} \\ 
 \end{diagram}
 $$
 
\addvspace{0.2cm}

\noindent
We consider subgroups which are relatively hyperbolic as abstract groups
and extend to this more general setting 
various characterizations of Anosov subgroups studied in our earlier papers.

The Beardon-Maskit condition has the most straightforward generalization,
namely by requiring that the discrete subgroup acts on its limit set like a relatively hyperbolic group.
This leads to the notions {\em relatively asymptotically embedded} and {\em relatively RCA},
see Definitions~\ref{def:relasembedded} and~\ref{dfn:rlrca},
which are equivalent in view of Yaman's dynamical characterization of relatively hyperbolic groups.
The intrinsic relatively hyperbolic structure of these discrete subgroups 
can be read off the dynamics on the limit set and is therefore uniquely determined. 

Also our coarse geometric characterization of Anosov subgroups as Morse subgroups,
that intrinsic geodesics in the subgroup are extrinsically perturbations of Finsler geodesics in the symmetric space,
generalizes naturally.
This leads to the notions of {\em relatively Morse} and {\em relatively Finsler-straight},
see Definitions~\ref{dfn:rlmrs} and~\ref{dfn:fstrghtctn}.

All these relative notions agree in rank one with geometric finiteness (see Corollary \ref{cor:rank1}).
The main result of the paper establishes relations (implications and equivalences) 
between them 
in higher rank.
It is summarized in Theorem~\ref{thm:main} and the diagram:

\medskip 
$$
 \begin{diagram}
 		                                              &         &\framebox{rel Morse} &              &                \\
 					    		    &\ldTo &				 &   \rdTo   & 																						 \\      								                    
\framebox{ \parbox{5.2cm} {\centering relatively uniformly\\ Finsler-straight  }}&        &\pile{\lTo\\ \rTo}             &              &\framebox{ \parbox{8cm} {\centering rel  asymptotically embedded with\\   uniformly regular peripheral subgroups  }}\\
 	\dTo                                                &	      &   				 &	        &\dTo		                                                                                                            							\\	
	\framebox{rel  Finsler-straight}    &\pile{\lTo\\ \rTo}&\framebox{rel RCA}&\pile{\lTo\\ \rTo}       &\framebox{  \parbox{3.25cm}{\centering rel  asymptotically \\embedded}}                                                                     \\
		                                            &	      &					&		&	\dTo	\uTo_{\hbox{if Zariski dense}}                                                      								 \\
	  						 &        &             			&		&	\framebox{rel  boundary embedded}                                                     									 \\                         
 \end{diagram}
 $$

The most difficult implications are 
between {\em relatively Finsler straight} and {\em relatively asymptotically embedded},
connecting coarse geometry and dynamics,
and their analogues in the uniformly regular case.
They are proven in section \ref{sec:str=ae}, 
which in turn relies on coarse geometric results about general (non-equivariant) Finsler straight maps
established in section \ref{sec:maps}. 

Examples of classes of discrete subgroups satisfying the relative conditions 
discussed in the paper are:

\begin{enumerate}
\item subgroups preserving a rank one symmetric subspace and acting on it in a geometrically finite fashion
(Theorem \ref{thm:RM=GF} and Example \ref{ex:geodembedding})

\item discrete groups of projective transformations acting with finite covolume on strictly convex solids in $\R^n$ 
(studied in \cite{CLT}) 

\item certain families of discrete subgroups of $PGL(3,\R)$ not preserving properly convex domains in $\R P^2$
(described in \cite{Sch} and \cite{Kim-Lee})

\item positive representations (into split semisimple Lie groups) of fundamental groups of punctured surfaces
(appearing in \cite{FG}) 

\item certain free products of opposite uniformly $\taumod$-regular elementary unipotent subgroups (see \cite{relmorse-2})

\item small relative deformations 
(see \cite{relmorse-2})
\end{enumerate}

The discussion of the relative notions introduced in this paper 
will be continued in \cite{relmorse-2}. {In particular, we will prove several forms of stability for 
subgroups satisfying these notions and combination theorems. We will also prove that relative uniform Finsler straightness is equivalent to the condition of {\em uniform regularity and relative non-distortion}: Partial results in this direction are obtained by Feng Zhu in \cite{Zhu}. Some things are unclear to us at this point, for instance: }

{1. Are the different relativizations of the Anosov condition (including {\em uniform regularity and relative non-distortion}) genuinely different or result in the same class of discrete subgroups?}  

{2. If these notions are different, which are the most important ones?}

\medskip 
{In view of this, we decided to discuss all the different notions in the present paper.}

\medskip  
{\bf Acknowledgements.} The first author was partly supported by the NSF grant  DMS-16-04241, by 
KIAS (the Korea Institute for Advanced Study) through the KIAS scholar program, by a Simons Foundation Fellowship, grant number 391602, 
and by Max Plank Institute for Mathematics in Bonn. Much of this work was done during our stay at KIAS and Oberwolfach and we are thankful to KIAS and Oberwolfach (MFO) for their  hospitality. We are grateful to Grisha Soifer for providing us with a reference to Prasad's work \cite{Prasad} 
and for writing the appendix to this paper. We are also grateful to Sungwoon Kim   and Jaejeong Lee for helpful conversations.

\section{Preliminaries and notation}

\subsection{Metric spaces}\label{sec:metric_spaces}

We will be using the notation $xy$ for a geodesic segment in a metric space connecting points $x$ and $y$. 
Similarly, in a geodesic metric space $Y$ which is Gromov hyperbolic or CAT(0),
we will use the notation $y\xi$ for a geodesic ray in $Y$ emanating from $y$ and asymptotic to a point $\xi$
in the visual (ideal) boundary $\geo Y$ of $Y$. For two distinct ideal boundary 
points $\eta_\pm\in \geo Y$ of a Gromov-hyperbolic space $Y$ we will use the notation $\eta_-\eta_+$ for a geodesic in $Y$ asymptotic to $\eta_\pm$.  A similar notation will be used for {\em Finsler geodesics} in symmetric spaces:  
 $x\tau$, $\tau_- \tau_+$ will denote a Finsler geodesic ray/line; see section \ref{sec:finsler}. 

We will use the notation $B(a, R)$ for the open $R$-ball with center $a$ in a metric space,
and the notation $N_R(A)$ for the open $R$-neighborhood of subsets $A$, where $R >0$. 
The subsets $N_R(A)$ are called {\em tubular neighborhoods} of $A$.

A metric space is called {\em taut} 
if every point lies at distance $\leq R$ from a geodesic line for some uniform constant $R$.

Two subsets in a metric space are called {\em $D$-separated} if their infimal distance is $\geq D$. 

We call a subset of a metric space {\em $s$-spaced} if its distinct points have pairwise distance $\geq s$,
and we call a map into a metric space {\em $s$-spaced} if it is injective and its image is $s$-spaced.

A sequence $(x_n)$ in a metric space is said to {\em diverge to infinity} if $\lim_{n\to\infty} d(x_1, x_n)= \infty$; 
we will refer to such $(x_n)$ as a {\em divergent sequence}.  

A map between metric spaces is called {\em metrically proper} if it sends divergent sequences to divergent sequences,
equivalently, if the preimages of bounded subsets are bounded.

\subsection{Group actions}

For an action $\Ga\acts X$ of a group $\Ga$ on a set $X$ 
we let $\Ga_x<\Ga$ denote the stabilizer of an element $x\in X$.
The associated orbit map is defined by 
$$ o_x: \Ga\to X, \quad \ga\mapsto\ga x. $$
If $\Ga\acts Y$ is another $\Ga$-action and $y\in Y$ is a point such that $\Ga_y\leq\Ga_x$,
then there is a well-defined $\Ga$-equivariant map of orbits
$$ o_{x,y}: \Ga y \to \Ga x , \quad \ga y \mapsto \ga x .$$

\subsection{Convergence actions}
\label{sec:conv}

A continuous action $\Ga\acts Z$ 
of a discrete group $\Ga$ 
on a compact metrizable topological space $Z$ 
is called a {\em (discrete) convergence action} 
if for each sequence $(\ga_n)$ of pairwise distinct elements in $\Ga$ 
there exists a pair of  points $z_-, z_+\in Z$ such that, after extraction, 
the sequence $(\ga_n)$ converges to $z_+$ uniformly on compacts in $Z- \{z_-\}$. 
Note that all actions on spaces with at most two points are convergence,
except actions of infinite groups on the empty space.
Also, all actions of finite groups are convergence. 

The {\em limit set} $\La=\LaGa\subseteq Z$ consists of all points which occur as such limits $z_+$.
The limit set is $\Ga$-invariant and compact.
If $|\La|\geq3$, then it is perfect and the action $\Ga\acts\La$ has finite kernel and is minimal.\footnote{I.e. every orbit is dense.}
If $|\Ga|=+\infty$, then $\La\neq\emptyset$,
and if $|\Ga|<+\infty$, then $\La=\emptyset$.

Elements of convergence groups fall into three classes:
An element is called {\em hyperbolic} if it has infinite order and exactly two fixed points,
{\em parabolic} if it has infinite order and exactly one fixed point,
and {\em elliptic} if it has finite order. 

A point $z\in Z$ is called a {\em parabolic fixed point} of $\Ga$ 
if it is the fixed point of some parabolic element in $\Ga$.
It then is a limit point of its stabilizer $\Ga_z$,
and it turns out that in fact $\La(\Ga_z)=\{z\}$, see \cite[Lemma 2F]{Tukia1994}.

The following types of limit points will be important for this paper
(given the nature of the actions of relatively hyperbolic groups on their ideal boundaries):

\begin{dfn}
\label{dfn:conlimbddpar}
A point $z\in \LaGa$ is called a 

(i) 
{\em conical limit point} for $\Ga$ 
if there exists a sequence $(\ga_n)$ of distinct elements in $\Ga$ 
and a point $w\in \La-\{z\}$ 
such that the sequence of pairs $(\ga_n^{-1} z, \ga_n^{-1} w)$ 
does not accumulate at the diagonal of $Z\times Z$. 

(ii) 
{\em bounded parabolic point} of $\Ga$ 
if its stabilizer $\Ga_z<\Ga$ acts on $\La(\Ga)- \{z\}$ properly discontinuously and cocompactly. 

(ii') 
{\em bounded parabolic fixed point} of $\Ga$ 
if it is both a bounded parabolic point and a parabolic fixed point.
\end{dfn}

Note that the stabilizer of a bounded parabolic point $z$ is necessarily infinite and $\La(\Ga_z)=\{z\}$.
Property (ii') is strictly stronger than (ii) 
because the stabilizer of a bounded parabolic point can be an infinite torsion group.

If $|\La(\Ga)|=1$ and $\Ga$ is not a torsion group,
then the limit point is a bounded parabolic fixed point
and not a conical limit point.
If $|\La(\Ga)|=2$, then both limit points are conical and not bounded parabolic.

The importance of convergence actions in our work is due primarily to two reasons:

\begin{itemize}
\item If $Y$ is a proper geodesic Gromov-hyperbolic metric space and $\Ga$ is a discrete subgroup of the isometry group of $Y$, then the natural action of $\Ga$  on the visual boundary $\geo Y$ of $Y$ is a convergence action, see \cite{Tukia1994}. 

\item If $\Ga$ is a $\taumod$-regular antipodal subgroup of the isometry group of a symmetric space of noncompact type, then the natural action of $\Ga$ on its $\taumod$-limit set in the flag-manifold 
$\Flagt$ is a convergence action, see \cite{anolec}. We  note that $\taumod$-RCA (regular antipodal and conical) subgroups are precisely the $\taumod$-Anosov subgroups. 

\end{itemize} 

We refer the reader to \cite{Bowditch_config, Tukia1994, Tukia_conical}  for in-depth discussions of convergence actions. {In particular, an interested reader can find there a 
 proof of the fact that, in the case when $Z$ 
is the visual boundary of a hyperbolic space $Y$,  the above definition of conical limit points is equivalent to the one in terms of conical subconvergence of $\Ga$-orbits 
in $Y$, see section \ref{sec:HypSpaces} for the definition of the conlical convergence. 
We also note that this proof is essentially the same as the one given for the hyperbolic 3-space by Beardon and Maskit in their pioneering paper \cite[Theorem 1]{BM}.}

\section{Some coarse hyperbolic geometry}
\subsection{Gromov hyperbolic spaces}\label{sec:HypSpaces}

Background material on hyperbolic spaces
can be found in \cite{BH}, \cite{Bowditch2012}, \cite{Bourdon}, \cite{Drutu-Kapovich} and \cite{Vaisala}.

Let $Y$ be a proper geodesic metric space which is $\de$-hyperbolic in the sense of Gromov 
for some $\de\geq0$.
We denote by $\ol Y=Y\sqcup\geo Y$ its visual compactification.  

Geodesics in $Y$ are roughly unique in the sense that
any two geodesic segments with the same endpoints have Hausdorff distance $\leq C\de$,
where $C$ is a uniform constant 
(depending on the definition of $\de$-hyperbolicity which is used).
The same holds for any two asymptotic geodesic rays with the same initial point,
and for any two (at both ends) asymptotic geodesic lines.

A family of geodesics in $Y$ is {\em bounded} if for some (any) point $y\in Y$ 
the distance of $y$ to the geodesics in this family is uniformly bounded.
The pairs $(\ol y,\ol y')$ of endpoints of geodesics $\ol y\ol y'$ in $Y$ lie in the set 
$(Y\sqcup\geo Y)^2 - \De_{\geo Y}$.
The boundedness of a family of geodesics in a Gromov hyperbolic space 
is an asymptotic property of its set of pairs of endpoints:
\begin{lem}
\label{lem:bddasyprp}
Let $E_Y\subset (Y\sqcup\geo Y)^2  - \De_{\geo Y}$.
Then the family of all geodesics in $Y$ with pair of endpoints in $E_Y$ is bounded
if and only if $E_Y$ is relatively compact in $(Y\sqcup\geo Y)^2 - \De_{\geo Y}$.
\end{lem}
\proof
Any bounded sequence of geodesics $\ol y_n\ol y'_n$ in $Y$ subconverges to a geodesic $\ol y\ol y'$,
and the pairs of endpoints $(\ol y_n,\ol y'_n)$ subconverge to $(\ol y,\ol y') \in (Y\sqcup\geo Y)^2 - \De_{\geo Y}$.
Thus the sequence of pairs $(\ol y_n,\ol y'_n)$ does not accumulate at $\De_{\geo Y}$.

On the other hand,
if a sequence of geodesics $\ol y_n\ol y'_n$ diverges,
i.e. their distances from some base point $y\in Y$ diverge to infinity,
then $\de$-hyperbolicity implies that there exists points $z_n\in y\ol y_n$ and $z'_n\in y\ol y'_n$
such that $z_n,z'_n\to\infty$ and the segments $yz_n$ and $yz'_n$ are $C\de$-Hausdorff close.
It follows that the pairs $(\ol y_n,\ol y'_n)$ accumulate at $\De_{\geo Y}$.
\qed

\medskip
A sequence $(y_n)$ in $Y$ is said to converge to $\eta\in \geo Y$ {\em conically} 
if $y_n\to \eta$ and $(y_n)$ is contained in a tubular neighborhood of a ray asymptotic to $\eta$.
This is independent of the ray since any two asymptotic rays have finite Hausdorff distance. 
For a subset $A\subset Y$,
the {\em conical accumulation set} $\geoc A\subset \geo Y$ 
consists of all points $\eta\in \geo Y$ for which there exists a sequence $(a_n)$ in $A$ converging to $\eta$ conically.

Given a discrete isometric group action $\Ga\acts Y$, we define its {\em limit set} $\La=\LaY$ as the accumulation set $\geo (\Ga y)$ in $\geo Y$ 
of one (equivalently, every) $\Ga$-orbit in $Y$. We will use the notation $\Lac=\LaYc$ for the {\em conical limit set} of this action, i.e. the set $\geoc (\Ga y)$ of conical limit points of the group $\Ga$.

\medskip 
{\em Straight triples.}
We denote by $T(Y):=Y^3$ the space of triples of points in $Y$ and by 
\begin{equation}
\label{eq:idtrpy}
T(Y, \geo Y):=(Y\sqcup\geo Y)\times Y\times (Y\sqcup\geo Y)
\end{equation}
the space of ideal triples in the visual compactification $\ol Y=Y\sqcup\geo Y$
with middle point in $Y$.

We first define straightness for (non-ideal) triples in $Y$:
\begin{definition}[Straight triple]
\label{def:strtrp}
A triple $(y_-,y,y_+)\in T(Y)$ is called {\em $D$-straight}, $D\geq0$, 
if the points $y_-,y$ and $y_+$ are $D$-close to points $y_-',y'$ and $y_+'$, respectively, 
which lie in this order on a geodesic (segment).
\end{definition}

This notion naturally extends to {\em ideal} triples in $T(Y, \geo Y)$:
We say that 
a triple $(y_-,y,\eta_+)\in Y^2\times\geo Y$ is $D$-straight 
if the points $y_-$ and $y$ are $D$-close to points $y_-'$ and $y'$, respectively,
such that $y'$ lies on a geodesic ray $y_-'\eta_+$.
Analogously for triples $(\eta_-,y,y_+)\in\geo Y\times Y^2$.
Similarly,
we say that 
a triple $(\eta_-,y,\eta_+)\in\geo Y\times Y\times\geo Y$ is $D$-straight 
if $\eta_-\neq\eta_+$ and the point $y$ lies within distance $D$ of a geodesic line $\eta_-\eta_+$.

$Y$ is taut, 
if every point $y$ is the middle point of a uniformly straight triple $(\eta_-,y,\eta_+)$.

\medskip 
{\em Straight holey lines.} 
We call a map $q:H\to Y$ from an arbitrary (``holey'') subset of $H\subset\R$ 
a {\em holey line}.
If $H$ has a minimal element, we also call $q$ a {\em holey ray}.
A {\em sequence} $(y_n)_{n\in\N}$ in $Y$ can be regarded as a holey ray $\N\to Y$.

We will consider extensions to infinity $\bar q:\ol H:=H\sqcup\{\pm\infty\}\to\ol Y=Y\sqcup\geo Y$ of holey lines $q:H\to Y$ 
by sending $\pm\infty$ to ideal points $\eta_{\pm}\in\geo Y$,
and refer to $\bar q$ as an {\em extended holey line}.
Similarly,
for holey rays $q:H_0\to Y$, we will consider extensions $\bar q:\ol H_0:=H_0\sqcup\{+\infty\}\to\ol Y$ 
by sending $+\infty$ to an ideal point $\eta\in\geo Y$, 
and refer to $\bar q$ as an {\em extended holey ray}.

We carry over the notion of straightness from triples to holey lines by requiring it for all triples in the image:
\begin{definition}[Straight holey line] 
\label{def:strspm}
A holey line $q:H\to Y$ is called {\em $D$-straight}
if the triples $(q(h_-), q(h), q(h_+))$ in $Y$ are $D$-straight for all $h_-\leq h\leq h_+$ in $H$. 
\end{definition}

Similarly, we say that an extended holey line $\bar q:\ol H\to\ol Y$ is $D$-straight
if the triples $(q(h_-), q(h), q(h_+))$ in $\ol Y$ are $D$-straight for all $h_-\leq h\leq h_+$ in $\ol H$ with $h\in H$,
and analogously in the ray case.

Straight holey lines are up to bounded perturbation monotonic maps into geodesics.
More precisely,
for a $D$-straight holey line $q:H\to Y$
there exists a geodesic $c\subset Y$ and a monotonic map $\bar q:H\to c$
which is $D'(D)$-close to $q$. 
The holey line $\bar q:\ol H\to\ol Y$
extended by $\bar q(\pm\infty)=\eta_{\pm}:=c(\pm\infty)$
is then $D'$-straight.
The ideal points $\eta_{\pm}$ are unique if $q$ is biinfinite.

\medskip 
Let $I\subseteq\R$ be an interval. 
We say that a function $f:I\to\R$ 

(i) 
has {\em $\eps$-coarsely slope $s$} if 
$$|f(t_1)-f(t_2)-s(t_1-t_2)| \leq \eps $$
for all $t_1,t_2\in I$.

(ii) 
is {\em $\eps$-coarsely convex} if 
$$\mu_1 f(t_1)+\mu_2 f(t_2)\leq f(\mu_1t_1+\mu_2t_2) +\eps $$
for all $t_1,t_2\in I$ and all $\mu_1,\mu_2\geq0$ with $\mu_1+\mu_2=1$.

Note: 
If there exists $t_0\in I$ such that $f|_{I\cap(-\infty,t_0]}$ has {\em $\eps$-coarsely slope $-1$ } 
and $f|_{I\cap[t_0,+\infty)}$ has {\em $\eps$-coarsely slope $+1$},
then $f$ is $2\eps$-coarsely convex. 

Transferring these notions,
we say that a function $f:Y\to\R$ has $\eps$-coarsely slope $s$ or is $\eps$-coarsely convex
along a geodesic $c:I\to Y$ 
if the composition $f\circ c$ has this property.

\medskip 
{\em Horofunctions and horoballs.}
Horofunctions coarsely measure relative distances from points at infinity.
They arise most naturally as limits of normalized distance functions.

Fix an ideal point $\eta\in\geo Y$.
Let $(y_n)$ be a sequence in $Y$ such that $y_n\to\eta$.
After passing to a subsequence, the sequence of distance functions $d(\cdot,y_n)$ converges up to additive constants,
i.e. there exists a sequence $a_n\to\infty$ of real numbers and a function $h:Y\to\R$ such that 
$$ d(\cdot,y_n)-a_n \to h $$
locally uniformly.
The function $h$ has the following properties:
It is 1-Lipschitz and for every point $y\in Y$ there exists a ray $\rho_y:[0,\infty)\to Y$ asymptotic to $\eta$ with initial point $y$
along which $h$ decays with slope $\equiv-1$,
i.e. $h\circ\rho_y|_{t_1}^{t_2}=t_1-t_2$ for all $t_1,t_2\geq0$. 
(Such a ray $\rho_y$ arises as a sublimit of the segments $yy_n$.)
For an arbitrary ray $\rho:[0,\infty)\to Y$ asymptotic to $\eta$,
it follows that 
\begin{equation}
\label{eq:hfcdec}
\bigl|(h\circ\rho|_{t_1}^{t_2})-(t_1-t_2)\bigr|\leq C\de
\end{equation} 
for all $t_1,t_2\geq0$ with a uniform constant $C$. 
We define a {\em horofunction at $\eta$}
as a function $h:Y\to\R$ which satisfies \eqref{eq:hfcdec} for all rays $\rho$ asymptotic to $\eta$.
Any two horofunctions $h,h'$ at $\eta$ coarsely differ by an additive constant,
i.e. 
\begin{equation}
\label{eq:hfceq}
|(h(y)-h(y'))-(h'(y)-h'(y'))| = |(h(y)-h'(y))-(h(y')-h'(y'))| 
\leq C\de
\end{equation} 
for all $y,y'\in Y$ (with a possibly different uniform constant $C$).\footnote{We will often use the same letter $C$ for a constant 
with the understanding that the constant may vary from inequality to inequality.}

Horofunctions are uniformly {\em coarsely convex}.
This is a consequence of the following stronger property:
For any horofunction $h$ and segment $zz'$ there 
exists a division point $y_0\in yy'$ such that $h$ has $C\de$-coarsely slope $1$ along the oriented segments $y_0y$ and $y_0y'$
with a uniform constant $C$.

We define horoballs as coarse sublevel sets of horofunctions. 
We say that a subset $Hb\subset Y$ is a {\em horoball at $\eta\in\geo Y$}
if there exists a horofunction $h$ at $\eta$ such that 
$$ \{h\leq 0\} \subseteq Hb \subseteq \{h\leq 10C\de\} $$
with the constant $C$ from formula \eqref{eq:hfcdec}.
Horoballs are uniformly quasiconvex; in particular, for every $y\in Hb$, the ray $y\eta$ is contained in $N_r(Hb)$ for some uniform constant $r$.  a
Moreover, the visual boundary of a horoball at $\eta$ equals $\{\eta\}$. 
It is mostly these two properties of horoballs which will be used in this paper. 

We will call the space $Y$ itself a {\em horoball} if $\geo Y$ consists of a single point 
and the horofunctions are bounded above.

\medskip 
{\em Quasiconvex subsets and hulls.}
We recall that the {\em quasiconvex hull} $\QCH(A)\subset Y$ of a subset $A\subset Y$ 
is the union of all geodesic segments with endpoints in $A$.
The subset $A$ is called {\em $r$-quasiconvex} if $\QCH(A)\subset N_r(A)$ 
and {\em quasiconvex} if this holds for some $r>0$.
Note that pairs of points in a quasiconvex subset can be connected by uniform quasigeodesics inside it.

As a consequence of the $\delta$-hyperbolicity of $Y$, quasiconvex hulls are $C\de$-quasiconvex subsets,
and $\geo\QCH(A)=\geo A$. 
Both properties follow from the fact that any geodesic segment with endpoints in $\QCH(A)$
is contained in the tubular $C'\de$-neighborhood of a geodesic segment with endpoints in $A$.
This in turn reduces to the case when $A$ is (at most) a quadruple and follows from the thinness of triangles. 

The {\em quasiconvex hull} $\QCH(B)\subset Y$ of a subset $B\subset \geo Y$ at infinity 
is defined accordingly as the union of all geodesic lines $l\subset Y$ asymptotic to (points in) $B$, $\geo l\subset B$. 
It is nonempty unless $|B|\leq 1$, 
and then again it is $C\de$-quasiconvex and $\geo \QCH(B)=B$,
which follows from the fact that any geodesic segment with endpoints in $\QCH(B)$
is contained in the tubular $C'\de$-neighborhood of a geodesic line asymptotic to $B$.
An analogous property holds for rays $y\eta$ with $y\in\QCH(B)$ and $\eta\in B$.

\subsection{Isometries}

For (proper geodesic) Gromov hyperbolic spaces 
there is a {\em rough classification of isometries} into three types (elliptic, parabolic and hyperbolic) 
as in the case of CAT(0) spaces.

For an isometry $\phi$ of a Gromov hyperbolic space $Y$ consider the orbit maps $\Z\to Y,n\mapsto\phi^ny$ 
of the cyclic group $\<\phi\>$ generated by $\phi$. The isometry $\phi$ is called 

{\em elliptic} if the orbits are bounded;

{\em hyperbolic} if the orbits are quasigeodesics;

{\em parabolic} if the orbits are unbounded and distorted.\footnote{i.e. the orbit map $n\mapsto \phi^n y$ is not a quasiisometric embedding $\Z\to Y$}

The {\em asymptotic displacement number} 
of an isometry $\phi$
is defined as 
$$ \tau_{\phi}:= \lim_{n\to\infty} \frac{1}{n} d(y,\phi^n y).$$ 
The limit exists (due to the subadditivity of $n\mapsto d(y,\phi^n y)$) and is independent of $y\in Y$.
Note that $\tau_{\phi}>0$ if $\phi$ is hyperbolic and $\tau_{\phi}=0$ otherwise.

Non-elliptic isometries have unbounded orbits.
In particular,
they have infinite order and no fixed points in $Y$.
They do have fixed points at infinity:

\begin{prop}
(i) If $\phi$ is hyperbolic,
then it has exactly two fixed points on $\geo Y$, an attractive fixed point $\eta_+$ and a repulsive fixed point $\eta_-$.
It holds that for $y\in Y$, $\phi^ny\to\eta_{\pm}$ as $n\to\pm\infty$. 

(ii) If $\phi$ is parabolic,
then it has exactly one fixed point $\eta$ on $\geo Y$
and $\phi^ny\to\eta$ as $n\to\pm\infty$. 
\end{prop}
\proof
Suppose that $\phi$ is not elliptic. 
Then the sublevel subsets $\{\de_{\phi}\leq c\}$ 
of the displacement function $\de_{\phi}(y)=d(y,\phi y)$ are unbounded if non-empty,
because they are $\phi$-invariant.
Their visual boundary is therefore non-empty, $\geo \{\de_{\phi}\leq c\}\neq\emptyset$. 
On the other hand, it is fixed point-wise by $\phi$  
and therefore can contain at most two ideal points
(because $\phi$ is non-elliptic). 
The two point case corresponds to $\phi$ being hyperbolic.
Hence, if $\phi$ is parabolic,
then it has a unique fixed point $\eta$ in $\geo Y$.
Furthermore, 
the orbits of $\phi$ are contained in sublevel sets of $\de_{\phi}$
and therefore must accumulate at $\eta$.
\qed

\medskip
As a consequence,
parabolic and hyperbolic isometries can be characterized in terms of their action at infinity
and in terms of the accumulation of their orbits at infinity.

\medskip
We turn our attention to the stabilizers of points at infinity.
The non-hyperbolic isometries in the stabilizers can shift horofunctions only by a bounded amount:
\begin{lem}
Let $\phi$ be a non-hyperbolic isometry fixing $\eta\in\geo Y$ 
and let $h$ be a horofunction at $\eta$.
Then $|h-h\circ\phi|\leq C\de$.
\end{lem}
\proof
Fix $\eps>0$.
Suppose there exists a point $y\in Y$ such that 
$|h(y)-h\circ\phi(y)|\geq(1+\eps)C\de$.
We may assume that $h(y)-h\circ\phi(y)\geq(1+\eps)C\de$.
(Otherwise, we replace $\phi$ with $\phi^{-1}$.)
In view of \eqref{eq:hfceq}, it follows that $h-h\circ\phi \geq \eps C\de$ on all of $Y$.
This implies that the orbit path $n\mapsto\phi^ny$ is a quasigeodesic
and hence $\phi$ must be hyperbolic,
a contradiction. 
Letting $\eps\searrow0$ the assertion follows. 
\qed

\medskip
As a consequence,
if $P<\Isom(Y)$ is a subgroup fixing $\eta\in \geo Y$ which contains no hyperbolic isometry,
then $P$ quasi-preserves the horoballs at $\eta$.
In fact,
every horoball at $\eta$ is uniformly Hausdorff close to a $P$-invariant one. 

A horoball cannot be preserved by a hyperbolic isometry since it contains no quasigeodesic.

\begin{rem}
By \cite[Thm. 2G]{Tukia1994}, 
a discrete group of isometries of $Y$ fixing a point in $\geo Y$ cannot contain both hyperbolic and parabolic isometries. 
Hence, the stabilizer of an ideal point then consists either only of non-hyperbolic isometries 
or only of non-parabolic isometries.
\end{rem}

\subsection{Relatively hyperbolic groups} 
\subsubsection{Gromov's definition}
\label{sec:grdfn}

Since the geometrically finite subgroups and their generalizations considered in this paper 
will be {\em relatively hyperbolic} as abstract groups, we need to review the notion of relative hyperbolicity. 
There are various ways of defining relatively hyperbolic groups 
(see \cite{Osin, Drutu-Sapir, Hruska, Bowditch2012, Farb, Gerasimov-Potyagailo, GerasimovPotyagailo2, Yaman}).
We will work essentially
with Gromov's original definition \cite[\S 8.6]{Gromov}
in terms of actions on hyperbolic spaces.
Motivating examples 
are non-uniform lattices acting on rank one symmetric spaces and, more generally, geometrically finite Kleinian groups. 

\begin{definition}\label{defn:RH}
A {\em relatively hyperbolic group} is a pair $(\Ga, {\mathcal P})$
consisting of a group $\Ga$ and a conjugation invariant collection ${\mathcal P}$ of subgroups $\Pi_i< \Ga$, $i\in I$,
such that there exists a properly discontinuous isometric action $\Ga\acts Y$ 
on a $\de$-hyperbolic proper geodesic space $Y$ 
satisfying the following: 

(i) $Y$ is either taut or a horoball.

(ii) $Y$ is equipped with a $\Ga$-invariant collection $\B=(B_i)_{i\in I}$ of disjoint open horoballs 
such that the stabilizer of each $B_i$ in $\Ga$ is $\Pi_i$.

(iii) The action $\Ga\acts Y^{th}:= Y - \bigcup_{i\in I} B_i$ is cocompact.

(iv) The subgroups $\Pi_i$ are infinite. 

(v) The subgroups $\Pi_i$ are finitely generated.
\end{definition}
The subgroups $\Pi_i$ in this definition are called the {\em peripheral subgroups} of $\Ga$
and their collection ${\mathcal P}$ the {\em peripheral structure} on $\Ga$; 
the pair $(Y,{\mathcal B})$ or simply the hyperbolic space $Y$, respectively, the $\Ga$-action on it
is called a {\em Gromov model} for $(\Ga,{\mathcal P})$; 
the horoballs $B_i$ are called the {\em peripheral horoballs},
the truncated hyperbolic space $Y^{th}$ is called the {\em thick part} of $Y$,
and the horospheres $\Si_i:=\D B_i$ are called the {\em boundary} or {\em peripheral} horospheres of $Y^{th}$.

\begin{figure}[tbh]
\centering
\includegraphics[width=55mm]{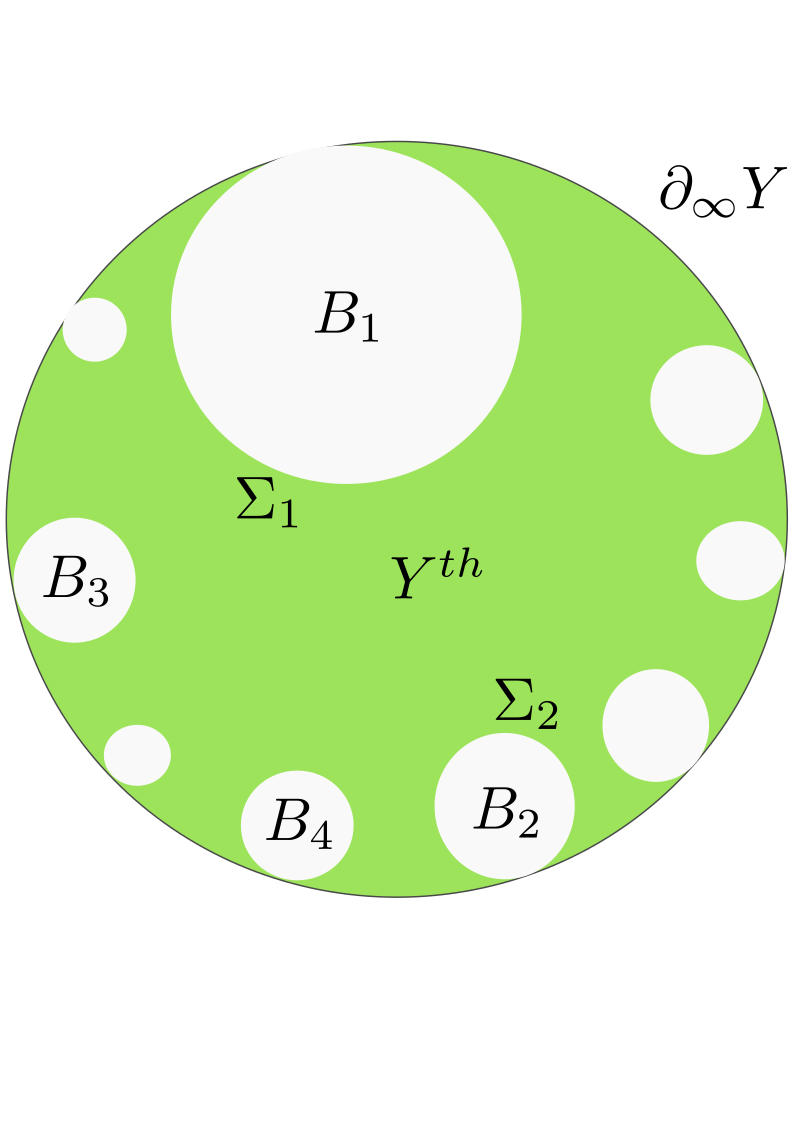}
\caption{A Gromov model.} 
\label{figure1.fig}
\end{figure}

For instance, 
in the case of nonuniform lattices acting on rank one symmetric spaces,
the natural Gromov model is the symmetric space itself.
For geometrically finite Kleinian groups, it is the closed convex hull of the limit set.\footnote{With two elementary exceptions: 
For finite groups the Gromov model is a singleton, while for Kleinian groups whose limit set is a single point $\zeta$ at infinity, 
the Gromov model is the intersection of a horoball with a certain convex subset, see the proof of Theorem \ref{thm:RM=GF}.} 

We call a Gromov model  {\em faithful} if there exists a point $y\in Y^{th}$ 
with trivial stabilizer in $\Ga$,
i.e. such that the action $\Ga\acts \Ga y$ is faithful. 
Every Gromov model can be modified to a faithful one by a slight enlargement,
e.g. by choosing a point $y\in Y^{th}$ and passing to the mapping cone $Y'$ of the orbit map 
$o_y: \Ga\to Y,\ga\mapsto\ga y$.
The geodesic segments $[\ga, \ga(y)]$ in $Y'$ 
are equipped with metrics as intervals of some fixed length $\la>0$; 
combined with the original path-metric on $Y$, this defines a path-metric on $Y'$ such that the natural inclusion map $Y\to Y'$ 
is a quasiisometry.  
The system of horoballs is kept the same.

\begin{rem}
\label{rem:rhdf}
1. When the peripheral structure is trivial, ${\mathcal P}=\emptyset$, then $\Ga$ is word hyperbolic.

2. $Y$ is a horoball if and only if $|{\mathcal P}|=1$.
Then the unique peripheral subgroup is $\Ga$ itself.

3. If $|\geo Y|=2$, then $\Ga$ is virtually cyclic and ${\mathcal P}=\emptyset$.
Indeed, by tautness $Y$ is then quasiisometric to a line and, in view of the infiniteness of peripheral subgroups,
$\Ga$ must be infinite, hence virtually cyclic,
and no infinite subgroup preserves a horoball.

4. By passing to subhoroballs, the peripheral horoballs can be made {\em $r$-separated}
for arbitrary $r>0$,
i.e. we may assume that any two distinct peripheral horoballs have distance $\geq r$.

5. In some treatments of the theory of relatively hyperbolic groups,
the peripheral subgroups $\Pi_i$ are not required to be infinite. 
However, if one omits this condition, 
then the peripheral horospheres $\D B_i$ with finite stabilizers $\Pi_i$ are compact,
the horoballs $B_i$ bounded by them are ends of $Y$ Hausdorff close to rays
and their centers are isolated points of $\geo Y$ which do not belong to the limit set of $\Ga$, 
compare Lemma~\ref{lem:prphsbgcchs} below whose proof uses only properties (i)-(iii) from our definition of relatively hyperbolic groups.
We do not want to allow this possibility. 

6. Gromov's original definition did not require the finite generation condition (v), only conditions (i-iv),
while other definitions discussed in the literature do require finite generation.
We added condition (v)
because under this assumption 
all known definitions of relatively hyperbolicity are equivalent (see \cite{Bowditch2012, Bumagin, Hruska} for proofs of the equivalences).\footnote{Specifically, Propositions 6.12 and 6.13 in Bowditch's paper \cite{Bowditch2012}
prove that Definition \ref{defn:RH} is equivalent to Definition 1 (and, hence, Definition 2) in   \cite{Bowditch2012}.}
Furthermore,  finite generation of the peripheral subgroups is a natural assumption in view of Theorem~\ref{thm:auslander}. 
\end{rem}

Several 
{\em finiteness} 
properties 
can be readily derived from the relatively hyperbolicity axioms:

The family ${\mathcal B}$ of peripheral horoballs is {\em locally finite},
i.e. every compact subset of $Y$ is intersected by only finitely many horoballs $B_i$.
This follows from the local compactness of $Y$ and since we can assume the peripheral horoballs to be $r$-separated for some $r>0$.

The cocompactness of the action $\Ga\acts Y^{th}$ further implies 
that there are {\em finitely many} $\Ga$-orbits of peripheral horoballs $B_i$, respectively, 
conjugacy classes of peripheral subgroups $\Pi_i$.

\begin{lem}
\label{lem:prphsbgcchs}
The actions $\Pi_i\acts\D B_i$ are cocompact.
\end{lem}
\proof 
Fix a base point $y\in Y$.
Since the action $\Ga\acts Y^{th}$ is cocompact,
there exists a subset $S\subset\Ga$ such that the subset 
$Sy\subset\Ga y$ 
is Hausdorff close to the horosphere $\D B_i$.
Then the horoballs $\ga^{-1}B_i$ for $\ga\in S$ intersect a compact subset and,
by the local finiteness of ${\mathcal B}$, 
only finitely many of them are different. 
It follows that $S$ is contained in a finite union of right cosets of $\Pi_i$,
and hence that $\Pi_iy$ and $\D B_i$ have finite Hausdorff distance. 
\qed

\medskip
Together with the finite generation of the peripheral subgroups this further implies:
\begin{lem}\label{lem:RH->fg} 
$\Ga$ is finitely generated.
\end{lem}
\proof
The previous lemma,  together with the finite generation of the peripheral subgroups $\Pi_i$, implies 
that the peripheral horospheres $\D B_i$ are coarsely connected (see \cite{Drutu-Kapovich}). 
Since $Y$ is path connected, it follows that the thick part $Y^{th}$ is coarsely connected.
The cocompactness of the action $\Ga\acts Y^{th}$ now implies that $\Ga$ is finitely generated.
\qed

\begin{cor}
Suppose that $\Ga$ satisfies parts (i)---(iii) of Definition \ref{defn:RH}. Then 
$\Ga$ is finitely generated if and only if all peripheral subgroups $\Pi_i, i\in I,$ are finitely generated. 
\end{cor}
\proof One direction is proven in Lemma \ref{lem:RH->fg}. The converse direction is proven in Part (d) of the main theorem 
of \cite{Gerasimov} (see also \cite{Drutu-Sapir} for the same implication with a different definition of relative hyperbolicity). \qed

\medskip
We describe next the {\em dynamics at infinity}.
The action $\Ga\acts\geo Y$ 
is a {\em convergence action} with certain characteristic features:

Since the horoballs $B_i$ are disjoint and form a $\Ga$-invariant family,
the stabilizers in $\Ga$ of the centers $\eta_i\in\geo Y$ of the horoballs $B_i$ are precisely the peripheral subgroups $\Pi_i$.
We can regard ${\mathcal P}$ as a subset of $\geo Y$
via the natural embedding $\Pi_i\mapsto\eta_i$.
The points $\eta_i$ are limit points of $\Ga$ 
due to our condition (iv) that the $\Pi_i$ are infinite,
however they are not conical.
They are called the {\em parabolic points} of $\Ga$ in $\geo Y$
and their stabilizers $\Pi_i$ the {\em maximal parabolic subgroups} of $\Ga$. 

All other points in $\geo Y$ are {\em conical limit points} of $\Ga$ 
as a consequence of the cocompactness of the action $\Ga\acts Y^{th}$.
(Recall that an ideal point $\eta\in\geo Y$ is a {\em conical limit point} in the dynamical sense of Definition~\ref{dfn:conlimbddpar}
if and only if the following geometric property is satisfied:
A(ny) geodesic ray in $Y$ asymptotic to $\eta$ has a tubular neighborhood 
which contains infinitely many points of a $\Ga$-orbit.)
In particular,
the limit set of $\Ga\acts Y$ is the entire $\geo Y$.

Note that the peripheral structure ${\mathcal P}$ can be read off the action $\Ga\acts\geo Y$
as the set of non-conical limit points and their stabilizers.

If $|\geo Y|\geq3$, 
then $\geo Y$ is a perfect metrizable compact topological space 
and the action $\Ga\acts\geo Y$ is a {\em minimal} convergence action.

The cocompactness of the actions $\Pi_i\acts\D B_i$ implies 
that also the actions $\Pi_i\acts\geo Y-\{\eta_i\}$ are properly discontinuous and cocompact,
i.e. the $\eta_i$ are {\em bounded parabolic points}, cf. Definition~\ref{dfn:conlimbddpar}.
Thus, the convergence action $\Ga\acts\geo Y$ is {\em geometrically finite} in the dynamical sense of Beardon-Maskit:

\begin{prop}
All points in $\geo Y$ are either conical limit points or bounded parabolic points
for the action of $\Ga$.
\end{prop}

In particular,
$\Ga$ is relatively hyperbolic in the sense of Bowditch's first definition in \cite{Bowditch2012}
which is formulated in terms of the dynamics at infinity of a properly discontinuous isometric action on a Gromov hyperbolic space.

\medskip
The Gromov model $Y$ is {\em not canonical} in the sense 
that its quasiisometry type is not determined by the pair $(\Ga, {\mathcal P})$.
This is because the action $\Ga\acts Y$ is cocompact only on the thick part $Y^{th}$
and the geometry of the peripheral horoballs is not controlled by the group.
Nevertheless, 
the {\em asymptotic} geometry of the Gromov models {\em is} determined.
By a remarkable result due to Bowditch \cite[Thm.~9.4]{Bowditch2012},
for any two Gromov models $\Ga\acts Y_1$ and $\Ga\acts Y_2$ 
there is a 
$\Ga$-equivariant homeomorphism $\geo Y_1\stackrel{\cong}{\to}\geo Y_2$. If $|\geo Y_1|\ge 3$, 
since the $\Ga$-actions on $\geo Y_1, \geo Y_2$ are minimal, this homeomorphism 
is necessarily unique. 

It follows that if the Gromov model $Y_1$ is faithful 
and the point $y_1\in Y_1$ has trivial stabilizer in $\Ga$,
then the $\Ga$-equivariant map of orbits $\Ga y_1\to\Ga y_2,\ga y_1\mapsto\ga y_2$ 
extends, by a homeomorphism at infinity, to an equivariant continuous map
\begin{equation}
\label{eq:orbcompgrmd}
\Ga y_1\sqcup\geo Y_1\stackrel{\cong}{\lra} \Ga y_2\sqcup\geo Y_2
\end{equation}
of the orbit closures inside the visual compactifications of the Gromov models.
Indeed, it can be read off the convergence dynamics of the action $\Ga\acts\geo Y_i$
whether a sequence $(\ga_ny_i)$ in the orbit $\Ga y_i$ converges in $\ol Y_i$ and, if yes, to which ideal point in $\geo Y_i$.
If also the point $y_2$ has trivial stabilizer in $\Ga$,
then \eqref{eq:orbcompgrmd} is a homeomorphism. 

This enables one to define a boundary at infinity and a compactification of a relatively hyperbolic group:

\begin{definition}[Ideal boundary]
The {\em ideal boundary} 
$\geo\Ga$ of a relatively hyperbolic group $(\Ga, {\mathcal P})$
is defined as the visual boundary $\geo Y$ of a Gromov model $Y$.
The {\em compactification} 
$$\ol\Ga=\Ga\sqcup\geo\Ga$$ of $\Ga$ is topologized at infinity by embedding it into the 
visual compactification $\ol Y=Y\sqcup\geo Y$ 
using an injective orbit map $\Ga\to Y$,
after enlarging $Y$ to a faithful Gromov model.
\end{definition}
Both $\geo \Ga$ and $\ol{\Ga}$ do not depend on the choice of the Gromov model and the orbit inside it.
To simplify notation, we will suppress the peripheral structure. 

\begin{rem}
Bowditch also constructed in \cite{Bowditch2012} 
a ``canonical'' Gromov model, unique up to (equivariant) quasiisometry,
with uniform {\em strict} exponential distortion of the peripheral horospheres. 
\end{rem}

The natural action $\Ga\acts\geo \Ga$ at infinity for a relatively hyperbolic group $\Ga$ 
is a minimal convergence action with finite kernel (unless $1\le |\geo \Ga|\le 2$) 
satisfying the Beardon-Maskit condition 
that every point $\eta\in \geo \Ga$ is either a conical limit point or a bounded parabolic fixed point. 
The stabilizers of the latter ones are the peripheral subgroups of $\Ga$.

Yaman showed that, conversely, 
the existence of an action with 
this kind of dynamics characterizes relatively hyperbolic groups,
thereby generalizing Bowditch's dynamical characterization of (absolutely) hyperbolic groups:

\begin{thm}
[Yaman \cite{Yaman}]
\label{thm:Yaman} 
Let $\Ga\acts Z$ be a convergence action on 
a nonempty perfect metrizable compact topological space $Z$.
Suppose that 

(i) each point in $Z$ is either a conical limit point or a bounded parabolic fixed point;

(ii) there are finitely many $\Ga$-orbits of bounded parabolic fixed points
and their $\Ga$-stabilizers 
are finitely generated. 

Then the family ${\mathcal P}$ of these stabilizers forms a relatively hyperbolic structure on $\Ga$,
and $Z$ is equivariantly homeomorphic to $\geo(\Ga,{\mathcal P})$. 
\end{thm}

\begin{rem}
1. As Yaman points out herself \cite[pp. 41-42]{Yaman},
the assumption that there are finitely many $\Ga$-orbits of bounded parabolic fixed points
can be dropped
as a consequence of a result by Tukia \cite[Thm 1B]{Tukia_conical}.

{2. The finite generation of the $\Ga$-stabilizers of bounded parabolic fixed points is needed in Yaman's paper \cite{Yaman} only indirectly, namely, because she verifies Definition 2 in \cite{Bowditch2012} and the latter requires finite generation of peripheral subgroups. }
\end{rem}

The following result on peripheral subgroups is also relevant for our paper.
Here, a space is said to have {\em coarsely bounded geometry} 
if there exists a scale $R_0>0$ and a function $\psi: [R_0,\infty)\to\N$ 
such that for all $R\geq R_0$ every $R$-ball in the space can be covered by at most $\psi(R)$ $R_0$-balls.

\begin{theorem}
[Dahmani, Yaman \cite{DY}] 
\label{thm:DY}
If a  relatively hyperbolic group 
admits a Gromov model 
with coarsely  bounded geometry, 
then all peripheral subgroups 
have polynomial growth.
\end{theorem}
Consequently,
by Gromov's theorem,
the peripheral subgroups then are virtually nilpotent. 

\begin{rem}
Dahmani and Yaman work with a stricter notion of bounded geometry:
They put $R_0=1$ and also require on the small scale that 
every $1$-ball can be covered by at most $\psi(R)$ balls of the radius $\frac{1}{R}$.
However, their proof only uses the assumption of coarsely bounded geometry.
\end{rem}

We observe that the property of coarsely bounded geometry behaves well under quasiisometric embeddings:
A space has coarsely bounded geometry as soon as it quasiisometrically embeds into a space with this property,
for instance, into a symmetric space. 
Therefore:
\begin{cor}\label{cor:DY}
If a  relatively hyperbolic group admits a Gromov model 
which quasiisometrically embeds into a symmetric space, 
then all peripheral subgroups are virtually nilpotent. 
\end{cor}

\subsubsection{Straight triples}
\label{sec:Gromov-straightness}

We carry over the notion of straightness of triples (defined in section~\ref{sec:HypSpaces}) 
from Gromov hyperbolic spaces to relatively hyperbolic groups as follows. 

For a relatively hyperbolic group $(\Ga,{\mathcal P})$
we consider the spaces of pairs
\begin{equation*}
(\Ga\sqcup\geo\Ga)^2 - \De_{\geo\Ga} \subset (\Ga\sqcup\geo\Ga)^2 
\end{equation*}
and triples
\begin{equation*}
T(\Ga,\geo\Ga) :=  (\Ga\sqcup\geo\Ga)\times\Ga\times(\Ga\sqcup\geo\Ga)
\end{equation*}
in $\ol\Ga$, 
compare \eqref{eq:idtrpy}.
If $(Y, {\mathcal B})$ is a 
Gromov model for $(\Ga,{\mathcal P})$ 
and $y\in Y$ is a point,
then the orbit map $o_y$ induces natural continuous maps 
$\ol\Ga^2\to\ol Y^2$ and $\ol\Ga^3\to\ol Y^3$ of pairs and triples
which restrict to maps 
$(\Ga\sqcup\geo\Ga)^2 - \De_{\geo\Ga} \to (Y\sqcup\geo Y)^2 - \De_{\geo Y}$ and $T(\Ga, \geo\Ga)\to T(Y, \geo Y)$.
Whether a family of triples in $\ol\Ga$ projects to a uniformly straight (in the sense of Definition~\ref{def:strtrp})
family of triples in $\ol Y$,
depends only on its asymptotics at infinity
and is therefore independent of the Gromov model $Y$ and the point $y$:

\begin{lemma}[Straightness is independent of Gromov model] 
\label{lem:instr}
The image in 
$T(Y, \geo Y)$
of a subset 
$T\subset T(\Ga, \geo\Ga)$ 
consists of $D$-straight triples for some uniform $D\geq0$
if and only if the subset of pairs 
$\{(\ga_2^{-1}\ol\ga_1,\ga_2^{-1}\ol\ga_3):(\ol\ga_1,\ga_2,\ol\ga_3)\in T\}$
is contained in $(\Ga\sqcup\geo\Ga)^2 - \De_{\geo\Ga}$
as a relatively compact subset.
\end{lemma}
\proof 
By $\Ga$-equivariance, 
the triples in the image of 
$T$ in $T(Y, \geo Y)$
are $D$-straight if and only if the triples in the image of the family 
$E:=\{(\ga_2^{-1}\ol\ga_1,e,\ga_2^{-1}\ol\ga_3): (\ol\ga_1,\ga_2,\ol\ga_3)\in T\}$ are. 
The latter holds for some uniform $D\geq0$ 
if and only if 
the family of all geodesics in $Y$ with pair of endpoints 
in the image $E_Y$ of $E$  
under the natural map 
$(\Ga\sqcup\geo\Ga)^2 - \De_{\geo\Ga} \to (Y\sqcup\geo Y)^2 - \De_{\geo Y}$
is bounded 
(in the sense defined in section~\ref{sec:HypSpaces}). 
This condition (by Lemma~\ref{lem:bddasyprp}) holds  
if and only if $E_Y$ is relatively compact in 
$(Y\sqcup\geo Y)^2 - \De_{\geo Y}$,
so in view of the continuity of the map $\ol\Ga^2\to\ol Y^2$,
if and only if $E$ is relatively compact in 
$(\Ga\sqcup\geo\Ga)^2 - \De_{\geo\Ga}$.
\qed

\medskip
It therefore makes sense to call a family of triples in 
$T(\Ga, \geo\Ga)$ 
{\em straight} if its image in 
$T(Y, \geo Y)$
consists of $D$-straight triples for some fixed data $(D,Y,y)$.

Note that straightness is a useful concept for triples in relatively hyperbolic groups $\Ga$ only if $|\geo \Ga|>1$.
If $\geo\Ga$ is a singleton, then 
a family of triples $(\ga_1, \ga_2, \ga_3)$ in $\Ga^3$ is straight if and only if
the corresponding subsets $\{\ga_2^{-1}\ga_1,\ga_2^{-1}\ga_3\}$ intersect some finite subset of $\Ga$.

\section{Some geometry of higher rank symmetric spaces}

\subsection{Basic notions and standing notation}

In this section we briefly discuss some 
basic definitions pertaining to symmetric spaces $X$ of noncompact type.
(We will call them simply {\em symmetric spaces}.)
We refer the reader to the books \cite{Helgason,Eberlein} and to \cite{habil} for the foundational material,
and to our earlier papers \cite{morse,mlem,coco15,anolec,bordif,manicures} for more specialized aspects of the theory,
developed specifically to study the asymptotic geometry and discrete isometry groups of symmetric spaces.  

The visual boundary $\geo X$ of a symmetric space $X$ admits a structure as a 
thick spherical building (the Tits building of $X$). Throughout the paper we will use the notation $\simod$ for the model spherical chamber of this building, 
$\De$ for the model euclidean Weyl chamber of $X$ and $\theta: \geo X\to \simod$ for the {\em type projection}.
The full isometry group of $X$ acts on $\simod$ isometrically; 
the map $\theta$ is equivariant with respect to this action. 
We will denote by $G< \Isom(X)$ the kernel of this action,
i.e. the subgroup of type preserving isometries.
It is a semisimple Lie group and has finite index in $\Isom(X)$.

We will freely use the notions introduced in our earlier papers,
such as 
the opposition involution $\iota$ of $\simod$, 
a type $\bar\theta\in\simod$, 
the face types $\taumod\subseteq\simod$ \cite[\S 2.2.2]{anolec}, 
the associated $\taumod$-flag manifolds $\Flagt$ \cite[\S 2.2.2, 2.2.3]{anolec},
the open Schubert cells $C(\tau)\subset \Flagt$ \cite[\S 2.4]{anolec}, 
the $\taumod$-boundary $\Dt\De$ of $\De$ \cite[\S 2.5.2]{anolec}, 
the $\Delta$-valued distance $d_\Delta$ on $X$ \cite[\S 2.6]{anolec}, 
$\Theta$-regular geodesic segments \cite[\S 2.5.3]{anolec}, 
parallel sets $P(\tau_-,\tau_+)$, stars $\st(\tau)$, open stars $\ost(\tau)$, $\Theta$-stars $\st_\Theta(\tau)$, 
Weyl cones $V(x, \st(\tau))$ and $\Theta$-cones $V(x, \st_\Theta(\tau))$, 
diamonds $\diamot(x,y)$ and $\Theta$-diamonds $\diamoTh(x,y)$ \cite[\S 2.5]{anolec}, 
$\taumod$-regular sequences and subgroups \cite[\S 4.2]{anolec}), 
uniformly $\taumod$-regular sequences and subgroups \cite[\S 4.6]{anolec}, 
$\taumod$-convergence subgroups, flag-convergence\footnote{ {Here we note (for readers familiar with Satake compactifications) that $\taumod$-convergence of sequences in $X$ to points in $\Flagt$ is equivalent to the convergence in the suitable Satake compactification of $X$, where the flag-manifold $\Flagt$ is the smallest (closed) stratum.}},  the Finsler interpretation of 
flag-convergence, see \cite[\S 4.5 and 5.2]{bordif} and \cite{anolec},
$\taumod$-limit sets $\LaXt=\Lat=\Lat(\Ga)\subset \Flagt$ \cite[\S 4.5]{anolec},
visual limit sets \cite[p.\ 4]{anolec}, 
Morse subgroups \cite[\S 5.4]{anolec}, Morse quasigeodesics and Morse maps \cite[Defs. 5.31, 5.33]{mlem}, 
antipodal limit sets \cite[Def. 5.1]{anolec} and antipodal maps to flag manifolds \cite[Def. 6.11]{mlem},
to name a few.
We review some of this material in sections \ref{sec:finsler} and \ref{sec:regularity}. We also refer to Appendix 1 for more a 
detailed discussion.

We will use the following conventions and standing notation.

Throughout the paper, $X$ will denote a symmetric space of noncompact type.
We will denote by $\Theta$ an $\iota$-invariant, compact, Weyl-convex (see \cite[Def. 2.7]{anolec}) 
subset of the open star $\ost(\taumod)\subset\simod$. 
For pairs $\Theta,\Theta'\subset\ost(\taumod)$ of such subsets, we will always assume that $\Theta\subset\inte(\Theta')$. 
Similarly, for pairs of positive constants $d, d'$ we will always assume that $d< d'$. 
Note that, when $X$ has rank one, the data $\taumod, \Theta$ are obsolete. 
In this case, we also have that $\D\simod=\emptyset$ and $\Theta=\inte(\simod)=\simod$ is clopen; in particular, $\Theta'=\inte(\Theta)=\Theta$.

\subsection{Finsler geometric notions}\label{sec:finsler}

In \cite{bordif}, see also \cite{anolec},
we considered a certain class of $G$-invariant ``polyhedral'' Finsler metrics on $X$.
Their geometric and asymptotic properties turned out to be well adapted 
to the study of geometric and dynamical properties of regular subgroups.
They provide a Finsler geodesic {\em combing} of $X$ which is, in many ways, more suitable for analyzing the asymptotic geometry of $X$ than the geodesic combing given by the standard Riemannian metric on $X$.
These Finsler metrics also play a basic role in the present paper.
We briefly recall their definition and some basic properties, 
and refer to \cite[\S 5.1]{bordif} for more details.

Let $\bar\theta\in\inte(\taumod)$ be a type spanning the face type $\taumod$. {Recall that $b_\xi$ denotes the {\em Busemann function} on $X$ associated with $\xi$.} 
The {\em $\bar\theta$-Finsler distance} $d^{\bar\theta}$ on $X$ is the $G$-invariant pseudo-metric defined by
\begin{equation*}
d^{\bar\theta}(x,y) := \max_{\theta(\xi)=\bar\theta} \bigl( b_{\xi}(x)-b_{\xi}(y) \bigr) 
\end{equation*}
for $x,y\in X$, 
where the maximum is taken over all ideal points $\xi\in\geo X$ with type $\theta(\xi)=\bar\theta$.
{Equivalently, $d^{\bar\theta}(x,y)$ can be defined as the inner product of $d_\Delta(x,y)$ and the vector $\bar\theta$.} 
The  $\bar\theta$-Finsler distance is positive, i.e. a (non-symmetric) metric, 
if and only if the radius of $\simod$ with respect to $\bar\theta$ is $<\pihalf$.
This is in turn equivalent to $\bar\theta$ not being contained in a factor of a nontrivial spherical join decomposition of $\simod$,
and is always satisfied e.g. if $X$ is irreducible.

If $d^{\bar\theta}$ is positive,
it is equivalent to the Riemannian metric. 
In general, if it is only a pseudo-metric,
it is still equivalent to the Riemannian metric $d$ on uniformly regular pairs of points.
More precisely,
if the pair of points $x,y$ is $\Theta$-regular, then 
$$
L^{-1} d(x,y) \le d^{\bar\theta}(x,y)\le L d(x,y)
$$ 
with a constant $L=L(\Theta)\ge  1$.

Regarding symmetry of the Finsler distance, one has the identity
\begin{equation*}
d^{\iota\bar\theta}(y,x)  = d^{\bar\theta}(x,y)
\end{equation*}
and hence $d^{\bar\theta}$ is symmetric if and only if 
$\iota\bar\theta=\bar\theta$.
We refer to $d^{\bar\theta}$ as a Finsler metric {\em of type $\taumod$}.

The $d^{\bar\theta}$-balls in $X$ are convex but not strictly convex.
(Their intersections with flats through their centers are polyhedra.)
Accordingly,
$d^{\bar\theta}$-geodesics connecting two given points $x,y$ are not unique. 
To simplify notation,
$xy$ will stand for {\em some} $d^{\bar\theta}$-geodesic connecting $x$ and $y$.
The union of all $d^{\bar\theta}$-geodesic $xy$ equals the $\taumod$-diamond $\diamot(x,y)$,
that is, a point lies on a $d^{\bar\theta}$-geodesic $xy$ if and only if it is contained in $\diamot(x,y)$, see \cite{anolec}.
Finsler geometry thus provides an alternative description of diamonds.  
Note that with this description, 
the diamond $\diamot(x,y)$ is also defined when the segment $xy$ is not $\taumod$-regular. 
Such a {\em degenerate} $\taumod$-diamond is contained in a smaller totally-geodesic subspace,
namely in the intersection of all $\taumod$-parallel sets containing the points $x,y$.
The description of geodesics and diamonds 
also implies that the unparameterized $d^{\bar\theta}$-geodesics depend only on the face type $\taumod$,
and not on $\bar\theta$.
We will refer to $d^{\bar\theta}$-geodesics as {\em $\taumod$-Finsler geodesics}.
Note that Riemannian geodesics are Finsler geodesics.

We will call a $\Theta$-regular $\taumod$-Finsler geodesic a {\em $\Theta$-Finsler geodesic}.
If $xy$ is a $\Theta$-regular (Riemannian) segment,
then the union of $\Theta$-Finsler geodesics $xy$ equals the $\Theta$-diamond $\diamoTh(x,y)$. 

Every $\taumod$-Finsler ray in $X$ is contained in a $\taumod$-Weyl cone,
and we will use the notation $x\tau$ for a $\taumod$-Finsler ray contained {in} $V(x,\st(\tau))$.
Similarly, 
every $\taumod$-Finsler line is contained in a $\taumod$-parallel set,
and we denote by $\tau_-\tau_+$ an oriented $\taumod$-Finsler line 
forward/backward asymptotic to two antipodal simplices $\tau_{\pm}\in\Flagt$ 
and contained in $P(\tau_-,\tau_+)$.

\subsection{Types of isometries}
\label{sec:isomsym}

Let $g\in \Isom(X)$.
The function $\de_g(x)=d(x,gx)$ on $X$ is called the {\em displacement function} of $g$
and the number 
\begin{equation*}
m_g:= \inf_X\de_g
\end{equation*}
is called the {\em infimal displacement} or {\em translation number} of $g$.

\medskip 
The isometry $g$ is called {\em semisimple} if $\de_g$ attains its infimum.
The semisimple isometries split into two subclasses:
A semisimple isometry $g$ is called 

(i) {\em elliptic} if $m_g=0$, i.e. if it fixes a point in $X$. 
Equivalently, the orbits in $X$ of the cyclic group $\<g\>$ are bounded. 

(ii) {\em axial} or {\em hyperbolic} if $m_g > 0$. 
In this case, the minimum set $\Min(g)$ of $\de_g$ is the union of the {\em axes} of $g$,
i.e. of the $g$-invariant geodesic lines. 
On each axis, $g$ acts as a translation by $m_g$.
The subset $\Min(g)$ is a symmetric subspace of $X$ and splits metrically as $\Min(g)\cong\R\times CS$,
the lines $\R\times pt$ being the axes of $g$
and the cross section $CS$ being a symmetric (sub)space.

A hyperbolic
isometry $g$ is 
a {\em transvection} 
if it preserves the parallel vector fields along some (and hence any) axis.
Then every line parallel to an axis of $g$ is itself an axis,
i.e. the minimum set is the full parallel set of a line. 
The transvections are the isometries of $X$ which can be written as the product of two distinct point reflections.

\medskip
The isometries $g$ for which $\de_g$ does not attain its infimum
are called {\em parabolic}.
A parabolic isometry $g$ has at least one fixed point in the visual boundary.
To see this,
consider a sequence $(x_n)$ in $X$ such that $\de_g(x_n)\searrow m_g$.
Then $x_n\to\infty$ and the accumulation points of $(x_n)$ in $\geo X$ are fixed by $g$. 
Moreover, at some of the fixed points at infinity also the horoballs are preserved by $g$.
Namely, 
choose a sequence $(x_n)$ more carefully,
by picking a base point $o\in X$ and a sequence {$\eps_n\searrow m_g$ } 
and letting $x_n$ be the nearest point projection of $o$ to $\{\de_g\leq\eps_n\}$.
Then for any accumulation point $\xi\in\geo X$ of $(x_n)$
the horoballs centered at $\xi$ are $g$-invariant,
see e.g. \cite[{Prop. 3.4}]{Ballmann}.

The parabolic isometries 
break up into several subclasses.
We will call a parabolic isometry $g$ {\em strictly parabolic} if $m_g=0$ 
and non-strictly parabolic otherwise. 
If $\rank X=1$ then all parabolic isometries are strictly parabolic,
but non-strictly parabolic isometries occur if $\rank X\geq2$.

An isometry $g\ne \id_X$ is called {\em unipotent} 
if the closure of its conjugacy class in $\Isom(X)$ contains $\id_X$,
i.e. if there exists a sequence of isometries $h_n\to\infty$ such that $h_ngh_n^{-1}\to\id_X$.
In this case,
there exists a transvection $h$ such that 
$$
\lim_{n\to\infty} h^n g h^{-n}= \id_X. 
$$
Unipotent isometries are strictly parabolic.

\medskip
Every isometry $g$ of $X$ has a unique {\em Jordan decomposition}
\begin{equation}\label{eq:Jordan}
g=g_sg_u= g_t g_e  g_u
\end{equation}
where $g_s= g_t g_e$ and 
$g_t, g_e, g_u$ are commuting isometries which are, respectively, a transvection, elliptic and unipotent.
The factor $g_s$ is semisimple. 
Note that 
$$
m_g = m_{g_s} = m_{g_t} 
$$
and $\Min(g_s)$ is preserved by $g_u$.

If $g$ is non-strictly parabolic, equivalently, if $g_t$ and $g_u$ are nontrivial,
then $g_u$ preserves the cross sections $\{t\}\times CS$ of $\Min(g_t)\cong\R\times CS$
and acts on them as a strictly parabolic isometry.

{We refer the reader to \cite{Eberlein} for further discussion. }

\medskip 
If $u\in G$ is a unipotent isometry,
then it is of the form $u=\exp(n)$ with $n\in{\mathfrak g}$ nilpotent. (Here ${\mathfrak g}$ is the Lie algebra of $G$.) 
According to the Morozov-Jacobson theorem
regarding nilpotent elements in semisimple Lie algebras, see \cite{Jacobson},
$n$ belongs to a 3-dimensional simple Lie subalgebra
${\mathfrak g}'\cong sl(2, \R)$.
Correspondingly,
$u$ lies in a rank one Lie subgroup $G'< G$
locally isomorphic to $SL(2, \R)$.
The subgroup $G'$ preserves a totally-geodesic hyperbolic plane $X'\subset X$
and $u$ acts on it as a parabolic element.
Consequently, 
$u$ fixes a unique ideal point $\xi\in\geo X'$ and preserves the horocycles in $X'$ centered at $\xi$.
It follows that its orbits accumulate in $\ol X$ at $\xi$, 
$\La(\<u\>)=\{\xi\}$.

More generally,
if $g\in G$ is strictly parabolic,
then $g_s$ is elliptic, $g_t=\id_X$, 
and $g$ has the Jordan decomposition $g=g_eg_u$ with $g_u\neq \id_X$.
The latter implies that the $g_u$-invariant fixed subspace $\Fix(g_e)\subset X$ is noncompact.
As above, it contains a $g_u$-invariant hyperbolic plane $X'$
on which $g_u$ acts as a parabolic isometry.
The $g$-orbits in $X$ have bounded distance from the $g_u$-orbits,
and it follows that they accumulate in $\ol X$ at a unique ideal point $\xi\in\geo X'\subseteq\geo\Fix(g_e)$,
\begin{equation}
\label{eq:vlimststrpb}
\La(\<g\>)=\{\xi\}.
\end{equation}
The hyperbolic plane $X'$ and the horocycles in it centered at $\xi$ are $g$-invariant. 

We define the {\em type} of the strictly parabolic isometry $g$ as $\theta(g):=\theta(\xi)\in\simod$
and its {\em face type} $\taumod(g)\subseteq\simod$ as the face of $\simod$ spanned by its type.
Note that both are $\iota$-invariant because the points in the visual boundary of a rank one symmetric subspace of $X$ (such as $X'$)
are pairwise antipodal and therefore have the same $\iota$-invariant type.

\medskip
Let $g$ be an isometry which fixes an ideal point $\xi\in\geo X$.
Then $g$ induces an isometry $g_{\xi}$ on the space $X_{\xi}$ of strong asymptote classes at $\xi$ ({see Appendix 1}), 
cf. e.g. \cite{habil,bordif}.
If $g$ preserves also the horoballs at $\xi$, then 
\begin{equation}
\label{eq:trnsnind}
m_{g_{\xi}}=m_g
\end{equation}

For every isometry $g$ of $X$ it holds that 
\begin{equation}
\label{eq:trnsnlin}
m_{g^n}=n m_g
\end{equation}
for $n\in\N_0$.
This is clear for semisimple isometries.
If $g$ is parabolic,
it follows by induction on the rank of $X$ using \eqref{eq:trnsnind},
or from the Jordan decomposition. 

\medskip 
As in the case of Gromov hyperbolic spaces,
one can relate the rough classification of isometries to the {\em distortion} of their {\em orbit paths}.
For an isometry $g$ of $X$ consider the orbit paths $\Z\to X,n\mapsto g^nx$ of the cyclic group $\<g\>$ generated by it.
One has the following:

If $m_g=0$, equivalently, if $g$ is elliptic or strictly parabolic,
then $n\mapsto\de_{g^n}(x)$ grows sublinearily as $n\to\infty$
and the orbit paths are distorted (not quasiisometrically embedded).
On the other hand,
if $m_g>0$, equivalently, if $g$ is hyperbolic or non-strictly parabolic,
then \eqref{eq:trnsnlin} or the Jordan decomposition implies 
that the orbit paths 
are undistorted.
Thus the orbits of $g$ are undistorted if and only if $m_g>0$.

We obtain a more precise picture by applying the Jordan decomposition:

If $g$ is a strictly parabolic isometry, 
then, as we noted earlier, $g$ preserves a hyperbolic plane $X'\subset X$ and acts on $X'$ as a parabolic isometry. 
In particular,
the orbits of $\<g\>$ in $X$ are {\em logarithmically distorted},
$\de_{g^n}(x)=O(\log |n|)$. 

It follows for an arbitrary isometry $g$ 
that its orbit paths deviate sublinearily, in fact logarithmically, 
from the orbit paths of its semisimple part $g_s$,
\begin{equation}
\label{eq:sblndistorb}
d(g^nx,g_s^nx) = O(\log |n|) 
\end{equation}
Thus, if $m_g>0$ and $l$ is an oriented axis of $g_s$,
then 
\begin{equation*}
g^nx \to l(\pm\infty)\in\geo X
\end{equation*}
in the visual compactification as $n\to\pm\infty$.

We conclude:

\begin{prop}[Distortion of orbits of isometries]
\label{prop:dstorbism}
Let $g$ be an isometry of $X$.

If $g$ is {\em elliptic}, then its orbits are bounded.

If $g$ is {\em strictly parabolic}, then its orbits are unbounded, but logarithmically distorted.
They accumulate in $\ol X$ at a single ideal point in $\geo X$.
(It lies in the visual boundary of a $g_u$-invariant totally geodesic hyperbolic plane.)

If $g$ is 
{\em hyperbolic}, 
then its orbits are undistorted.
They are Hausdorff close to an(y) axis $l$ of $g$ and accumulate 
in $\ol X$ 
at the pair of antipodes $\geo l\subset\geo X$.

If $g$ is {\em non-strictly parabolic}, then its orbits are undistorted.
They deviate sublinearily, in fact, logarithmically, from an(y) axis $l$ of the semisimple part $g_s$ 
but they are not Hausdorff close to any line.
They accumulate in $\ol X$ 
at the pair of antipodes $\geo l\subset\geo X$.
\end{prop}

{In particular, the vanishing resp. positivity of $m_g$ can be read off {\em coarse} properties of the $\<g\>$-orbits: $m_g>0$ if and only if each $\<g\>$-orbit is undistorted in $X$.}

\subsection{Regularity}\label{sec:regularity}

\subsubsection{Notions of regularity and limit sets}
As in our earlier papers,
we will be imposing certain {\em regularity assumptions} on discrete subgroups $\Ga< G$. 
In this section, we go through some variations of the notion of regularity.

\begin{rem}
It is imperative to note here that 
notions of $\taumod$-regularity and $\taumod$-limit sets, 
as well as the relation to convergence-type dynamics
and many other indispensable related concepts and theorems 
have their origin in the foundational paper by Yves Benoist \cite{Benoist}, 
in fact, 
even earlier in the work of Tits \cite{Tits} and Guivarc'h \cite{Guivarch}. 
\end{rem}

A subset of $X$ is called {\em $\taumod$-regular} 
if all divergent sequences in it are $\taumod$-regular.
A map into $X$ is called {\em $\taumod$-regular} 
if its image is $\taumod$-regular.

The following  strengthening of  regularity occurs naturally in equivariant settings:

\begin{dfn}[Weakly uniformly regular]\label{defn:WUR}
We say that an (unbounded) subset $W\subset X$ is {\em $(\taumod,\phi)$-regular} 
if for $x, x'\in W$ 
$$d(d_{\De}(x,x'),\Dt\De)\geq \phi(d(x,x'))$$
where $\phi:[0,+\infty)\to\R$ is a monotonic function with $\lim_{d\to+\infty}\phi(d)=+\infty$.
We say that a subset $W\subset X$ is {\em weakly uniformly $\taumod$-regular} 
if it is $(\taumod,\phi)$-regular for some $\phi$.

Accordingly, we say that a map into $X$ is {\em $(\taumod,\phi)$-regular} or {\em weakly uniformly $\taumod$-regular} 
if its image in $X$ is.
\end{dfn}

The orbits $\Ga x\subset X$ of $\taumod$-regular actions $\Ga\acts X$ are weakly uniformly $\taumod$-regular subsets. 

Weak uniform regularity is stable under bounded perturbation:
If $W'\subset X$ is $d$-Hausdorff close to $W\subset X$ and $W$ is $(\taumod,\phi)$-regular,
then $W'$ is $(\taumod,\phi(\cdot-2d)-2d)$-regular,
as follows from the (weak) $\De$-triangle inequality
{
$$
||d_\Delta(x,y)- d_\Delta(y,z)||\le d(x,z),
$$
see \cite{KLM}. 
}

Note that 
a subset of $X$ is {\em uniformly $\taumod$-regular} (see \cite[\S 4.6]{anolec}) 
if and only if it is $(\taumod,\phi)$-regular for some (affine) linear function $\phi$. 

\medskip
For a $\taumod$-regular subset $W\subset X$, we define $\geot W\subset\Flagt$ as its {\em $\taumod$-accumulation set}.
Similarly, we define the $\taumod$-{\em conical accumulation set} $\geotc W\subset\Flagt$ 
as the set of conical $\taumod$-limits of sequences in $W$ (see \cite[Def. 5.33]{anolec}).

For a $\taumod$-regular subgroup $\Ga< G$, 
besides the limit set $\Lat=\LaXt=\geot(\Ga x)$ we will also consider the {\em conical $\taumod$-limit set} 
$$
\LaXtc:=\geotc(\Ga x)\subset \LaXt.$$

A $\taumod$-regular  subgroup $\Ga< G$ is said to be $\taumod$-antipodal if its limit set $\Lat$ is antipodal, 
i.e. if any two distinct points in $\Lat$ are antipodal.
A $\taumod$-regular subgroup $\Ga< G$ is called {\em $\taumod$-elementary} if $|\Lat|\le 2$.

It is a basic fact connecting the theory of regular discrete subgroups of $G$ to the classical theory of Kleinian groups, that  
each $\taumod$-regular antipodal subgroup $\Ga< G$ acts as a convergence group on its $\taumod$-limit set,  
see  \cite[\S 5.1]{anolec} or \cite[Corollary 3.16]{manicures}. 
In particular, for a nonelementary $\taumod$-regular antipodal subgroup $\Ga< G$, 
its $\taumod$-limit set $\Lat$  is perfect and every $\Ga$-orbit is dense in $\Lat$.

\begin{example}\label{ex:rank1embedding}
Let $G_1<G$ be a connected rank one simple Lie subgroup.
By the Karpelevich--Mostow theorem,
there exists a rank one symmetric subspace $X_1\subset X$
which is a $G_1$-orbit. 
Its visual boundary $\geo X_1\subset \geo X$ is an antipodal subset. 
Hence, it consists of ideal points of the same $\iota$-invariant type $\bar\xi\in\simod$,
$\theta(\geo X_1)=\{\bar\xi\}$.
We call $\theta(\geo X_1):=\bar\xi$ the {\em type} of the rank one subspace $X_1$,
and the $\iota$-invariant face $\taumod(X_1):=\taumod(\bar\xi)\subseteq\simod$ spanned by $\bar\xi$ its {\em face type}. 
All non-degenerate segments in $X_1$ have type $\bar\xi$.
We thus have a map 
$$
\geo X_1\to \Flag_{\taumod(X_1)}(X) 
$$
sending $\xi_1\in \geo X_1$ to the simplex $\tau_{\xi_1}\in \Flag_{\taumod(X_1)}$ spanned by 
$\xi_1$. 
Also, for every $\iota$-invariant face $\taumod\subseteq\taumod(X_1)$, by composing this map with the projection
$\Flag_{\taumod(X_1)}\to \Flagt$, one obtains a natural antipodal embedding $\beta: \geo X_1\to \Flagt$. 
All divergent sequences in $X_1$ 
are uniformly $\taumod$-regular and the $\taumod$-accumulation set
of $X_1$ in $\Flagt$ equals $\beta(\geo X_1)$. 
Every discrete subgroup $\Ga_1< G_1$ is uniformly $\taumod$-regular as a subgroup of $G$.
Moreover, $\Lat(\Ga_1)=\beta(\La(\Ga_1))$, where $\La(\Ga_1)\subset \geo X_1$ 
is the visual limit set of $\Ga_1$. 
\end{example}

\subsubsection{Zariski dense subgroups}

In general, verifying the (uniform) regularity of a subgroup is not an easy task. However, it is simpler 
for Zariski dense subgroups (see our paper \cite[Theorem 9.6]{bordif}):

\begin{thm}\label{thm:Zardense}
Let $\Ga<G$ be Zariski dense.
Suppose that $Z$ is a compact metrizable space, $\Ga\acts Z$ is a convergence action 
and $\beta: Z\to \Flagt$ is a $\Ga$-equivariant antipodal continuous map. 
Then $\Ga$ is $\taumod$-regular.  
\end{thm}

Moreover, we can say about the relation of the image to the limit set of $\Ga$:
\begin{add}
\label{add:Zarden}
If, in addition, the action $\Ga\acts Z$ is minimal,
then $\beta(Z)=\Lat$.
\end{add}
\proof 
Let $\la_+\in \Lat$.
Then, 
since regular subgroups act on flag manifolds as discrete convergence groups (see \cite[Lemma 4.19]{anolec}), 
there exist a sequence $(\ga_n)$ in $\Ga$ and a point $\la_-\in \Flagt$ 
such that $\ga_n|_{C(\la_-)}\to\la_+$ uniformly on compacts ($C(\la_-)\subset\Flagt$ being the open Schubert cell).
The complement $\Flagt - C(\la_-)$ is a proper subvariety of $\Flagt$.
Hence, 
by the Zariski density, $\beta(Z)\cap C(\la_-)\neq\emptyset$.
For any point $\tau$ in the intersection,
it holds that $\ga_n(\tau)\to \la_+$. 
Since $\beta(Z)$ is closed and $\Ga$-invariant, it follows that $\la_+\in\beta(Z)$. 
Thus $\Lat\subseteq\beta(Z)$. 
The minimality of the action $\Ga\acts Z$ implies equality.
\qed

\subsubsection{Accumulation sets of regular sequences}

We collect some facts needed later in the paper.

\begin{lem}
\label{lem:smtmlim}
Suppose that $(x_n)$ and $(y_n)$ are uniformly $\taumod$-regular sequences in $X$ 
such that 
$\frac{d(x_n,y_n)}{d(x_n,o)}\to0$,
where $o\in X$ is a base point.
Then their $\taumod$-accumulation sets in $\Flagt$ coincide. 
\end{lem}
\proof
It suffices to consider the case when $(x_n)$ $\taumod$-flag converges, $x_n\to\tau\in\Flagt$,
and to show that then also $y_n\to\tau$.

We extend the Riemannian segments $ox_n$ and $oy_n$ to Riemannian rays $o\xi_n$ and $o\eta_n$.
By uniform regularity,
we may assume that the types $\theta(\xi_n),\theta(\eta_n)$ of the ideal points $\xi_n,\eta_n\in\geo X$ 
are contained in a compact subset $\Theta\subset\ost(\taumod)$.
Let $\tau_{\xi_n},\tau_{\eta_n}\in\Flagt$ denote the simplices in $\geo X$ spanned by them.
Then $\tau_{\xi_n}\to\tau$ in $\Flagt$ and we must show that also $\tau_{\eta_n}\to\tau$.

After extraction, we may assume that also $(\xi_n)$ converges in $\geo X$, $\xi_n\to\xi$.
Then $\xi\in\tau$. 
By our assumption, $\angle_o(\xi_n,\eta_n)=\angle_o(x_n,y_n)\to0$.
It follows that also $\eta_n\to\xi$. 
In view of uniform regularity, $\theta(\eta_n)\in\Theta$, this implies that $\tau_{\eta_n}\to\tau$.
\qed

\medskip

Let $K<G$ denote a maximal compact subgroup.
We denote by $o\in X$ its fixed point. 
\begin{lem}
\label{lem:nbrgsqcs}
Let $x_n\to\infty$ be a uniformly $\taumod$-regular sequence in $X$ which flag converges, $x_n\to\tau\in\Flagt$.
Let $(k_n)$ be a sequence in $K$ such that $d(x_n,k_nx_n)$ is uniformly bounded. 

Then $(k_n)$ accumulates at $\Stab_K(\tau)<K$.
\end{lem}
\proof
Let $o\in X$ be the fixed point of $K$.
After passing to a subsequence,
we may assume that the Riemannian segments $ox_n$ converge to a $\Theta$-regular ray $o\xi$, $\xi\in\ost(\tau)$,
and that $k_n\to k$. 
Since $d(x_n,k_nx_n)$ is uniformly bounded, $k$ fixes $\rho$ and hence $\tau$.
(Compare Lemma~\ref{lem:smtmlim}.)
\qed

\medskip
The visual and flag accumulation sets of regular sequences are related as follows:

\begin{lem}
\label{lem:vsaccsrgsq}
Let  $(x_n)$ be a 
$\taumod$-regular sequence in $X$ 
which accumulates in $\ol X$ at the (compact) subset $A\subseteq\geo X$.
Then the accumulation set of $(x_n)$ in $\Flagt$ consists only of simplices
which are faces of chambers containing a point of $A$.
\end{lem}
\proof
We may assume that $A$ consists only of one ideal point $\xi$.
We fix a base point $o\in X$ and extend the segments $ox_n$ to rays $o\xi_n$.
Then $\xi_n\to\xi$ in $\geo X$.
Let $\si_n\subset\geo X$ be chambers containing the ideal points $\xi_n$
and let $\tau_n\in\Flagt$ be their faces of type $\taumod$. 
By the definition of flag convergence,
the accumulation set of the $\taumod$-re\-gu\-lar sequence $(x_n)$ in $\Flagt$ equals the accumulation set of the sequence $(\tau_n)$.
Its elements are faces of chambers in the accumulation set of the sequence $(\si_n)$ in $\Flags$.
The chambers in the latter accumulation set contain $\xi$. 
\qed

\subsubsection{A continuity property for Weyl cones}

From Lemma~\ref{lem:nbrgsqcs}, we deduce a continuity property for Weyl cones.
Let again $K<G$ denote a maximal compact subgroup and $o\in X$ its fixed point. 

\begin{lem}
\label{lem:nbctffcns}
For $\Theta,d,r,\eps$ there exists $R$ such that the following holds. 

Let $\tau,\tau'\in\Flagt$ 
and let $x\in V(o,\st(\tau))$ and $x'\in V(o,\st(\tau'))$ be points 
such that the pairs $(o,x)$ and $(o,x')$ are $\Theta$-regular with distance $\geq R$.
Suppose that $d(x,x')\leq d$.

Then $V(o,\st(\tau))\cap B(o,r)$ and $V(o,\st(\tau'))\cap B(o,r)$ have Hausdorff distance $\leq\eps$.
\end{lem}
\proof
We can write $\tau'=k\tau$ with $k\in K$ such that the points $kx$ and $x'$ lie in the same euclidean Weyl chamber with tip at $o$.
Then $d(kx,x')\leq d(x,x')\leq d$
and hence $d(x,kx)\leq 2d$.

The elements $\tilde k\in K$,
for which $V(o,\st(\tau))\cap B(o,r)$ and $V(o,\st(\tilde k\tau))\cap B(o,r)$ have Hausdorff distance $\leq\eps$,
form a neighborhood $U$ of $\Stab_K(\tau)$. 
Lemma~\ref{lem:nbrgsqcs} implies that,
if $R$ is sufficiently large, $k$ must lie in $U$.
\qed

\section{Elementary and unipotent subgroups}

In our relativizations of the Anosov condition a prominent role is played by the stabilizers 
of bounded parabolic limit points.
They are the peripheral subgroups for an induced relatively hyperbolic structure.
It is the presence of such subgroups that distinguishes the {\em relative } from the {\em absolute} case. 
They are regular subgroups with a unique limit point.

In this section we collect geometric and algebraic information about and discuss some examples of subgroups with a unique limit point. 
We will see that they tend to consist of elements with zero infimal displacement.
This leads us to also discussing subgroups with zero infimal displacement,
in particular unipotent subgroups.

\subsection{Cyclic subgroups}
\label{sec:cyc}

We first establish some properties  of cyclic subgroups and their limit sets. 

Let $g\in G$ be non-elliptic and consider the (discrete and free) cyclic subgroup $\<g\><G$.

We first look at the case $m_g>0$  
and let $l$ be an oriented axis of the semisimple part $g_s$  (see section \ref{sec:isomsym}).
The orbits of $\<g\>$ deviate sublinearily from $l$ and their visual limit set is 
$\La(\<g\>)=\geo l$,
cf. \eqref{eq:sblndistorb} and Proposition~\ref{prop:dstorbism}.
If $l$ is $\taumod$-regular,
then $l(\pm\infty)\in\ost(\tau_\pm)$
for a pair of antipodal simplices $\tau_{\pm}\in\Flagt$ 
and $\geot l=\{\tau_-, \tau_+\}$. 

\begin{lem}
\label{lem:cclsbgmgp}
Suppose that $m_g>0$ and $l$ is an axis of $g_s$.
Then:

(i) $\<g\>$ is uniformly $\taumod$-regular if and only if $l$ is $\taumod$-regular.
In this case, 
$\Lat(\<g\>)=\geot l$ is a pair of antipodes.

(ii) If $g$ is hyperbolic and $\<g\>$ is $\taumod$-regular, then $\<g\>$ is uniformly $\taumod$-regular.
\end{lem}

\proof
(i) 
The equivalence follows from \eqref{eq:sblndistorb}, since $l$ contains $g_s$-orbits. 
If the $g$- and $g_s$-orbits are uniformly $\taumod$-regular,
then in view of Lemma~\ref{lem:smtmlim}
they have the same $\taumod$-flag limits, 
$g^{\pm n}\to\tau_{\pm}$ and $g_s^{\pm n}\to\tau_{\pm}$,
and thus $\Lat(\<g\>)=\Lat(\<g_s\>)=\geot l$.

(ii)
If $g$ is hyperbolic, then 
$l$ is an axis of $g$.
Hence, if $\<g\>$ is $\taumod$-regular, then so is $l$.
\qed

\medskip 
If $g$ is non-strictly parabolic,
$\taumod$-regularity of $\<g\>$ does not imply uniform $\taumod$-regularity.
If $\<g\>$ is non-uniformly $\taumod$-regular,
then the axis $l$ of $g_s$ is not $\taumod$-regular.
By Lemma~\ref{lem:vsaccsrgsq},
the limit set $\Lat(\<g\>)$ then consists of simplices $\tau$ in $\Flagt$
which are faces of chambers containing one of the ideal points $l(\pm\infty)$;
both points $l(\pm\infty)$ must be covered.

In the $\simod$-regular case we obtain, supplementing the previous lemma:
\begin{lem}
\label{lem:lstptrn}
If $g$ is non-strictly parabolic and $\<g\>$ is $\simod$-regular, then $|\Las(\<g\>)|\geq2$. 
\end{lem}
\proof 
By Lemma~\ref{lem:vsaccsrgsq}, 
both ideal points $l(\pm\infty)$ lie in a chamber contained in $\Las(\<g\>)$.
\qed

\medskip
We are left with the case when $m_g=0$ and $g$ is strictly parabolic. Then,
according to the discussion in section~\ref{sec:isomsym} leading to \eqref{eq:vlimststrpb},
the orbits of $\<g\>$ are Hausdorff close to horocycles in a totally-geodesic hyperbolic plane $X'\subset X$
and $|\La(\<g\>)|=1.$
We also had defined there  the type $\theta(g)\in \simod$ and the face type 
$\taumod(g)$ of $g$; both are $\iota$-invariant.
\begin{lem}
\label{lem:strpburg}
If $g\in G$ is strictly parabolic, 
then the subgroup $\<g\>$ is uniformly $\taumod$-regular precisely for the face types $\taumod\subseteq\taumod(g)$,
and $|\La_{\taumod}(\<g\>)|=1$ for these $\taumod$.
\end{lem}
\proof
The unique visual limit point $\xi$ of $\<g\>$ spans a simplex $\tau(g)$ of type $\taumod(g)$
and therefore is $\taumod$-regular if and only if $\taumod\subseteq\taumod(g)$.
Hence $\<g\>$ is uniformly $\taumod$-regular precisely for these $\taumod$.
From the definition of flag convergence
it follows that $\La_{\taumod}(\<g\>)$ consists of the type $\taumod$ face of $\tau(g)$.
\qed

\subsection{Elementary subgroups}

Recall that an antipodal $\taumod$-regular subgroup $\Ga<G$ is {\em $\taumod$-elementary} if $|\LatGa|\leq2$. 
Since $\Ga$ is discrete, it holds that $\LatGa=\emptyset$ if and only if $\Ga$ is finite. 

For cyclic subgroups, we saw in section~\ref{sec:cyc}:
\begin{example}
Uniformly $\taumod$-regular cyclic subgroups $\<g\><G$ are $\taumod$-antipodal elementary.
Moreover,
$|\Lat(\<g\>)|=1$ if and only if $g$ is strictly parabolic. 
\end{example}

Further examples are provided by rank one symmetric subspaces and products of rank one symmetric spaces:

\begin{example}\label{ex:elementary2}
Let $G_1<G$ and $X_1\subset X$ be as in Example~\ref{ex:rank1embedding}.
Suppose that $\Ga_1< G_1$ is a discrete subgroup 
which consists entirely of parabolic and elliptic elements, equivalently, which preserves a horosphere 
$Hs_1\subset X_1$,
whose center we denote by $\zeta_1\in\geo X_1\subset\geo X$.
Then $\Ga_1$ has visual limit set $\La(\Ga_1)=\{\zeta_1\}$ 
and is uniformly $\taumod$-regular with $\LatGa=\{\tau_{\zeta_1}\}$ 
for the face types $\taumod\subseteq\taumod(X_1)$. 
\end{example}

\begin{example}\label{ex:elementary4}
Consider the product space $X=X_1\times X_2=\H^2\times \H^2$. 
In this case, $\simod$ is an arc of length $\pihalf$,
and we denote by $\taumod^i$ the vertex of $\simod$ corresponding to the hyperbolic plane factor $X_i$.
Then $\Flags\cong\Flag_{\taumod^1}\times\Flag_{\taumod^2}$ with $\Flag_{\taumod^i}=\geo X_i\cong S^1$.

A {\em non-strictly parabolic isometry} of $X$ has, up to switching the factors, the form $g=(g_1, g_2)$ 
with $g_1\in \Isom(X_1)$ hyperbolic (with two ideal fixed points $\la_{\pm}\in\geo X_1$) 
and $g_2\in \Isom(X_2)$ parabolic (with unique ideal fixed point $\mu\in\geo X_2$). 
The subgroup $\<g\><G$ is $\taumod$-regular and $\taumod$-elementary for all face types $\taumod\subseteq\simod$, 
but its uniformity and antipodality depend on $\taumod$: 
It is non-uniformly $\simod$-regular
with $\Las(\<g\>)=\{(\la_-,\mu), (\la_+,\mu)\}$
and hence {\em not $\simod$-antipodal}.
It is uniformly $\taumod^1$-regular with $\La_{\taumod^1}(\<g\>)=\{\la_-,\la_+\}$
and hence {\em $\taumod^1$-antipodal}.
And it is non-uniformly $\taumod^2$-regular with $\La_{\taumod^2}(\<g\>)=\{\mu\}$
and hence also {\em $\taumod^2$-antipodal}.
\end{example}

\medskip
In the case of two antipodal limit points and $\taumod=\simod$, we can say in general:
\begin{prop}
If $\Ga<G$ is $\simod$-regular antipodal with $|\Las|=2$.
Then $\Ga$ is virtually cyclic and contains only semisimple elements. 
\end{prop}
\proof
There exists a $\Ga$-invariant maximal flat $F\subset X$ on which $\Ga$ acts by translations. 
Since $|\Las|=2$, $\Ga$ must be virtually cyclic.
\qed

\medskip
In the $\taumod$-regular case,
the algebraic conclusion (virtually cyclic) still holds 
and $\Ga$ must be uniformly $\taumod$-regular, 
see \cite[Lemma 5.45]{anolec}.  

We now turn to discussing subgroups with a unique limit point.

\subsection{Unique limit point versus zero infimal displacement}

We begin with a geometric property of subgroups with a unique limit point:

\begin{lem}
\label{lem:parabolic}
(i)
If $\Ga< G$ is $\simod$-regular and $|\LasGa|=1$,
then all elements of $\Ga$ are have zero infimal displacement number,
equivalently, are elliptic or strictly parabolic.

(ii)
If $\Ga< G$ is uniformly $\taumod$-regular and $|\LatGa|=1$,
then the same conclusion holds.
\end{lem}
\proof 
This is a direct consequence of Lemmas~\ref{lem:cclsbgmgp} and~\ref{lem:lstptrn}.
\qed

\medskip
We now discuss consequences of zero infimal displacement. {For this discussion we will use some of the material from \cite{anolec}. 
Given a spherical chamber $\si\subset \geo X$, recall that the stabilizer of a chamber $\si\in\Flags$ is a minimal parabolic subgroup $P_{\si}<G$. We refer the reader to  \cite[sections 2.10 and 2.11, Remark 2.28]{anolec} for the discussion of horocycles in $X$, the horocyclic foliation of $X$ by the {\em horocycles at a spherical chamber} $\si\subset \geo X$ and the {\em horocyclic subgroup} $N_{\si}\triangleleft P_{\si}$, which is the common stabilizer of all horocycles at $\si$. Algebraically speaking,  $N_\si$ decomposes as the semidirect product $N_{\si}=U_{\si}\rtimes K_{\si,\hat\si}$,  where $U_{\si} \triangleleft P_{\si}$ is the unipotent radical of $P_{\si}$, 
$\hat\si$ is a chamber in $\geo X$ opposite to $\si$, 
and $K_{\si,\hat\si}$ is the pointwise stabilizer in $G$ of the maximal flat $F\subset X$ containing $\si, \hat\si$ in its visual boundary. 
}

{With these preliminaries out of the way, now turn to arbitrary subgroups with a fixed point on $\Flags$.} 

\begin{prop}
If a (not necessarily discrete) subgroup $\Ga<P_{\si}$ consists of elements with zero infimal displacement,
then it preserves every horocycle based at $\si$,
i.e. $\Ga<N_{\si}$.
\end{prop}

\proof
$P_{\si}$ preserves the (transversely Riemannian) foliation of $X$ by horocycles based at $\si$.
Every maximal flat $F\subset X$ asymptotic to $\si$ is a cross section to this foliation 
and hence is naturally isometric to the leaf space.
The action of $P_{\si}$ on the leaf space is by translations. 
It follows that elements with zero infimal displacement number act trivially on it,
i.e. preserve every horocycle at $\si$.
\qed

\medskip
This proposition has the following algebraic consequence for discrete subgroups:

\begin{cor}
If in addition $\Ga$ is discrete, then it is finitely generated and virtually nilpotent.
\end{cor}
\proof
This follows from Auslander's theorem (see Theorem \ref{thm:auslander} in the appendix).
\qed

\medskip
We conclude for $\simod$-regular subgroups with unique limit point:
\begin{cor}
If $\Ga<G$ is  $\simod$-regular and $\LasGa=\{\si\}$,
then $\Ga<N_{\si}$ and therefore $\Ga$ is finitely generated and virtually nilpotent.
\end{cor}

\medskip
More generally, without a fixed point assumption,
one can deduce from the work of Prasad \cite{Prasad} or already from Tits \cite{Tits}:

\begin{thm}\label{thm:purezero}
Suppose that $\Ga< G$ is a subgroup 
consisting only of elements with zero infimal displacement.
Then there exists a unipotent Lie subgroup $N<G$
and a compact subgroup $K_N< G$ normalizing $N$,
such that $\Ga$ is contained in $N \rtimes K_N$.
\end{thm}
\proof
To relate our condition of zero infimal displacement to the one used by Prasad,
we note that $m_g=0$ for $g\in G$ 
if and only if the transvection component $g_t$ in the Jordan decomposition \eqref{eq:Jordan} is trivial,
equivalently, if the adjoint action of $g$ has all eigenvalues in $S^1$.

Consider now the Zariski closure $\ol{\Ga}<G$  of $\Ga$.
Let $N$ denote the unipotent radical of the identity component $\ol{\Ga}_0$ of
$\ol{\Ga}$.
Then the projection $\Ga'$ of $\Ga$ to $G'=\ol{\Ga}/N$ still consists only of
elements of zero displacement and is Zariski dense in $G'$. 

We claim that $G'$ is compact. 
If not, then a theorem by Prasad \cite{Prasad} implies that 
$\Ga'$ contains elements $g$ whose adjoint action has eigenvalues outside the unit circle,
a contradiction. 

Hence, $G'$ is compact and we obtain that $\Ga< N \rtimes K_N$ with $K_N\cong G'$.
\qed

\medskip
We conclude for uniformly $\taumod$-regular subgroups with unique limit point:
\begin{cor}\label{cor:nilfg}
If $\Ga< G$ is $\taumod$-uniformly regular and $|\Lat(\Ga)|=1$,
then there exists a unipotent Lie subgroup $N<G$
and a compact subgroup $K_N< G$ normalizing $N$,
such that $\Ga$ is contained in $N \rtimes K_N$. 
In particular, $\Ga$ is finitely generated and virtually nilpotent.  
\end{cor}

\subsection{Unipotent subgroups}\label{sec:unipotent}

A natural class of zero displacement subgroups is provided by unipotent subgroups.
We now look at their limit sets in some examples.

In rank one, unipotent subgroups always have a unique limit point. 
We also know that, in higher rank, unipotent {\em one-parameter} subgroups 
are $\taumod$-regular\footnote{I.e. their orbits are $\taumod$-regular subsets.} 
for some type $\taumod$ depending on the subgroup
and have a single $\taumod$-limit point. 
In contrast, as we will see,
unipotent subgroups of dimension $\geq2$ in higher rank are not necessarily $\taumod$-regular for any $\taumod$,
and even if they are uniformly $\taumod$-regular, they may fail to be $\taumod$-elementary.

We now discuss this in the case of $G=SL(3,\R)$.

We begin with one-parameter unipotent subgroups.
There are two conjugacy classes of such subgroups.
Each subgroup of either type is contained in a Lie subgroup locally isomorphic to $SL(2,\R)$
and preserves a totally geodesic hyperbolic plane of $\theta$-type equal to the midpoint $\bar\mu$ of the Weyl arc $\simod$.
The subgroups conjugate to the group $U_1$ consisting of the elements
$$
\left ( \begin{array}{ccc}
1&&t\\
&1&\\
&&1\end{array}
\right) 
$$
are contained in $SL(2,\R)\subset SL(3,\R)$ (reducibly embedded).
The subgroups conjugate to the group $V_1$ consisting of the elements
$$
\left ( \begin{array}{ccc}
1&t&\frac{t^2}{2}\\
&1&t\\
&&1\end{array}
\right) 
$$
are contained in $SO(2,1)\subset SL(3,\R)$ (irreducibly embedded).

The unique limit flags of these subgroups can be determined as follows.
Consider the unipotent subgroup $\exp(\R\cdot n)$ for a nilpotent element $n\in sl(3,\R)$.
If $\rank(n)=1$, then the limit flag equals $\im(n)\subset\ker(n)$,
and if $\rank(n)=2$, it equals $\im(n^2)\subset\ker(n^2)$.

The geometry of the $U_1$-orbit foliation of $X$ is particularly nice:
Since the normalizer of $U_1$ contains a minimal parabolic subgroup
and therefore acts transitively on $X$,
this foliation is homogeneous, i.e. any two $U_1$-orbits are congruent.
As a consequence, the $\De$-distance of any pair of points in any $U_1$-orbit has type $\bar\mu$,
i.e. lies on the bisector of $\De$. {(Recall that $\bar\mu$ is the midpoint of the arc $\simod$.) } 
The $V_1$-orbit foliation does not have either of these properties.

Now we turn to two-parameter unipotent subgroups.
There are three conjugacy classes 
represented by the subgroups $U_2^{\pm}$ and $V_2$ consisting of the elements
$$
\left( \begin{array}{ccc}
1&\star &\star \\
&1&\\
&&1\end{array}
\right) 
\quad \hbox{ , }  \quad 
\left( \begin{array}{ccc}
1&&\star \\
&1&\star \\
&&1\end{array}
\right)
\quad \hbox{ and }  \quad 
\left( \begin{array}{ccc}
1&t&s \\
&1&t \\
&&1\end{array}
\right),
$$
respectively.
Note that $U_2^{\pm}$ are conjugate inside the full isometry group of $X$. 

Again the foliations of $X$ by $U_2^{\pm}$-orbits have nice geometry, 
even though they are no longer homogeneous:
Since all one-parameter subgroups of $U_2^{\pm}$ are conjugate to $U_1$,
the $\De$-distance of any pair of points in any $U_2^{\pm}$-orbit still has type $\bar\mu$.
In particular, $U_2^{\pm}$ is {\em uniformly $\simod$-regular}.

However, the subgroups $U_2^{\pm}$ have {\em large limit sets}:
One verifies that they consist of the limit points of their one-parameter subgroups.
In the case of $U_2^+$ these are the flags of the form $\<e_1\>\subset E^2$,
and in the case of $U_2^-$ the flags of the form $L^1\subset \<e_1,e_2\>$.

We note that the same discussion applies to unipotent subgroups of $SL(n,\R)$ of the form
$$ 
\left( \begin{array}{cccc}
1&\star& \ldots &\star\\
&\ddots&  &\\
&          & 1& \\
&          &    &    1\end{array}
\right)
\quad \hbox{ and }  \quad 
\left( \begin{array}{cccc}
1&         & &\star\\
&\ddots&  &\vdots\\
&          & 1& \star\\
&          &    &    1\end{array}
\right).
$$

Returning to $SL(3,\R)$,
in contrast, the subgroup $V_2$ is {\em not $\simod$-regular} (and hence neither are its lattices).
This is a consequence of the following fact about the non-regularity of sequences in the full horocyclic subgroup:
Any diverging sequence of elements
$$
\left( \begin{array}{ccc}
1&x_n&\\
&1&y_n\\
&&1\end{array}
\right)
$$
where $0< c\leq\frac{|x_n|}{|y_n|}\leq C$ 
is not $\simod$-regular,
as one can see from its dynamics on $\R P^2$.

In conclusion,
$SL(3,\R)$ contains no non-cyclic discrete $\simod$-regular elementary unipotent subgroups.

\section{Finsler-straight paths and maps in symmetric spaces} 
\label{sec:finsstraight}

In this section we introduce a notion of Finsler-straightness 
which adapts the notion of straightness in Gromov hyperbolic spaces 
discussed earlier in section~\ref{sec:HypSpaces}
to the geometry of higher rank symmetric spaces. 
This notion of straightness can be regarded as a regularity condition and 
is implicit in our earlier work on Morse quasigeodesics \cite{morse}. 
We will use it later on to define relative versions of our notions of Morse (equivalently, Anosov) subgroups,
namely the notions of {\em relatively Morse}, see section~\ref{sec:relmo},
and {\em relatively straight} subgroups, see section~\ref{sec:str=ae}.

The main results of this sections are Propositions~\ref{prop:bdmcnc} and~\ref{prop:fllbdmpay}
dealing with extensions of straight maps to infinity.
They will be used in section~\ref{sec:fstrctns} to construct boundary embeddings for straight subgroups.

\subsection{Triples}

We denote by $T(X):=X^3$ the space of triples of points in $X$ and by 
\begin{equation*}
T(X,\Flagt):=(X\sqcup\Flagt)\times X\times (X\sqcup\Flagt)
\end{equation*}
the space of ideal triples in the Finsler bordification $X\sqcup\Flagt$
with middle point in $X$.

We first define straightness for (non-ideal) triples in $X$:

\begin{definition}[Finsler-straight triple] 
\label{defn:diamondstraight}
{Fix an $\iota$-invariant, compact, Weyl-convex subset $\Theta$ of the open star $\ost(\taumod)\subset\simod$.} 
A triple $(x_-,x,x_+)\in T(X)$ is called

(i) 
{\em $(\Theta,d)$-straight}, 
$d\geq0$, 
if the points $x_-,x$ and $x_+$ are $d$-close to points $x_-',x'$ and $x_+'$,
respectively, 
which lie in this order on a $\Theta$-Finsler geodesic.

(ii) {\em $(\taumod,d)$-straight} if the same property holds with $\Theta$ replaced by $\taumod$.
\end{definition}

In particular, 
a triple $(x_-,x,x_+)$ is $(\Theta,0)$-straight if and only if 
the points $x_-,x,x_+$ lie in this order on a $\Theta$-Finsler geodesic.

Finsler-straightness is stable under perturbation:
Any triple $(\hat x_-,\hat x,\hat x_+)$ which is $r$-close to a $(\Theta,d)$-straight triple $(x_-,x,x_+)$, 
i.e. $d(x_-,\hat x_-),d(x,\hat x),d(x_+,\hat x_+)\leq r$,
is $(\Theta,d+r)$-straight.

\medskip
It is useful to note that (modulo doubling the constant $d$)
the nearby Finsler geodesic in the definition can be chosen through one of the endpoints of the triple:
\begin{lem}
\label{lem:thrghnndpt}
If $(x_-,x,x_+)$ is $(\Theta,d)$-straight,
then the points $x$ and $x_+$ are $2d$-close to points $x''$ and $x_+''$,
respectively, 
such that $x_-,x''$ and $x_+''$ lie in this order on a $\Theta$-Finsler geodesic.

The same assertion holds with $\Theta$ replaced by $\taumod$.
\end{lem}
Proof:
The points $x'$ and $x_+'$ in the definition of Finsler-straightness 
are contained in a $\taumod$-Weyl cone $V(x_-',\st(\tau_+))$.
The Weyl cone $V(x_-,\st(\tau_+))$ asymptotic to it has Hausdorff distance $\leq d(x_-,x_-')\leq d$.
It can be represented as the image $V(x_-,\st(\tau_+))=gV(x_-',\st(\tau_+))$ by an isometry $g\in G$
fixing $\tau$ at infinity and mapping $x_-'\mapsto x_-$.
The isometry $g$ has displacement $\leq d$ on the entire cone $V(x_-',\st(\tau_+))$.
We put $x''=gx'$, $x_+''=gx_+'$ and choose the $\Theta$-Finsler geodesic through $x_-$ as the $g$-image of the one through $x_-'$
given by the definition.\footnote{The points $x''$ and $x_+''$ can also be described as follows:
The point $x'$ lies on a Riemannian ray $x_-'\xi$ asymptotic to $\xi\in\st(\tau_+)$.
We choose $x''\in x_-\xi$ with $d(x_-,x'')=d(x_-',x')$.
The point $x_+''$ is constructed similarly.}
\qed

\medskip
Note that the Finsler geodesic in the conclusion of the lemma can be chosen as a Finsler segment $x_-x_+''$
through $x''$ and is then contained in a Weyl cone $V(x_-,\st(\tau_+))$, $\tau_+\in\Flagt$.

To show that the nearby Finsler geodesic can be chosen through both endpoints of the triple
and to control its distance from the middle point,
takes more effort.\footnote{Weyl cones vary 1-Lipschitz continuously with their tips,
whereas we do not have such a result for diamonds at our disposal in full generality.}

\medskip
The notion of Finsler-straightness naturally extends to {\em ideal} triples in 
$T(X,\Flagt)$: 

We say that a triple $(x_-,x,\tau_+)\in X^2\times\Flagt$ is $(\Theta,d)$-straight 
if the points $x_-$ and $x$ are $d$-close to points $x_-'$ and $x'$, respectively, 
such that $x'$ lies on a $\Theta$-Finsler ray $x_-'\tau_+$.
This ray is then $d$-Hausdorff close to a $\Theta$-Finsler ray $x_-\tau_+$,
compare the proof of Lemma~\ref{lem:thrghnndpt},
and the latter passes within distance $2d$ from $x$.
We say that $(x_-,x,\tau_+)$ is 
$(\taumod,d)$-straight if the same property holds with $\Theta$ replaced by $\taumod$.
Analogously for triples $(\tau_-,x,x_+)\in \Flagt\times X^2$.

Similarly,
we say that a triple $(\tau_-,x, \tau_+)\in\Flagt\times X\times\Flagt$ is $(\taumod,d)$-straight 
if the simplices $\tau_{\pm}$ are antipodal 
and $x$ lies within distance $d$ of a $\taumod$-Finsler line $\tau_-\tau_+$.
Note that $(\Theta,d)$-straightness would be an equivalent property for triples $(\tau_-,x, \tau_+)$, 
because $\tau_-\tau_+$ can be chosen $\Theta$-regular,
which is why we do not introduce it.

\subsection{Paths}
\subsubsection{Holey rays and lines}\label{sec:Holey rays and line}

As for Gromov hyperbolic spaces,
we call a map 
$q:H\to X$ from a subset of $H\subset\R$ 
a {\em holey line}.
If $H$ has a minimal element, we also call $q$ a {\em holey ray}.
(The domains of holey rays will usually be denoted $H_0$ below.)  
A {\em sequence} $(x_n)_{n\in\N}$ in $X$ can be regarded as a holey ray $\N\to X$.

{In the next section we will frequently use {\em linear interpolation} of holey lines and rays. Namely, assume that $H$ is a closed and discrete  subset of $\R$. For each pair of consecutive elements $h_i, h_{i+1}\in H$ we extend the map $q$ to the interval $[h_i, h_{i+1}]$ by the constant speed parameterization of the Riemannian geodesic $q(h_i)q(h_{i+1})$. We will use the notation ${\mathfrak q}$ for the resulting map.}  

{Suppose that there exists a simplex $\tau$ in $\geo X$ of type $\taumod$ such that for each 
pair of points $h_i< h_j$ in $H$ we have that $q(h_j)$ lies in the cone $V(q(h_i), \st(\tau))$. Then the 
 {\em nested cones property} (see Appendix 1) implies that the map ${\mathfrak q}$ is a $\taumod$-regular Finsler geodesic in $X$. Furthermore, if each segment $q(h_i)q(h_{i+1})$ is $\Theta$-regular for some Weyl-convex compact $\Theta\subset \ost(\taumod)$, then the entire ${\mathfrak q}$ is 
 $\Theta$-regular. }

We will consider extensions to infinity $$\ol q:\ol H:=H\sqcup\{\pm\infty\}\to X\sqcup\Flagt$$ 
of holey lines $q:H\to X$
by sending $\pm\infty$ to simplices $\tau_{\pm}\in\Flagt$, 
and refer to $\ol q$ as an {\em extended holey line}.
In the case of holey rays $q:H_0\to X$,
we consider extensions $\ol q:\ol H_0:=H_0\sqcup\{+\infty\}\to X\sqcup\Flagt$ 
by sending $+\infty$ to a simplex $\tau\in\Flagt$,
and refer to $\ol q$ as an {\em extended holey ray}.

We carry over the notion of Finsler-straightness from triples to holey lines 
by requiring it for all triples in the image:
\begin{definition}[Finsler-straight holey line] 
\label{def:fstrspm}
A holey line $q:H\to X$ is called 

(i) {\em $(\Theta,d)$-straight} if the triples $(q(h_-), q(h), q(h_+))$ 
are $(\Theta,d)$-straight for all $h_-\leq h\leq h_+$.

(ii) {\em $(\taumod,d)$-straight} if the same property holds with $\Theta$ replaced by $\taumod$.
\end{definition}

We say that $q$ is {\em $\Theta$-straight} if it is $(\Theta,d)$-straight for some $d$,
analogously for {\em $\taumod$-straight},
and that $q$ is {\em uniformly $\taumod$-straight} if it is $\Theta$-straight for some $\Theta$.

\begin{rem}
\label{rem:fstr}
(i) 
Finsler-straightness is preserved under restriction to subsets of $H$. 

(ii) 
Finsler-straightness is stable under perturbation:
If two holey lines $q,q':H\to X$ are $r$-close,
$d(q(h),q'(h))\leq r$ for all $h\in H$,
and $q$ is $(\Theta,d)$-straight,
then $q'$ is $(\Theta,d+r)$-straight.

(iii)
A holey line $q:H\to X$ is $(\Theta,0)$-straight
if and only if $q$ maps monotonically into a $\Theta$-Finsler geodesic.
\end{rem}

Similarly, we say that 
an extended holey line $\ol q:\ol H\to X\sqcup\Flagt$ is {\em $(\Theta,d)$-straight} 
if all triples $(\ol q(h_-),q(h),\ol q(h_+))$ in $X\sqcup\Flagt$ 
for $-\infty\leq h_-\leq h\leq h_+\leq +\infty$ in $\ol H$ with 
$h\in H$ 
are $(\Theta,d)$-straight,
and analogously in the ray case.
The properties 
{\em $(\taumod,d)$-straight}, {\em $\Theta$-straight}, {\em $\taumod$-straight} and {\em uniformly $\taumod$-straight}
are then defined in the obvious way.

\medskip
The ($\taumod$-)straightness of an extended holey ray $\ol q:\ol H_0 \to X\sqcup\Flagt$ implies 
that $q(H_0)$ lies within distance $2d$ of the Weyl cone $V(x_-,\st(\tau_+))$
with $x_-=q(\min H_0)$ and $\tau_+=\ol q(+\infty)$,
however it does not imply flag-convergence $q(h)\to\tau_+$ as $h\to\sup H_0$
due to possible lack of regularity.
If also $q$ is $\taumod$-regular, then $q(h)\to\tau_+$ conically. 
The straightness of an extended holey line $\ol q:\ol H \to X\sqcup\Flagt$ implies 
that the simplices $\tau_{\pm}=\ol q(\pm\infty)\in\Flagt$ are antipodal and 
the image $q(H)$ lies within distance $d$ from the parallel set $P(\tau_-,\tau_+)$.

\subsubsection{Asymptotics at infinity}

We now discuss the {\em convergence at infinity} of Finsler straight holey rays. {As before, we fix an $\iota$-invariant, compact, Weyl-convex subset $\Theta$ of the open star $\ost(\taumod)\subset\simod$.} 

To obtain flag-convergence, one needs to impose in addition regularity.

\begin{lem}
\label{lem:smccstx}
Let $x_n,x'_n\to\infty$ be $\taumod$-regular sequences in $X$ 
such that the triples $(x,x_n,x'_n)$ are $(\taumod,d)$-straight for some base point $x\in X$ and $d\geq0$.

Then the sequences $(x_n)$ and $(x'_n)$ have the same accumulation set in $\Flagt$.
\end{lem}
\proof
By straightness,
there exists a sequence $(\tau_n)$ in $\Flagt$ 
so that the points $x_n,x'_n$ are contained in the $2d$-neighborhood of the Weyl cone $V(x,\st(\tau_n))$ for all $n$, 
see Lemma~\ref{lem:thrghnndpt}.
It follows that the $\taumod$-flag accumulation sets of both sequences $(x_n)$ and $(x'_n)$ in $\Flagt$
coincide with the accumulation set of $(\tau_n)$. 
\qed

\medskip
It follows that regular Finsler-straight holey rays converge at infinity,
as long as they are unbounded. 
(Note that we allow ``infinite holes'' and put no restriction on the ``speed''.) 
Here, we call a holey ray or line {\em $\taumod$-regular} if its image in $X$ is a $\taumod$-regular subset. 

\begin{cor}
\label{cor:rgstrhrcnv}
If $q:H_0\to X$ is a $\taumod$-regular $\taumod$-straight holey ray
with unbounded image $q(H_0)$,
then it $\taumod$-flag converges at infinity,
$q(h)\to\tau\in\Flagt$  as $h\to\sup H_0$.
\end{cor}
\proof
Since $q(H_0)$ is unbounded, 
there exists a sequence $h_n\nearrow\sup H_0$ in $H_0$ so that the sequence $(q(h_n))$ in $X$ diverges and hence is $\taumod$-regular.
By the compactness of $\Flagt$,
after passing to a subsequence, 
it flag converges, $q(h_n)\to\tau\in\Flagt$.

If $h'_n\to\sup H_0$ is another sequence in $H_0$,
there exists a sequence of indices $m_n\to+\infty$ in $\N$ growing slowly enough so that $h_{m_n}\leq h'_n$ for large $n$.
The triples $(x,q(h_{m_n}),q(h'_n))$ are then $(\taumod,d)$-straight for some base point $x\in X$ and $d>0$.
By Lemma~\ref{lem:smccstx},
also $q(h'_n)\to\tau$.
\qed

\medskip
In the {\em uniformly} regular case, we can show that the flag-convergence is {\em conical}:
\begin{lem}
\label{lem:rgstrhrcnvcncl}
If $q:H_0\to X$ is a $(\Theta,d)$-straight holey ray
with unbounded image $q(H_0)$,
then it conically $\taumod$-flag converges at infinity,
$q(h)\to\tau\in\Flagt$  as $h\to\sup H_0$.
More precisely,
the extended holey ray $\ol q:\ol H_0\to X\sqcup\Flagt$ with $\ol q(+\infty)=\tau$ 
is still $(\Theta,2d)$-straight.
\end{lem}
\proof
Let $h_0=\min H_0$ and $o=q(h_0)$.

By straightness,
for any $h<h'$ in $H_0$ there exists a $\taumod$-Weyl cone with tip at $o$
which intersects both balls $\ol B(q(h),2d)$ and $\ol B(q(h'),2d)$.

With Lemma~\ref{lem:nbctffcns} it follows that
for every $\eps>0$ and every $h_1\in H_0$ there exists $h_2>h_1$ in $H_0$ with the property:
If $h\leq h_1<h_2\leq h'$, then {\em every} $\taumod$-Weyl cone with tip at $o$,
which intersects $\ol B(q(h'),2d)$, also intersects $B(q(h),2d+\eps)$.

Now we take a sequence $h'_n\to\sup H_0$
and let $\tau_n\in\Flagt$ so that $q(h'_n)\in V(o,\st(\tau_n))$.
Then for every $h\in H_0$ and $\eps>0$ the Weyl cone $V(o,\st(\tau_n))$ 
intersects $B(q(h),2d+\eps)$ for all sufficiently large $n$.
It follows that $(\tau_n)$ converges, $\tau_n\to\tau\in\Flagt$,
and that $q(H_0)$ is contained in $\ol N_{2d}(V(o,\st(\tau)))$.
\qed

\medskip
{From now on, we will only consider holey lines and rays $q:H_0\to X$ such 
that $H_0\subset \R$ is closed and discrete.}

For {\em extended} Finsler-straight holey rays,
already a weaker uniformity assumption implies conical flag-convergence,
more precisely, closeness to a Finsler geodesic:

\begin{claim}
\label{claim:strghtphrgclfns}
For $d,\phi$ there exists $r$ such that:

If $q:H_0\to X$ is a holey ray
which admits a 
$(\taumod,d)$-straight extension $\ol q$
and is $(\taumod,\phi)$-regular\footnote{see Definition \ref{defn:WUR}}, 
then there exists a $\taumod$-Finsler ray $q(0)\tau$ and a monotonic map $q':H_0\to q(0)\tau$ {whose image} 
{\mini which} is $r$-close to {$q(H_0)$},
where $\tau=\ol q(+\infty)$.
\end{claim}
\proof
By straightness, $q(H_0)\subset\ol N_{2d}(V(q(0),\st(\tau)))$.

Let $q'':H_0\to V(q(0),\st(\tau))$ be a map  $2d$-close to $q(H_0)$,
e.g. the {composition of $q$ with the} nearest point projection to the Weyl cone {$V(q(0),\st(\tau))$}.
We extend $q''$ to infinity by $\ol q''(+\infty):=\tau$.
Then $\ol q''$ is $(\taumod,3d)$-straight and $(\taumod,\phi-4d)$-regular.

The straightness of $\ol q''$ implies that, 
for $h_1<h_2$ in $H_0$,
the point $q''(h_2)$ lies within distance $6d$ of the subcone $V(q''(h_1),\st(\tau))\subset V(q(0),\st(\tau))$.
We wish to show that it is contained in the subcone, 
provided that its distance from the tip $q''(h_1)$ is sufficiently large.
To do so,
let $s>0$ so that $\phi(s)>10d$.
Then,
if $d(q(h_1),q(h_2))>s$, the $(\taumod,\phi)$-regularity of $q$ and the choice $\phi(s)>10d$ imply 
 that 
$$
d(d_{\De}(q''(h_1),q''(h_2)),\Dt\De)\geq\phi(s)-4d>6d.
$$
It follows that $q''(h_2)$ has distance $>6d$ from the boundary of the cone $V(q''(h_1),\st(\tau))$
which forces it to lie inside it,
$$q''(h_2)\in V(q''(h_1),\st(\tau)) .$$

Now let $H_0^s\subset H_0$ be a maximal subset containing $0$ such that $q(H_0^s)$ is $s$-spaced\footnote{see section \ref{sec:metric_spaces} for the definition}.
Then $q(H_0)$ is $s$-Hausdorff close to $q(H_0^s)$.
By the above,
for $h_1<h_2$ in $H_0^s$, it holds that $q''(h_2)\in V(q''(h_1),\st(\tau))$.
{It follows (see the discussion at the beginning of section \ref{sec:Holey rays and line}) 
 that $q''|_{H_0^s}$ maps monotonically into some $\taumod$-Finsler ray 
 $q(0)\tau\subset V(q(0),\st(\tau))$, the image of the linear interpolation ${\mathfrak q}''$ of $q''$.  
 Restricting ${\mathfrak q}''$ to $H_0$ we obtain 
 a monotonic map $q':H_0\to q(0)\tau$ whose image is is $(2d+s)$-close to that of $q$.}
Thus, the assertion holds with $r=2d+s$. 
\qed

\medskip
A corresponding result for holey lines $q:H\to X$ is readily derived:
\begin{add}
\label{add:strghtphrgclfns}
For $d,\phi$ there exists $r$ such that:

If $q:H\to X$ is a holey line
which admits a 
$(\taumod,d)$-straight extension $\ol q$
and is $(\taumod,\phi)$-regular, 
then there exists a $\taumod$-Finsler line {$\varphi=\tau_-\tau_+$} and a monotonic map $q':H\to \tau_-\tau_+$ {\mini which} {whose image 
 is $r$-close to that of $q$,}  where $\tau_{\pm}=\ol q(\pm\infty)$.
\end{add}
\proof
Pick some $h_0\in H$ and a point $q'(h_0)\in P(\tau_-,\tau_+)$ within distance $d$ from $q(h_0)$.
Applying Claim~\ref{claim:strghtphrgclfns} to the two holey subrays of $q$ starting in $q(h_0)$
yields monotonic maps into suitable $\taumod$-Finsler rays $q(h_0)\tau_{\pm}$.
The latter are $d$-Hausdorff close to two $\taumod$-Finsler rays $q'(h_0)\tau_{\pm}$ 
which together form a $\taumod$-Finsler line $\tau_-\tau_+$ with the desired property.
\qed

\medskip
When there are no arbitrarily large holes {in $q(H)$}, regularity can be promoted to {\em uniform} regularity.

{\mini We say that a holey line $q: H\to X$ is {\em coarsely $l$-connected}, $l\geq0$,
if for any $h<\tilde h$ in $H$ there exists a sequence $h_0=h<h_1<\ldots<h_n=\tilde h$ in $H$ 
such that $d(q(h_{i-1}),q(h_i))\leq l$ for all $i$.}

{We say that for $l\geq0$ a holey line $q: H\to X$ is {\em coarsely $l$-connected}, 
if for every pair of consecutive elements $h_i, h_{i+1}\in H$, 
$d(q(h_{i}),q(h_{i+1}))\leq l$. The following claim is straightforward and is left to the reader:}

\begin{claim}
{Suppose that $q$ is coarsely $l$-connected and $H_s\subset H$ is a maximal subset such that $q(H_s)$ is $s$-spaced. 
Then the restriction $q|_{H_s}$ is coarsely $(2s+l)$-connected.
}
\end{claim}

\begin{claim}
\label{claim:lngsgmnfrg}
For $d,\phi,l$ there exist $\Theta, r'$ such that:

If $q$ is as in Addendum~\ref{add:strghtphrgclfns}
and moreover coarsely $l$-connected,
then the $\taumod$-Finsler line $\tau_-\tau_+$ {in Addendum~\ref{add:strghtphrgclfns}} can be chosen to be $\Theta$-regular.
\end{claim}
\proof {Let $q'$ be as in as in Addendum~\ref{add:strghtphrgclfns}.} 
Take $r=r(d,\phi)$ as in Addendum~\ref{add:strghtphrgclfns},
and choose $s,a>0$ so that $\phi(s) \geq 2r+a$.

{Consider any pair $h, \tilde h\in H$ satisfying $h<\tilde h$ and $d(q(h),q(\tilde h))\geq s$. Recall that $d(q(h), q'(h))\le r$ and  $d(q(\tilde h), q'(\tilde h))\le r$. Thus, by the triangle inequality for $\Delta$-distances,
$$
||d_{\De}(q'(h),q'(\tilde h))- d_{\De}(q(h),q(\tilde h))||\le 2r. 
$$
Since the map $q$ is assumed to be $(\taumod,\phi)$-regular, the distance from $d_{\De}(q(h),q(\tilde h))$ to $\Dt\De$ is $\ge \phi(s)\ge 2r+a$. 
It follows that the distance from 
$d_{\De}(q'(h),q'(\tilde h))\in \De$ to $\Dt\De$ is $\geq a$. Assuming, in addition, that $d(q'(h), q'(\tilde h))\le R$, we obtain that the vector $d_{\De}(q'(h),q'(\tilde h))$ is $\Theta$-regular for some 
Weyl-convex compact subset $\Theta=\Theta(R,a)$ in $\ost(\tau)$}.

Let $H_s\subset H$ {be a maximal subset such} that $q(H_s)$ is $s$-spaced. 
{As we noted above}, $q|_{H_s}$ is coarsely $(2s+l)$-connected,
and, consequently (since the distance between $q$ and $q'$ is $\le r$), $q'|_{H_s}$ is coarsely $R=2s+l+2r$-connected.  
It follows that for every pair of consecutive points $h, \tilde h$ in $H_s$ the segment $q'(h) q'(\tilde h)$ is $\Theta$-regular. 

According to Addendum~\ref{add:strghtphrgclfns}, $q'$ is a monotonic map to the $\taumod$-Finsler line $\varphi=\tau_-\tau_+$. We define a new Finsler geodesic $\varphi'$ by linear interpolation ${\mathfrak q}': \R\to X$ of  $q'|_{H_s}$. Then $\varphi'$ is a $\Theta$-Finsler geodesic $\varphi'$ in $X$, just as in the proof of Claim \ref{claim:strghtphrgclfns}. Since $\varphi'$ contains the biinfinite sequence $q'(H_s)$ which (forward/backward) converges to $\tau_\pm$, it follows that the Finsler geodesic $\varphi'$ is asymptotic to $\tau_\pm$. 

By the construction of $\varphi'$, for each pair of consecutive points $h_i, h_{i+1}\in H_s$ the image 
${\mathfrak q}'([t_i, t_{i+1}])$ is $R$-Hausdorff close to $\{q'(h_i), q'(h_{i+1})$. It follows that 
$\varphi'$ is $R$-Hausdorff close to $q(H_s)$. Thus, $q(H)$ is $r':=R+l$-Hausdorf close to $\varphi'$. 
%Since $q'$ maps monotonically into the $\taumod$-Finsler line $\varphi=\tau_-\tau_+$ as in Addendum~\ref{add:strghtphrgclfns},
%it follows further that $q'|_{H_s}$ maps monotonically into a (different) $\Theta$-Finsler line $\tau_-\tau_+$.
%The map $q'|_{H_s}$ can, by interpolation, be extended to a monotonic map $q': H\to\tau_-\tau_+$
%which, after suitably increasing $r=r(d,\phi,l)$, is $r$-close to $q$.
\qed

\medskip
If one allows arbitrarily large holes,
one needs an extra assumption to ensure uniform regularity.
We will consider the following condition:
\begin{dfn}[Uniformly regular large holes]
\label{dfn:nfrglhls}
A holey line $q:H\to X$ with $H\subset\R$ closed and discrete 
is said to have {\em $(\hat\Theta,l)$-regular large holes}
if for any two consecutive elements $h<\tilde h$ in $H$ with $d(q(h),q(\tilde h))>l$,
the pair $(q(h),q(\tilde h))$ is $\hat\Theta$-regular.
We say that $q$ has {\em uniformly $\taumod$-regular large holes}
if it has $(\hat\Theta,l)$-regular large holes for some data $\hat\Theta,l$.
\end{dfn}

Then a very similar argument as for Claim~\ref{claim:lngsgmnfrg} yields
that uniformly regular large holes imply uniform regularity 
for the holey lines under consideration:

\begin{claim}
\label{claim:lngsgmnbgh}
For $d,\phi,l,\hat\Theta$ there exist $\Theta,r$ such that:

Let $q$ be as in Addendum~\ref{add:strghtphrgclfns}
and suppose moreover that it has $(\hat\Theta,l)$-regular large holes.
Then the $\taumod$-Finsler line $\tau_-\tau_+$ can be chosen to be $\Theta$-regular.
\end{claim}

\subsubsection{Morse quasigeodesics}

Morse quasigeodesics are a particular class of uniformly Finsler-straight holey lines (with ``bounded holes'') 
which play a prominent role in our earlier work,
see \cite{morse, mlem, anolec}.

\begin{definition}[Morse quasigeodesic]
A {\em $(\Theta,d,L,A)$-Morse quasigeodesic} in $X$ is a $(\Theta,d)$-straight holey line $q:I\to X$ 
which is defined on an interval $I\subseteq\R$ and is an $(L,A)$-qua\-si\-iso\-met\-ric embedding. 
\end{definition}

We will call a $(\Theta,d,L,A)$-Morse quasigeodesic also briefly a {\em $\taumod$-Morse quasigeodesic}.

One can show that $\taumod$-Morse quasigeodesics are,
up to quasiisometric reparameterization, 
uniformly close to 
$\taumod$-Finsler geodesics.
For quasirays and quasilines, this is a consequence of Claim~\ref{claim:strghtphrgclfns} and Addendum~\ref{add:strghtphrgclfns}.
The definition therefore agrees with the definition given in our earlier papers.

The main result of \cite{mlem} is that uniformly $\taumod$-regular quasigeodesics are $\taumod$-Morse.
(Note that the converse holds trivially.)

\subsection{Maps}\label{sec:maps}
\subsubsection{Straight maps}

We now generalize the notion of Finsler-straight holey line
and introduce and study a class of maps from subsets of Gromov hyperbolic spaces into symmetric spaces 
which preserve straightness.
Here we use in the hyperbolic spaces the notion of straightness in terms of closeness to geodesics (cf.\ Definition~\ref{def:strtrp})
and in the symmetric spaces the notion of straightness in terms of closeness to {\em Finsler} geodesics (cf.\ Definition~\ref{defn:diamondstraight})
which is well-adapted to the higher rank geometry.
Straight maps are coarse analogues of {\em projective maps} in Riemannian geometry 
which are smooth maps sending unparameterized geodesics to unparameterized geodesics. 

Let $Y$ be a $\de$-hyperbolic proper geodesic space
and $X$ a symmetric space. 
In the following, $$f:A\to X$$ will always denote a {\em metrically proper} map defined on a subset $A\subset Y$.

\begin{dfn}
[Finsler-straight map]
\label{dfn:fstrghtmp}
The map $f:A\to X$ is called 

(i) {\em $\Theta$-straight} if for every $D$ there exists $d=d(\Theta,D)$ 
such that $f$ sends $D$-straight triples in $A$ to $(\Theta,d)$-straight triples in $X$. 

(ii) {\em $\taumod$-straight} if the same property holds with $\Theta$ replaced by $\taumod$.

(iii) {\em uniformly $\taumod$-straight} 
if it is $\Theta$-straight for some $\Theta$.
\end{dfn}

Note that uniform Finsler-straightness implies (coarse) uniform regularity. 

Finsler-straight maps carry straight holey lines to Finsler-straight holey lines.

Uniformly regular quasiisometric embeddings $Y\to X$ are uniformly straight \cite{mlem}.

\medskip
We will use the notion of straightness also for {\em extensions} of maps to infinity: 
If $\beta:B\to\Flagt$ is a map defined on a subset $B\subset\geo A$, 
we say that the 
combined 
map 
$$\bar f=f\sqcup\beta:A\sqcup B\to X\sqcup\Flagt $$
is {\em $\Theta$-straight} if 
for every $D$ there exists $d=d(\Theta,D)$ 
such that the induced map on triples 
$$ T(f\sqcup\beta) : T(A,B)\to T(X,\Flagt) $$
sends $D$-straight triples to $(\Theta,d)$-straight triples,
where we use the notation (cf. section \ref{sec:HypSpaces})  
$$T(A,B):=(A\sqcup B)\times A\times(A\sqcup B).$$

The properties 
{\em $\taumod$-straight} and {\em uniformly $\taumod$-straight} are defined accordingly.

\subsubsection{Morse quasiisometric embeddings}

Especially important are straight maps which are quasiisometric embeddings.

\begin{definition}[Morse quasiisometric embedding]\label{defn:Morse-embedding}
A map $Y\to X$ from a Gromov hyperbolic geodesic metric space $Y$ is called a {\em $\taumod$-Morse quasiisometric embedding}
if it sends geodesics to uniform $\taumod$-Morse quasigeodesics. 
\end{definition}

This is equivalent to the definitions given in our earlier papers (see \cite[Def.~7.23]{morse} and \cite[Def.~5.29]{mlem}).
Note that there we allowed more generally for quasigeodesic metric spaces as domains.
However, we showed in \cite[Thm.~6.13]{mlem} that such domains are necessarily Gromov hyperbolic.

Reformulating the above definition,
a quasiisometric embedding $Y\to X$ 
is $\taumod$-Morse if and only if it is uniformly $\taumod$-straight.

The main result of \cite{mlem} implies that 
a quasiisometric embedding $Y\to X$ 
is $\taumod$-Morse if and only if it is uniformly $\taumod$-regular.

\subsubsection{Asymptotics at infinity}
\label{sec:frstrbd}

It is plausible that straightness is related to good asymptotic behavior and the existence of boundary maps.

We first address {\em continuity at infinity}. 
We assume in the following that the map $f:A\to X$ is $\taumod$-regular and $\taumod$-straight.
Our tool is the following direct consequence of Lemma \ref{lem:smccstx}:
\begin{lem}
\label{lem:smccst1}
Suppose that $y_n, y'_n\to\infty$ are divergent sequences contained in $A$ 
so that the triples $(y,y'_n,y_n)$ are $D$-straight for some base point $y\in Y$ and $D\geq0$.

Then the sequences $(f(y_n))$ and $(f(y'_n))$ in $X$ have the same accumulation set in $\Flagt$.
\end{lem}

\begin{dfn}[Shadowing]
A subset $S\subset A$ is {\em shadowing} a subset $\Si\subset\geo A$ at infinity 
(in $A$)
if for every sequence $y_n\to\infty$ in $A$ accumulating at $\Si$ 
there exists
a sequence $y'_n\to\infty$ in $S$
such that the triples $(y,y'_n,y_n)$ are $D$-straight for some base point $y\in Y$ and $D\geq0$.
\end{dfn}

\begin{exa}
(i) A conical accumulation point $\eta\in\geoc A$ is shadowed (in $Y$) by a sequence in $A$ conically converging to it.

(ii) 
Suppose that $A$ is disjoint from a horoball $B\subset Y$ centered at $\zeta\in\geo Y$
and contains a subset $S$ 
which has finite Hausdorff distance from the horosphere $\D B$.
Then $\zeta\in\geo A$ is shadowed (in $A$) by $S$. {Indeed, suppose that $S$ is at Hausdorff distance $D$ from $\D B$. 
Fix $y\in Y$ and consider a sequence $y_n\in A$ converging to $\zeta$. Then for all sufficiently large $n$, the segment 
$yy_n$ will cross the horosphere $\D B$ at a point $z_n$. For $y_n'\in S$ within distance $D$ from $z_n$, the triple $(y, y'_n, y_n)$ will be $D$-straight.}  
\end{exa}

The Lemma \ref{lem:smccst1} immediately yields:
\begin{lem}
\label{lem:smccst}
Suppose that $S\subset A$ is shadowing $\Si\subset\geo A$ in $A$.
Then for every 
subset $W\subset A$ with $\geo W\subset\Si$ it holds that $\geot(f(W))\subseteq\geot (f(S))$. 
\end{lem}
In particular:
\begin{cor}
\label{cor:dffsbshsp}
Suppose that $S,S'\subset A$ have the same visual accumulation set $\Si\subset\geo A$, 
and they are both shadowing $\Si$ in $A$,
then $\geot(f(S))=\geot (f(S'))$. 
\end{cor}

We will apply the last lemma as follows:

\begin{cor}
\label{cor:smccst}
If $S\subset A$ is shadowing the ideal point $\eta\in\geo A$ in $A$
and $\geot (f(S))$ consists of a single simplex $\tau$,
then {the extension of $f$ defined by $f(\eta)=\tau$, defines a map $S\cup \{\eta\}$ continuous at $\eta$}. 
%\mini{is continuously extended to $\eta$ by mapping $\eta\mapsto\tau$.}
\end{cor}

Our second observation concerns {\em antipodality} at infinity.

\begin{lem}
\label{lem:ntpdccst}
Suppose that $y^{\pm}_n\to\infty$ are sequences in $A$ 
whose accumulation sets in $\geo A$ are disjoint. 
Then the accumulation sets of the sequences $(f(y^{\pm}_n))$ 
in $\Flagt$ are antipodal.\footnote{I.e.
every accumulation point of $(f(y^+_n))$ is antipodal to every accumulation point of $(f(y^-_n))$.}
\end{lem}
\proof
We may assume that the sequences $(f(y^+_n))$ $\taumod$-converge, $(f(y^\pm_n))\to\tau_{\pm}$.
Let $y\in A$ be a base point. 
The assumption implies that the triples $(y^-_n,y,y^+_n)$ in $Y$ are $D$-straight for some $D$,
cf.\ Lemma~\ref{lem:bddasyprp}.
By the Finsler-straightness of $f$,
the triples $(f(y^-_n),f(y),f(y^+_n))$ in $X$ are then 
$(\taumod,d)$-straight for some $d$.
This means that there exists a bounded sequence of $\taumod$-parallel sets $P(\tau^-_n,\tau^+_n)$
such that $f(y^{\pm}_n)$ has uniformly bounded distance from $V(f(y),\st(\tau^{\pm}_n))$.
Then $\tau^{\pm}_n\to\tau_{\pm}$ by the definition of flag convergence.
The antipodality of $\tau_{\pm}$ follows from the boundedness of the sequence of parallel sets.
\qed

\medskip 
We next apply these observations to show the existence of a {\em partial boundary map} at infinity.

For a map $\beta:B\to \Flagt$ defined on a subset $B\subseteq \geo A$ 
we say that the combined map
$$
f\sqcup \beta: A \sqcup B\to X\sqcup \Flagt
$$
is {\em continuous at infinity} if for every sequence $(y_n)$ in $A$ with $y_n\to\eta\in B$ it holds that 
$f(y_n)\to \beta(\eta)$
in the sense of flag convergence. 
Note that then $\beta$ must necessarily be continuous. 

We obtain that the map $f: A\to X$ extends continuously to the conical accumulation set:
\begin{prop}
\label{prop:bdmcnc}
There exists an antipodal continuous map $\geoc f:\geoc A\to\geot(f(A))\subset\Flagt$ such that the extended map 
$$ f\sqcup\geoc f:A\sqcup\geoc A \to X\sqcup\Flagt $$
is continuous at infinity.

If $f$ is uniformly $\taumod$-straight, then $\geoc f(\geoc A)\subseteq\geotc(f(A))$. 
\end{prop}
\proof 
Given a point $\eta\in \geoc A$, we pick a sequence $(y_n)$ in $A$ converging to $\eta$ conically. 
After extraction, this sequence moves to infinity ``monotonically''
in the sense that it is $D$-straight for some $D$. \footnote{{Indeed, in view of the conical convergence,  
there exists a constant $D$ and a geodesic ray $\rho: \R_+\to Y$ asymptotic to $\eta$ and a sequence $t_n\in \R_+$ 
diverging to $\infty$ such that 
$d(y_n, \rho(t_n))\le D$. After extraction, we can assume that the sequence $(t_n)$ is increasing. Hence, for $y:=\rho(0)$, each triple $(y, y_m, y_n)$ is $D$-straight. }      
}

Hence, the image sequence $(f(y_n))$ in $X$ is $(\taumod,d)$-straight for some $d$.
(It is also $\taumod$-regular due to our assumption that $f$ is $\taumod$-regular.)
By Corollary~\ref{cor:rgstrhrcnv},
it $\taumod$-converges at infinity,
$$ f(y_n) \to \tau\in\geot(f(A))\subset\Flagt .$$

Since the sequence $(y_n)$ converges to $\eta$ conically,
it is shadowing $\eta$ (in $Y$). 
Corollary~\ref{cor:smccst} therefore implies that 
$f$ is continuously extended to $\eta$ by mapping $\eta\mapsto\tau$.
This shows that 
there exists a well-defined map at infinity $\geoc f:\geoc A\to\geot(f(A))\subset\Flagt$ 
so that the extension $f\sqcup\geoc f$ is continuous at infinity.

The antipodality of $\geoc f$ is a consequence of Lemma~\ref{lem:ntpdccst}.

The last part follows from Lemma~\ref{lem:rgstrhrcnvcncl}.
\qed

\medskip
In particular, we recover (cf. \cite[Theorem 6.14]{mlem}):

\begin{cor}
\label{cor:mrsxtnd}
For every Morse quasiisometric embedding $f:Y\to X$ 
there exists an antipodal continuous map $\geo f:\geo Y\to\geotc(f(Y))\subset\Flagt$ such that the extended map 
$$ f\sqcup\geo f:Y\sqcup\geo Y \to X\sqcup\Flagt $$
is continuous at infinity.
\end{cor}

We now specialize the discussion to a setting motivated by relatively hyperbolic groups. 
Here we can show the existence of {\em full} boundary maps:
\begin{prop}
\label{prop:fllbdmpay}
Suppose that 

(i) $A\subset Y$ 
has finite Hausdorff distance from the complement of a family 
${\mathcal B}=(B_i)_{i\in I}$ of disjoint open horoballs;

(ii)
for some, equivalently, every subset $S_i\subset A$ which has finite Hausdorff distance from a horosphere $\D B_i$,
the $\taumod$-accumulation set $\geot(f(S_i))$ consists of a single simplex $\tau_i$.

Then there exists an antipodal continuous map $\geo f:\geo A\to\geot(f(A))\subset\Flagt$,
sending the center $\zeta_i\in\geo A$ of each horoball $B_i$ to $\tau_i$,
such that the combined map 
$$ f\sqcup\geo f:A\sqcup\geo A \to X\sqcup\Flagt $$
is continuous at infinity.

If $f$ is uniformly $\taumod$-straight, then $\geo f(\geoc A)\subseteq\geotc(f(A))$.
\end{prop}
Regarding condition (ii),
note that according to Corollary~\ref{cor:dffsbshsp} the simplex $\tau_i$ is independent of the choice of $S_i$.

\proof
We continue the argument in the proof of the last proposition. 

In order to further extend the boundary map $\geoc f$ to the non-conical part $(\geo -\geoc)A$ of the accumulation set,
we note that the latter consists of the centers $\zeta_i$ of the horoballs $B_i\in{\mathcal B}$.

Subsets $S_i\subset A$ as in hypothesis (ii) exist by hypothesis (i).
Each subset $S_i$ is shadowing the ideal point $\zeta_i$ (in $A$).
We can therefore apply Corollary~\ref{cor:smccst} again to obtain the desired continuous extension $\geo f$ of $\geoc f$
by mapping $\zeta_i\mapsto\tau_i$ for all $i\in I$.

The antipodality of $\geo f$ follows again from Lemma~\ref{lem:ntpdccst}.
\qed

\section{Asymptotic conditions for subgroups}
\label{sec:BM}

\subsection{Relative asymptotic and boundary embeddedness}
\label{sec:hrbm}

We start with characterizations of Anosov subgroups in terms of their topological dynamics on associated flag manifolds. 
The first such notion given in \cite[Def.~5.12]{anolec} is {\em asymptotic embeddedness}.
The relative version is as follows:\footnote{It will be extended further beyond geometrically finite subgroups in Definition~\ref{def:relembedded}.}

\begin{definition}[Relatively asymptotically embedded]
\label{def:relasembedded}
A subgroup $\Ga<G$ is called {\em relatively $\taumod$-asymptotically embedded} 
if it is $\taumod$-regular, antipodal and 
admits a structure as a relatively hyperbolic group $(\Ga,{\mathcal P})$
such that there exist a $\Ga$-equivariant homeomorphism
$$ \al:\geo\Ga \stackrel{\cong}{\lra} \Lat \subset\Flagt $$
from its ideal boundary to its $\taumod$-limit set.
\end{definition}

\begin{rem}\label{rem:RAE}
{\mini made it a remark} 
This definition can be phrased purely dynamically in terms of the $\Ga$-action on $\Flagt$
by replacing the $\taumod$-regularity with the $\taumod$-convergence condition, see \cite{anolec}.
Note that 
the peripheral structure is uniquely determined {by the action of $\Ga$ on $\Lat$}  
because it can be read off the dynamics on the limit set: 
the peripheral subgroups are the maximal ones with exactly one limit point in $\Lat$.\end{rem}

Since relatively hyperbolic groups act as convergence groups on their ideal boundaries,
so do asymptotically embedded subgroups on their limit sets,
and notions from the theory of abstract convergence groups apply to our setting,
such as conical limit points, bounded parabolic points and bounded parabolic fixed points,
see Definition~\ref{dfn:conlimbddpar}.
As explained in section~\ref{sec:grdfn}, 
every limit point is either conical or bounded parabolic.
The peripheral subgroups $\Pi_i< \Ga$ are precisely the 
stabilizers of the bounded parabolic points $\tau_i$, 
and $\La(\Pi_i)=\{\tau_i\}$. 
For general relatively hyperbolic groups, the $\Pi_i$ can be infinite torsion groups,
however for asymptotically embedded subgroups this cannot occur, because they are linear,
as follows from Schur's theorem or Selberg's lemma. 
Thus all bounded parabolic points are bounded parabolic fixed points. 

\begin{rem}
For  antipodal $\taumod$-regular subgroups $\Ga< G$ with at least two limit points,  
intrinsic conicality (defined in terms of the dynamics of $\Ga\acts \Lat$) 
is equivalent to {\em extrinsic conicality}
(defined in terms of the conical convergence of sequences in $\Ga x\subset X$), see \cite[Proposition 5.41 and Lemma 5.38]{anolec}. 
\end{rem}

For relatively asymptotically embedded subgroups,
the orbit maps extend continuously to infinity by an asymptotic embedding
(which is unique if the limit set has at least three points):
\begin{lemma}
\label{lem:cont}
If $\Ga< G$  is relatively $\taumod$-asymptotically embedded and $x\in X$,
then there is a continuous extension 
$$ \bar o_x = o_x\sqcup\al : \ol\Ga=\Ga\sqcup\geo\Ga\lra X\sqcup \Lat$$
of the orbit map $o_x$ by an asymptotic embedding $\al$.
\end{lemma} 
\proof 
Suppose first  that $|\geo \Ga|\ge 3$. 
Let $(\ga_n)$ be a sequence in $\Ga$ converging to $\xi_+\in \geo \Ga$ in $\ol{\Ga}$. {We claim that 
$(\ga_n)$ flag-converges to $\al(\xi_+)$. Suppose this is not the case. Then, in view of the $\taumod$-regularity of $\Ga$ and the convergence property of 
the action of $\Ga$ on $\geo \Ga$, there exists $\xi_-\in \geo \Ga$ and points $\la_{\pm}\in\Lat, \la_+\ne \al(\xi_+)$, such that, 
after extraction, the sequence $(\ga_n)$ converges, as a sequence of maps, 
to $\xi_+$ uniformly on compacts in $\geo \Ga - \{\xi_-\}$, while $(\ga_n)$ converges to $\la_+$ uniformly on compacts in the open Schubert cell 
$C(\la_-)$. Since $\Lat$ is antipodal, it follows that $(\ga_n)$  flag-converges to $\la_+$ 
uniformly on compacts in $\Lat -\{\la_-\}$. Since the limit set $\Lat$ contains a third point $\tau$ besides 
 $\al(\xi_-)$ and $\la_-$, continuity and equivariance of $\al$ imply that $\ga_n(\tau)\to  \al(\xi_+)$. A contradiction. }

%The above paragraph in RED replaces the old paragraph below: 
%{\mini 
%It follows that there exists a point $\xi_-\in \geo \Ga$ such that 
%We claim that 
%$(\ga_n)$ flag-converges to $\al(\xi_+)$. Assume that this is not the case: 
%In view of the $\taumod$-regularity of $\Ga$, 
%there exist $\la_{\pm}\in\Lat$
%such that, after extraction, 
%$(\ga_n)$ converges to $\la_+$ uniformly on compacts in the open Schubert cell 
%$C(\la_-)$, where $\la_+\ne\al(\xi_+)$. 
%Since $\Lat$ is antipodal, it follows that $(\ga_n)$  flag-converges to $\la_+$ 
%uniformly on compacts in $\Lat -\{\la_-\}$. 
%Since the limit set contains a third point beyond $\al(\xi_-)$ and $\la_-$,
%we conclude that $\tau=\la_+$,
%a contradiction. }

If $|\geo \Ga|\le 1$, there is nothing to prove. 
If $|\geo\Ga|=2$ then $\P=\emptyset$, $\Ga$ is virtually cyclic (see Remark~\ref{rem:rhdf}) and the claim follows from 
\cite[Lemma 5.38]{anolec}.  
\qed

\medskip
The following related condition is weaker than relative asymptotic embeddedness, 
but easier to verify, since there is no need to check regularity and to identify the limit set: 
\begin{definition}[Relatively boundary embedded]
\label{def:relboundaryembedded}
A discrete subgroup $\Ga<G$ is called {\em relatively $\taumod$-boundary embedded} 
if it admits a structure as a relatively hyperbolic group $(\Ga,{\mathcal P})$
such that there exist an antipodal $\Ga$-equivariant embedding
$$ \beta:\geo\Ga \lra\Flagt $$
called a {\em boundary embedding}.
\end{definition}

\medskip 
In the {\em Zariski dense} case, relative boundary embeddedness implies relative asymptotic boundary embeddedness 
(see \cite[Corollary 5.14]{anolec} for the absolute case): 

\begin{thm}
\label{thm:BdZarden}
If $\Ga< G$ is relatively $\taumod$-boundary embedded and Zariski dense,
then it is relatively $\taumod$-asymptotically embedded. 
\end{thm} 
\proof By Zariski density, $\geo \Ga$ is infinite and hence the action $\Ga\acts\geo\Ga$ is minimal.
Applying Theorem~\ref{thm:Zardense} and Addendum~\ref{add:Zarden} to the given boundary embedding,
it follows that $\Ga$ is $\taumod$-regular and $\Lat=\beta(\geo\Ga)$.
\qed

\subsection{A higher rank Beardon-Maskit condition}

The actions of relatively hyperbolic groups on their ideal boundaries are convergence actions characterized by the Beardon-Maskit property
\cite{Yaman}.
We use this characterization to translate relative asymptotic embeddedness into a higher rank Beardon-Maskit condition
for the action on the limit set.
We formulate this condition for antipodal regular subgroups because their actions on the limit set are convergence:

\begin{definition}[Relatively RCA]
\label{dfn:rlrca}
An antipodal $\taumod$-regular subgroup $\Ga< G$ is called {\em relatively $\taumod$-RCA}
if each $\taumod$-limit point is either a conical limit point or a bounded parabolic point (for the action $\Ga\acts\Lat$)
and, moreover, the stabilizers of the bounded parabolic points are finitely generated. 
\end{definition}

We recall, see section~\ref{sec:conv},
that for the stabilizer $\Ga_{\tau}<\Ga$ of a bounded parabolic  point $\tau\in\LatGa$ 
it holds that $\Lat(\Ga_\tau)=\{\tau\}$. 

{
\begin{rem}
The terminology RCA was first introduced in \cite{morse} (and, in published form, in  \cite{anosov}): This abbreviation stands for {\em regular, conical, antipodal}: Such subgroups are required to be $\taumod$-regular, antipodal (i.e. 
any two distinct points of $\Lat$ are antipodal) and conical, i.e. every limit point is conical.  The notion of {\em relatively RCA groups} 
is a relativization of the RCA condition: The difference is that we are now allowing certain non-conical limit points.  
\end{rem}
}

All finitely generated regular subgroups with one point limit set are relatively RCA.
In general, we know little about the structure of such subgroups. 
On the other hand, for {\em uniformly} regular subgroups with one point limit set 
we have the following information on their geometric and algebraic properties,
which are well-known in rank one: 

\begin{lem}\label{lem:ur->fgn}
A uniformly $\taumod$-regular subgroup of $G$ with one point $\taumod$-limit set
consists of elements with zero translation number, 
is finitely generated and virtually nilpotent.  
\end{lem}
\proof This follows from Lemma \ref{lem:cclsbgmgp} and Corollary \ref{cor:nilfg}.  
\qed

\medskip
Clearly, relative asymptotic embeddedness implies relatively RCA.
The converse is a consequence of Yaman's theorem.
Thus:

\begin{thm}\label{thm:RRCA=RAE}
A subgroup $\Ga<G$ is relatively $\taumod$-RCA if and only if it is relatively $\taumod$-asymptotically embedded.
\end{thm}

In fact,
Yaman's theorem applies only if $|\Lat|\geq3$.
If $|\Lat|=2$, then $\P=\emptyset$ and $\Ga$ is (absolutely) $\taumod$-asymptotically embedded
by \cite[Lemma 5.38]{anolec}.

\subsection{More general relative settings}
\label{sec:mgrlstt}

Let $\Ga<G$ be a discrete subgroup.
We equip $\Ga$ as an abstract discrete group 
with an additional {\em intrinsic} geometric structure 
in the form of a properly discontinuous isometric action $\Ga\acts Y$ on a Gromov hyperbolic proper geodesic space $Y$.
This action is {\em not} required to be cocompact.
(If it is cocompact or, more generally, undistorted,
then $\Ga$ is word hyperbolic 
and the additional intrinsic structure amounts to the choice of a word metric.
This is the context in which we mostly worked in our earlier papers.)
We are interested in 
geometric and dynamical properties of the action $\Ga\acts X$ 
{\em relative} to the action $\Ga\acts Y$.

To relate the actions $\Ga\acts X$ and $\Ga\acts Y$,
we fix base points $x\in X$ and $y\in Y$ so that $\Ga_y\leq\Ga_x$ 
and consider the $\Ga$-equivariant map of orbits 
$$ o_{x,y}: \Ga y \to \Ga x , \quad \ga y \mapsto \ga x .$$
Note that for any point $y\in Y$ there exists a point $x\in X$ fixed by $\Ga_y$, 
because $\Ga_y$ is finite and finite groups acting isometrically on symmetric spaces have fixed points. 

We further extend the relative versions of asymptotic and boundary embeddedness 
(see Definitions~\ref{def:relasembedded} and~\ref{def:relboundaryembedded})
to our present more general setting. 
\begin{dfn}[Relatively boundary and asymptotically embedded II]
\label{def:relembedded}
A discrete subgroup $\Ga<G$ is called 

(i) {\em $\taumod$-boundary embedded relative $\Ga\acts Y$}
if there exists a $\Ga$-equivariant antipodal embedding,
called a {\em boundary embedding},
$$\beta:\LaY\to\Flagt .$$

(ii) {\em $\taumod$-asymptotically embedded relative $\Ga\acts Y$}
if it is $\taumod$-regular, antipodal
and there exists a $\Ga$-equivariant homeomorphism,
called an {\em asymptotic embedding},
$$\al:\LaY\stackrel{\cong}{\to}\LaXt \subset\Flagt .$$ 
\end{dfn}

As before,
in the non-degenerate case when $|\LaY|\geq3$,
an asymptotic embedding continuously extends the maps of orbits to infinity.
Lemma~\ref{lem:cont} and its proof directly generalize:
\begin{lem}
\label{lem:asmbdcntxt}
If $\Ga<G$ is $\taumod$-asymptotically embedded
relative $\Ga\acts Y$ 
and if $|\LaY|\geq3$, 
then the combined map 
$$ \bar o_{x,y}=o_{x,y}\sqcup\al : \ol{\Ga y}=\Ga y\sqcup\LaY \to \ol{\Ga x}^{\,\taumod}= \Ga x\sqcup\LaXt \subset X\sqcup\Flagt $$
is continuous.
\end{lem}

%The continuity on the orbit $\Ga y$ is trivial by discreteness. 
{Similarly, the proof of Theorem \ref{thm:BdZarden} goes through as well and we obtain:}

\begin{thm}
\label{thm:BdZarden2}
{If $\Ga< G$ is relatively $\taumod$-boundary embedded and Zariski dense,
then it is relatively $\taumod$-asymptotically embedded. } 
\end{thm}

\section{Coarse geometric conditions for subgroups}

We introduce two coarse geometric conditions 
which are a priori stronger than the asymptotic conditions discussed above.
The advantage of these coarse geometric properties is that they allow for a {\em local-to-global principle} 
similar to the one for Morse subgroups (cf. \cite[\S 7]{morse}) 
and hence {\em both} define classes of discrete subgroups which are {\em all structurally stable} ({these stability results will be proven elsewhere}). These conditions are also sometimes easier to verify in concrete situations. 
The main results in this section compare the coarse geometric conditions to asymptotic embeddedness
(see Theorems~\ref{thm:RM->UR}, \ref{thm:RM=GF}, \ref{thm:str-ae} and \ref{thm:ae->str}).

\subsection{Relatively Morse subgroups}
\label{sec:relmo}

In our earlier paper \cite{morse} 
we defined Morse subgroups as finitely generated word-hyperbolic subgroups 
whose orbit maps are Morse quasiisometric embeddings.
We relativize this as follows:

\begin{definition}[Relatively Morse]
\label{dfn:rlmrs}
A subgroup $\Ga<G$ is called {\em relatively $\taumod$-Morse}
if there exists a relatively hyperbolic structure ${\mathcal P}$ on $\Ga$ with a Gromov model $Y$
and a $\Ga$-equivariant $\taumod$-Morse quasiisometric embedding $f:Y\to X$.
\end{definition} 

The peripheral subgroups have to be virtually nilpotent, as follows from Corollary \ref{cor:DY}.
We will see that the peripheral structure ${\mathcal P}$ is uniquely determined
because relatively Morse implies relatively asymptotically embedded (see Corollary~\ref{cor:uniqueness}).

{\mini Relatively Morse subgroups are uniformly regular, since Morse quasiisometric embeddings are.
The latter continuously extend to infinity (see Corollary~\ref{cor:mrsxtnd}).}

In the equivariant situation (for general discrete subgroups which need not be relatively Morse)
one can relate the limit sets:
\begin{lem}\label{lem:RM-RAE}
Let $\Ga<G$ be a discrete subgroup. 
Suppose that $\Ga\acts Y$ is a properly discontinuous isometric action on a proper geodesic hyperbolic space $Y$ 
and that $f: Y\to X$ is a $\Ga$-equivariant $\taumod$-Morse quasiisometric embedding.

Then $\Ga$ is $\taumod$-uniformly regular 
and $\LaXt=\geo f(\La_Y)$ is antipodal.
\end{lem}
\proof
This follows from the uniform regularity of Morse quasiisometric embeddings ({the images of geodesic segments are uniformly $\taumod$-regular quasigeodesics in $X$, which is the content of Definition \ref{defn:Morse-embedding}}), 
the continuity of $f\sqcup\geo f$ at $\geo Y$ and the antipodality of $\geo f$ ({the last two properties are established in Corollary \ref{cor:mrsxtnd}; see also \cite[Theorem 6.14]{mlem}}).
\qed

\medskip
In the relatively Morse setting, we obtain:
\begin{thm}
\label{thm:RM->UR}
Every relatively $\taumod$-Morse subgroup $\Ga< G$ is relatively $\taumod$-asymptotically embedded. 
If $f: Y\to X$ is an equivariant $\taumod$-Morse quasiisometric embedding as in the definition of relatively Morse subgroups,
then $\geo f$ is an asymptotic embedding.
\end{thm}
\proof 
In this situation 
$\La_Y=\geo Y$ and the lemma yields the assertion. 
\qed

\begin{cor}
\label{cor:uniqueness} 
The relatively hyperbolic structure on a relatively Morse subgroup is unique. 
\end{cor}
\proof This follows from the uniqueness of the relatively hyperbolic structure on relatively asymptotically embedded subgroups, {see Remark \ref{rem:RAE}}. \qed

\begin{thm}
\label{thm:RM=GF}
If $X$ has rank one, then relatively Morse is equivalent to geometrically finite.
\end{thm}
\proof 
According to Theorem~\ref{thm:RM->UR}, 
a relatively Morse subgroup is relatively asymptotically embedded 
and hence its action on its limit set satisfies the Beardon-Maskit condition.
In rank one, this is classically known to be equivalent to geometric finiteness (see \cite{Bowditch}). 

Conversely, let $\Ga<G$ be geometrically finite. 

If $|\La|\geq2$,
the closed convex hull of the limit set serves as a Gromov model, $Y:= \CH(\La)$.
The subgroups $\Pi_i$ are the stabilizers of the bounded parabolic fixed points of $\Ga$. 
There exists a $\Ga$-invariant family of pairwise disjoint horoballs $B_i$ 
in $X$ such that $\Ga$ acts cocompactly on 
$$
Y^{th}= Y- \bigcup_{i} B_i, 
$$
and we use this as a Gromov model of $(\Ga,\P)$. 
The Morse quasiisometric embedding is the inclusion $Y\embed X$,
and we see that $\Ga< G$ is a relatively Morse subgroup. 

If $\La$ consists of a single ideal point $\la$,
we equip it with the trivial relatively hyperbolic structure $\P=\{\Ga\}$. 
Let $B\subset X$ be a horoball centered at $\la$. 
It is preserved by $\Ga$, 
and according to \cite[sect. 4]{Bowditch}
there exists a $\Ga$-invariant closed convex subset $C\subseteq X$
such that the action $\Ga\acts \D B\cap C$ is cocompact.
($C$ can be obtained as the closed convex hull of a $\Ga$-orbit in $\geo X-\{\la\}$.)
One can then take $B\cap C$ as a Gromov model.
\qed 

\begin{cor}
\label{cor:rlmqvsmbemb}
In rank one,
relatively Morse is equivalent to relatively asymptotically embedded.
\end{cor}

\begin{example}\label{ex:geodembedding}
Let $\rank(X)\geq2$ 
and let $X_1\subset X$ be a totally geodesic subspace of rank 1. 
It is $\taumod$-regular for $\iota$-invariant face-types $\taumod\subset \taumod(X_1)$, 
see Example \ref{ex:rank1embedding}. Let $\Ga<G$ be a subgroup which preserves $X_1$ 
and acts on it as a geometrically finite group.
Then $\Ga$ is $\taumod$-Morse,
the Gromov model being a convex subset of $X_1$ as described in the proof of Theorem~\ref{thm:RM=GF}.
\end{example}

\subsection{Relatively Finsler-straight subgroups}
\label{sec:str=ae}

This section is the heart of the paper. 
We define the notion of relatively Finsler-straight groups of isometries of symmetric spaces.
For actions of relatively hyperbolic groups, we prove its equivalence to relative asymptotic embeddedness.

\subsubsection{From straightness to boundary maps}
\label{sec:fstrctns}

We take up the discussion of Finsler-straight maps in an equivariant situation.
We deduce from the results in section~\ref{sec:frstrbd}
that Finsler-straightness of maps of orbits 
implies the existence of partial and, under suitable assumptions,
of full boundary maps.

We work in the general relative setting of section~\ref{sec:mgrlstt}.
The notion of Finsler-straightness for maps (see Definition~\ref{dfn:fstrghtmp}) carries over to subgroups:
\begin{definition}[Finsler-straight subgroup]
\label{dfn:fstrghtctn}
A discrete subgroup $\Ga<G$ is said to be 

(i) {\em $\taumod$-straight} rel $\Ga\acts Y$ 
if the map $o_{x,y}$ is $\taumod$-straight. 

(ii) {\em uniformly $\taumod$-straight} rel $\Ga\acts Y$ if $o_{x,y}$ is uniformly $\taumod$-straight. 
\end{definition}

Note that a $\taumod$-straight subgroup $\Ga<G$ is uniformly $\taumod$-straight if and only if it is uniformly $\taumod$-regular.

We will consider the notion of relative straightness only in the context of regular subgroups.

\medskip
Now we use the results from section~\ref{sec:frstrbd} 
in order to obtain boundary maps for Finsler-straight subgroups.
Proposition~\ref{prop:bdmcnc}, applied to the maps of orbits $o_{x,y}$, yields a {\em partial} asymptotic embedding:
\begin{cor}
\label{cor:bdmcnc}
If $\Ga<G$ is $\taumod$-straight relative $\Ga\acts Y$,
then there exists an antipodal map $\geoc o_{x,y}:\LaYc\to\LaXt\subset\Flagt$ such that the extended map 
$$ o_{x,y}\sqcup\geoc o_{x,y}:\Ga y\sqcup\LaYc \to \Ga x\sqcup\LaXt \subset X\sqcup\Flagt $$
is continuous (at infinity).

If $\Ga<G$ is uniformly $\taumod$-straight relative $\Ga\acts Y$,
then $\geoc o_{x,y}(\LaYc)\subseteq\LaXtc$.
\end{cor}

Note that 
the boundary map $\D^{con}_\infty o_{x,y}$ is independent of the choice of the base points $y,x$.

Note also that, if $|\LaY|\geq2$, 
then $\LaYc$ is nonempty \cite[Thm.\ 2R]{Tukia1994}
and hence dense in $\LaY$. 

\medskip
In the relatively hyperbolic setting,
i.e. when $\Ga\acts Y$ is the action on a Gromov model, 
we obtain a {\em full} asymptotic embedding
under an additional assumption on the actions of the peripheral subgroups.
Namely, we consider the following condition:
\begin{definition}[Tied-up horospheres]
\label{def:tphrsphs}
A $\taumod$-regular subgroup $\Ga<G$ is said to have {\em tied-up horospheres}
with respect to a relatively hyperbolic structure ${\mathcal P}$ 
if the limit set $\LaXt(\Pi_i)\subset\Flagt$ of each peripheral subgroup $\Pi_i< \Ga$ is a singleton. 
\end{definition}

We adapt Finsler-straightness as follows to relatively hyperbolic subgroups:
\begin{definition}[Relatively Finsler-straight]
A $\taumod$-regular subgroup $\Ga<G$ is called 

(i) {\em relatively $\taumod$-straight}
if there exists a relatively hyperbolic structure ${\mathcal P}$ on $\Ga$ with a Gromov model $Y$
such that $\Ga$ is {\em $\taumod$-straight} relative $\Ga\acts Y$
and has tied-up horospheres.

(ii) {\em relatively uniformly $\taumod$-straight} if in addition $\Ga$ is uniformly $\taumod$-regular.
\end{definition} 

In view of Lemma~\ref{lem:instr},
relative Finsler-straightness does not depend on the choice of the Gromov model $Y$.

Applying Proposition~\ref{prop:fllbdmpay} now yields: 
\begin{thm}
\label{thm:str-ae}
If $\Ga<G$ is relatively $\taumod$-straight,
then it is relatively $\taumod$-asymptotically embedded. 
\end{thm}
\proof
We apply Proposition~\ref{prop:fllbdmpay} with $A=\Ga y$ and $f=o_{x,y}$.
Hypothesis (i) of the proposition is satisfied, 
because $\Ga$ acts cocompactly on the thick part $Y^{th}\subset Y$ of the Gromov model,
which equals the complement of the family of peripheral horoballs $B_i$.
In hypothesis (ii) we can take $S_i=\Pi_i y$,
because $\Pi_i$ acts cocompactly on horospheres at $\zeta_i$ (cf.\ Lemma~\ref{lem:prphsbgcchs}),
and the condition is satisfied because $\Ga$ has tied-up horospheres.
Since $\geo A=\geo(\Ga y)=\geo Y$ and $\geot(f(A))=\geot(\Ga x)=\LaXt$, 
the proposition yields an antipodal continuous map $\geo o_{x,y}:\geo Y\to\LaXt$ sending $\zeta_i\mapsto\tau_i$
so that the extension $\bar o_{x,y}=o_{x,y}\sqcup\geo o_{x,y}$ is continuous at infinity.
The latter implies that the image of $\geo o_{x,y}$ equals $\LaXt$.
\qed

\medskip
Note that, as a consequence of the theorem, 
the relatively hyperbolic structure in the definition of relative Finsler-straightness is {\em unique}.

\subsubsection{From boundary maps to straightness}
\label{sec:main}

We now, conversely, 
explore what the existence of boundary maps implies for the orbit geometry of actions $\Ga\acts X$.
We show that asymptotic embeddedness implies Finsler-straightness. 
Our discussion follows \cite[\S 5.3]{anolec}, generalizing it.

Suppose first that $\beta:\LaY\to\Flagt$ is a {\em boundary embedding} relative $\Ga\acts Y$.

Our preliminary step concerns the position of the $\Ga$-orbits in $X$ 
relative to the parallel sets spanned by pairs of simplices in the image of $\beta$
(cf. \cite[Lemma 5.3]{anolec}):
\begin{lem}
\label{lem:contrdstprst}
For every $D$ there exists $d$ (also depending on the $\Ga$-actions and -orbits) so that: 

If a triple $(\eta_-,\ga y,\eta_+)$ with $\eta_{\pm}\in\La_Y$ is $D$-straight, 
then the triple $(\beta(\eta_-),\ga x,\beta(\eta_+))$ is 
$(\taumod,d)$-straight.
\end{lem}
\proof 
By equivariance, we may assume that $\ga=e$.

The set of 
pairs $(\eta_-,\eta_+)\in (\geo Y\times \geo Y)-\De_{\geo Y}$,
for which the triple $(\eta_-,y,\eta_+)$ is $D$-straight, is compact.
It follows that the set $C$ of their images $(\beta(\eta_-),\beta(\eta_+))$ 
in the space $(\Flagt\times\Flagt)^{opp}\subset\Flagt\times\Flagt$ is also compact.
Since $(\Flagt\times\Flagt)^{opp}$ is a homogeneous $G$-space,
it is of the form $C=C'\cdot(\tau_0^-,\tau_0^+)$ with a compact subset $C'\subset G$ and some antipodal pair $(\tau_0^-,\tau_0^+)$. 
The set of triples 
$(\beta(\eta_-),x,\beta(\eta_+))=(g\tau_0^-,x,g\tau_0^+)=g(\tau_0^-,g^{-1}x,\tau_0^+)$ for $g\in C'$
is $(\taumod,d)$-straight for some $d=d(D)$,
because the set ${C'}^{-1}x\subset X$ is compact and depends on $D$.
\qed

\medskip
The conclusion can be rephrased as follows:
If $\ga y$ lies within distance $D$ of a line $\eta_-\eta_+\subset Y$, $\eta_{\pm}\in\LaY$,
then $\ga x$ lies within distance $d(D)$ of the parallel set $P(\beta(\eta_-),\beta(\eta_+))\subset X$.

\medskip
In order to get more control,
we strengthen our assumptions for the rest of this section:
\begin{ass}
\label{ass:rlsmbd}
$\Ga<G$ is {\em $\taumod$-asymptotically embedded} relative $\Ga\acts Y$ 
with asymptotic embedding $\al:\LaY\stackrel{\cong}{\to}\LaXt$ 
and $|\LaY|\geq3$.
\end{ass}

Then $\Ga<G$ is $\taumod$-regular
and $\al$ continuously extends the map of orbits $o_{x,y}$ (Lemma~\ref{lem:asmbdcntxt}).

Now we can relate the position of the $\Ga$-orbits in $X$ to Weyl cones:
\begin{lemma}\label{lem:cntrldstwcn}
For every $D$ there exists $d$ such that: 

If a triple $(\ga_-y,\ga y,\eta_+)$, $\eta_+\in\La_Y$, is $D$-straight, 
then $(\ga_-x,\ga x,\al(\eta_+))$ is $(\taumod,d)$-straight.
\end{lemma}
\proof Let $$R:=d(y, \QCH(\LaY)).$$
We may assume that $\ga_-=e$, and denote $\eta_+=:\eta$.

Suppose that the triple $(y,\ga y,\eta)$ with $\eta\in\La_Y$ is $D$-straight.
Due to the quasiconvexity of $\QCH(\LaY)$,
the ray $y\eta$ lies within distance $R+C\de$
of a geodesic line $\hat\eta\eta\subset Y$ with $\hat\eta\in\La_Y$,
and it follows that the triple $(\hat\eta,\ga y,\eta)$ is $(D+R+C\delta)$-straight. 

By Lemma~\ref{lem:contrdstprst},
the triple $(\al(\hat\eta),\ga x,\al(\eta))$ is $d$-straight for some $d=d(D+R+C\delta)$.
This means that $\ga x$ lies within distance $d$ of the parallel set $P(\al(\hat\eta),\al(\eta))$.
It applies in particular to $x=ex$. 
It follows (compare \cite[Dichotomy Lemma 5.5 and Proposition 5.16]{anolec}) 
that $\ga x$ lies within distance $d$ of a Weyl cone $V(\bar x,\st(\tau'))$
for a point $\bar x\in P(\al(\hat\eta),\al(\eta))$ with $d(x,\bar x)\leq d$
and a type $\taumod$ simplex $\tau'\subset\geo P(\al(\hat\eta),\al(\eta))$,
that is,
the triple $(x,\ga x,\tau')$ is $(\taumod,d)$-straight.
It is important to note that $\al(\eta)$ is the only simplex contained in $\geo P(\al(\hat\eta),\al(\eta))$ which is antipodal to $\al(\hat\eta)$.
Hence either $\tau'=\al(\eta)$ or $\tau'$ is {\em not} antipodal to $\al(\hat\eta)$.

Consider a sequence of $D$-straight triples $(y,\ga_n y,\eta_n)$ with $\eta_n\in\La_Y$
and $\ga_n\to\infty$ in $\Ga$,
and corresponding sequences of ideal points $\hat\eta_n\in\LaY$, points $\bar x_n\in P(\al(\hat\eta_n),\al(\eta_n))$ 
and simplices $\tau'_n\subset\geo P(\al(\hat\eta_n),\al(\eta_n))$ 
as above.
Then $\ga_n x$ lies within distance $d$ of the Weyl cone $V(\bar x_n,\st(\tau'_n))$.

Suppose that the simplices $\tau'_n$ are not antipodal to the simplices $\al(\hat\eta_n)$ for all $n$.
After extraction, we may assume that $\eta_n\to\eta$, $\hat\eta_n\to\hat\eta$, $\bar x_n\to\bar x$ and $\tau'_n\to\tau'$.
Then $\al(\eta_n)\to\al(\eta)$, $\al(\hat\eta_n)\to\al(\hat\eta)$, $\ga_n x\to\tau'$
and $\tau'\subset\geo P(\al(\hat\eta),\al(\eta))$.
Moreover, since the relation of being non-antipodal is closed with respect to the visual topology,
$\tau'$ is not antipodal to $\al(\hat\eta)$,
and hence $\tau'\neq\al(\eta)$.
On the other hand, $\ga_ny\to\eta$ and hence $\ga_nx\to\al(\eta)$
due to the continuity of $\bar o_{x,y}$, 
a contradiction. 

It follows that, for all $\ga\in\Ga$ outside a finite subset of $\Ga$ depending on $D$,
a triple $(x,\ga x,\al(\eta))$ is $(\taumod,d)$-straight
whenever the triple $(y,\ga y,\eta)$ with $\eta\in\LaY$ is $D$-straight.
After suitably enlarging $d$,
the implication holds for all $\ga\in\Ga$.
\qed 

\medskip
We rephrase the conclusion:
If $\ga y$ lies within controlled distance of a ray $(\ga_-y)\eta_+\subset Y$, $\eta_+\in\LaY$,
then $\ga x$ lies within controlled distance of the Weyl cone $V(\ga_-x,\al(\eta_+))\subset X$.

In the next step,
we control the straightness of triples in $\Ga$-orbits:
\begin{lem}
\label{lem:rbdstrghtn}
For every $D$ there exists $d$ such that: 

If a triple $(\ga_-y,\ga y,\ga_+y)$ is $D$-straight, 
then $(\ga_-x,\ga x,\ga_+x)$ is $(\taumod,d)$-straight.
\end{lem}
\proof 
Suppose that the triple $(\ga_-y,\ga y,\ga_+y)$ is $D$-straight. 
By the quasiconvexity of $\QCH(\LaY)$,
it is $C\cdot (D+R+\de)$-Hausdorff close 
to a triple of points lying in the same order on a
geodesic line $\eta_-\eta_+\subset Y$ with $\eta_\pm\in\LaY$. 
Hence its middle point $\ga y$ lies within distance $C'\cdot (D+R+\de)$ of geodesic rays $(\ga_{\mp}y)\eta_{\pm}$.

By Lemmas~\ref{lem:contrdstprst} and~\ref{lem:cntrldstwcn}, 
the points $\ga_{\pm}x,\ga x$ lie within 
distance $d$ 
from the parallel set $P=P(\al(\eta_-), \al(\eta_+))$,
and 
$\ga x$ lies within 
distance $d$ 
from the two Weyl cones $V_{\pm}=V(\ga_{\mp}x,\st(\al(\eta_{\pm})))$ for some $d=d(D)$.
(We suppress the dependence on the actions and orbits.)

Let $\bar x_{\pm},\bar x\in P$ denote the nearest-point projections of $\ga_{\pm} x,\ga x$.
Then the Weyl cones $V_{\pm}$
and $\ol V_{\pm}=V(\bar x_{\mp},\st(\al(\eta_{\pm})))\subset P$
have Hausdorff distance $\leq d(\ga_{\mp}x,\bar x_{\mp})\leq d$. 
Hence $\ga x$ lies within distance $2d$ 
from both Weyl cones 
$\ol V_{\pm}$, 
and $\bar x$ lies within distance $3d$ from them. 

Now we invoke again the $\taumod$-regularity of $\Ga$.
It implies the existence of a finite subset $F\subset\Ga$ 
depending on $D$ such that:
If $\ga^{-1}\ga_{\mp}\not\in F$,
then $\bar x\in\ol V_{\pm}$.
If both conditions are satisfied, then 
$\bar x$ lies on a $\taumod$-Finsler geodesic $\bar x_-\bar x_+$, 
and hence 
the triple $(\ga_- x, \ga x, \ga_+ x)$ is $(\taumod,d)$-straight.

If one of the elements $\ga^{-1}\ga_{\mp}$ lies in $F$,
then the corresponding distance $d(\ga_{\mp}x,\ga x)$ is bounded 
and the conclusion holds trivially after increasing $d$ sufficiently.
\qed

\medskip
Lemmas~\ref{lem:contrdstprst}, \ref{lem:cntrldstwcn} and~\ref{lem:rbdstrghtn} yield together:
\begin{prop}
\label{prop:rbdstrght}
The extension $\ol o_{x,y}=o_{x,y}\sqcup\al$ of the map of orbits $o_{x,y}$ is $\taumod$-straight.
\end{prop}

Specializing to the relatively hyperbolic setting,
we obtain the converse to Theorem~\ref{thm:str-ae}:

\begin{thm}
\label{thm:rlsmbrlstr}
If $\Ga<G$ is relatively $\taumod$-asymptotically embedded, then it is relatively $\taumod$-straight.
\end{thm}
\proof
In the case when $|\geo Y|\geq3$,
i.e. when Assumption~\ref{ass:rlsmbd} is satisfied,
the straightness of $\Ga$ relative to the action on the Gromov model is the content of the proposition.
That $\Ga$ has tied-up horospheres is immediate from asymptotic embeddedness.

If $|\geo Y|=2$, then ${\mathcal P}=\emptyset$ and we are in the absolute case.
There, asymptotic embeddedness implies Morse, and this in turn straightness \cite{anolec}.

If $|\geo Y|=1$, then ${\mathcal P}=\{\Ga\}$ 
and relative straightness amounts to the $\taumod$-regularity of $\Ga$ and $|\LaXt|=1$.
Both properties follow from asymptotic embeddedness.
\qed

\subsubsection{Further to uniform straightness}\label{sec:mor-uniform}

We return to the more general setting of a subgroup $\Ga<G$ which is asymptotically embedded relative to an action $\Ga\acts Y$
as in Assumption~\ref{ass:rlsmbd}.
We show that under suitable assumptions the subgroup $\Ga$ is {\em uniformly} regular.

We continue the discussion of the previous section 
and further promote the control on the position of triples in orbits of $\Ga\acts X$
to control on the position of {\em holey lines}.
By Proposition \ref{prop:rbdstrght}, 
we know that straight holey lines $q:H\to \Ga y$ go to 
Finsler-straight holey lines $o_{x,y}\circ q:H\to\Ga x$.
We establish next that the latter lie near Finsler geodesics in $X$. 
To achieve this,
we use that the holey lines $o_{x,y}\circ q$ are asymptotically embedded 
and satisfy, as parts of orbits of the regular action $\Ga\acts X$, a weak form of {\em uniform} regularity.
This allows us to apply  Addendum~\ref{add:strghtphrgclfns}, 
and we obtain:

\begin{prop}
\label{prop:hlflfns}
For $D$ there exists $d$ such that: 

If $q: H\to\Ga y$ is a $D$-straight holey line,
then there exists a $\taumod$-Finsler line $\tau_-\tau_+$, $\tau_{\pm}\in\LaXt$,
and a monotonic map $q':H\to \tau_-\tau_+$
which is $d$-close to the holey line $o_{x,y}\circ q:H\to\Ga x$. 
\end{prop}
\proof
By straightness,
$q(H)$ lies within distance $D'=D'(D)$ of a line in $Y$,
and within distance $D''=D''(D,R)$ of a line $\eta_-\eta_+\subset Y$ asymptotic to $\La_Y$, $\eta_{\pm}\in\LaY$.
The extended holey line $\bar q:\ol H=H\sqcup\{\pm\infty\}\to\ol Y=\Ga y\sqcup\LaY$
with $\bar q(\pm\infty)=\eta_{\pm}$ is $D''$-straight.

Since $\bar o_{x,y}=o_{x,y}\sqcup\al$ is 
$\taumod$-straight 
by Proposition~\ref{prop:rbdstrght}, 
it follows that the extended holey line 
$\ol{o_{x,y}\circ q}=\bar o_{x,y}\circ\bar q:\ol H\to\Ga x\sqcup\LaXt$ 
mapping $\pm\infty\mapsto\al(\eta_{\pm})=:\tau_{\pm}$
is $(\taumod,d)$-straight for some $d=d(D)$.
Furthermore, the holey lines $o_{x,y}\circ q$ are weakly uniformly $\taumod$-regular 
as a consequence of the $\taumod$-regularity of the action $\Ga\acts X$.
Applying Addendum~\ref{add:strghtphrgclfns} yields the assertion.
\qed

\medskip
Since the image holey lines $o_{x,y}\circ q$ in $X$ follow Finsler geodesics,
their weak uniform regularity turns into (strong) uniform regularity 
when there are no arbitrarily large holes,
i.e.\ when the holey lines $o_{x,y}\circ q$, equivalently, the straight holey lines $q$ are coarsely connected:
\begin{claim}
\label{claim:crscnhlnfrg}
For $D,L$ there exist $\Theta,d$ such that:

If $q:H\to\Ga y$ is a $D$-straight holey line
which is coarsely $L$-connected,
then the $\taumod$-Finsler line $\tau_-\tau_+$ in Proposition \ref{prop:hlflfns}
can be chosen to be $\Theta$-regular.
\end{claim}
\proof
The holey line $o_{x,y}\circ q:H\to\Ga x$ is then coarsely $l$-connected with $l=l(L)$
and the assertion is a consequence of Claim \ref{claim:lngsgmnfrg}.
\qed

\medskip
Straight holey lines in $\Ga y$ with holes of bounded size are, up to reparameterization, {\em quasigeodesics}.
We conclude that $o_{x,y}$ sends uniform quasigeodesics in $\Ga y$ 
to uniformly $\taumod$-regular uniform quasigeodesics in $\Ga x$:
\begin{thm}
Suppose that $\Ga<G$ is $\taumod$-asymptotically embedded rel $\Ga\acts Y$. 
Then for $L,A$ there exist $l,a,\Theta,d$ such that:

If $q:I\to\Ga y\subset Y$ is an $(L,A)$-quasigeodesic,
then $o_{x,y}\circ q:I\to\Ga x\subset X$ is an $(l,a)$-quasigeodesic 
which is contained in the $d$-neighborhood of a $\Theta$-Finsler geodesic.
\end{thm}

\begin{rem}
In the ``absolute'' case, that is, when $\Ga\acts Y$ is cocompact (or undistorted) and hence $\Ga$ is Gromov hyperbolic,
this recovers our earlier result 
that $\taumod$-asymptotically embedded subgroups $\Ga<G$ are $\taumod$-Morse
\cite{anolec}.
\end{rem}

The theorem yields some {\em partial} uniform regularity for the map of orbits $o_{x,y}$.
In order to obtain full uniform regularity,
we impose an additional assumption, cf.\ Definition~\ref{dfn:nfrglhls}.
We then can extend Claim \ref{claim:crscnhlnfrg} as follows:

\begin{claim}
\label{claim:nfrglhls}
For $D,\hat\Theta,l$ there exist $\Theta,d$ such that:

If $q:H\to\Ga y$ is a $D$-straight holey line 
and $o_{x,y}\circ q:H\to\Ga x$ has $(\hat\Theta,l)$-regular large holes,
then the $\taumod$-Finsler line $\tau_-\tau_+$ in Proposition \ref{prop:hlflfns} 
can be chosen to be $\Theta$-regular.
\end{claim}
\proof The assertion follows from Claim \ref{claim:lngsgmnbgh}. \qed 

\medskip
If any two orbit points can be connected by such a holey line,
we obtain uniform regularity:

\begin{cor}
\label{cor:cntblnfrglhls}
Let $\Ga<G$ be $\taumod$-asymptotically embedded rel $\Ga\acts Y$. 
Suppose that there exist data $D,\hat\Theta,l$ and $y\in Y$ 
such that for each $\ga\in\Ga$
the points $y$ and $\ga y$ can be connected by a $D$-straight holey line $q:H\to\Ga y$
so that $o_{x,y}\circ q:H\to\Ga x$ has $(\hat\Theta,l)$-regular large holes.

Then $\Ga$ is uniformly $\taumod$-regular
and hence uniformly $\taumod$-straight rel $\Ga\acts Y$.
\end{cor}
\proof
The uniform straightness follows from uniform regularity together with the straightness proven earlier in Proposition \ref{prop:rbdstrght}.
\qed

\medskip
In the next section, we will apply this result in the relatively hyperbolic setting.

\subsubsection{Relatively hyperbolic subgroups}\label{sec:ae->unstr} 

We now restrict to relatively hyperbolic subgroups $\Ga<G$ 
and the case when $\Ga\acts Y$ is the action on a Gromov model. 
In this setting we obtain the following criterion for uniform regularity:

\begin{thm}
\label{thm:ae->str}
Suppose that $\Ga< G$ is relatively $\taumod$-asymptotically embedded 
and that each peripheral subgroup $\Pi_i< \Ga$ is uniformly $\taumod$-regular. 

Then $\Ga< G$ is uniformly $\taumod$-regular and hence relatively uniformly $\taumod$-straight.
\end{thm}
\proof First consider the case when $|\geo Y|\geq3$,
i.e. when Assumption~\ref{ass:rlsmbd} is satisfied.
In order to apply Corollary \ref{cor:cntblnfrglhls}, 
we need to check the connectability condition there 
for the orbits of the action $\Ga\acts Y$ on a Gromov model $(Y,{\mathcal B})$ of $(\Ga,\P)$.

We may assume that $y\in Y^{th}$.
Given $\ga\in\Ga$,
we connect $y$ and $\ga y$ by a geodesic in $Y$.
Along this geodesic, we choose a monotonic sequence of points $y_0=y,y_1,\ldots,y_n=\ga y$ in $Y^{th}$ 
so that any two successive points $y_{k-1}$ and $y_k$ have distance $\leq d_0$ for some fixed constant $d_0>0$
or lie on the same peripheral horosphere $\D B$, $B\in{\mathcal B}$.

Since the action $\Ga\acts Y^{th}$ is cocompact,
we may choose orbit points $\ga_k y$ 
at uniformly bounded distance from the points $y_k$. 
The sequence $(\ga_k y)$, viewed as a holey line $q:H=\{0,\ldots,n\}\to\Ga y$,
is $D$-straight with a constant $D$ independent of $\ga$.

Since also the actions $\Pi_i\acts\D B_i$ of the peripheral subgroups on the corresponding horospheres are cocompact,
and since there are finitely many conjugacy classes of peripheral subgroups,
there exists a finite subset $\Phi\subset\Ga$ independent of $\ga$ such that 
$$\ga_{k-1}^{-1}\ga_k\in \Phi(\bigcup_i\Pi_i)\Phi$$ for all $k$.

Due to our assumption that the subgroups $\Pi_i< G$ are uniformly $\taumod$-regular,
there exist $\hat\Theta$ and another finite subset $\Phi'\subset\Ga$, both independent of $\ga$, 
such that the pair of points $(\ga' x,\ga'' x)$ in $\Ga x$ is $\hat\Theta$-regular 
whenever ${\ga'}^{-1}\ga''\in \Phi(\bigcup_i\Pi_i)\Phi - \Phi'$.
This means that the holey line $o_{x,y}\circ q$, which corresponds to the sequence $(\ga_k x)$ in $\Ga x$,
has $(\hat\Theta,l)$-large holes 
for a sufficiently large constant $l$ independent of $\ga$.
Hence Corollary \ref{cor:cntblnfrglhls} implies the assertion.

If $|\geo Y|=2$, then ${\mathcal P}=\emptyset$ and we are in the absolute case.
There, asymptotic embeddedness implies Morse, and this in turn uniform regularity \cite{anolec}.

If $|\geo Y|=1$, then ${\mathcal P}=\{\Ga\}$ 
and the hypothesis of the theorem implies that $\Ga$ is uniformly $\taumod$-regular. 
\qed

\section{Comparing conditions for subgroups}
\label{sec:relation}

The following is the main theorem.
It summarizes the relation between different conditions on relatively hyperbolic subgroups established in the paper:

\begin{thm}\label{thm:main}
For subgroups $\Ga< G$, the following implications hold:

(i) relatively Morse $\Ra$ relatively uniformly Finsler-straight

(ii) relatively Finsler-straight $\Leftrightarrow$ relatively asymptotically embedded $\Leftrightarrow$ relatively RCA

(iii) relatively asymptotically embedded with uniformly regular peripheral subgroups $\Ra$ relatively uniformly Finsler-straight

(iv) relatively boundary embedded and Zariski dense $\Ra$ relatively asymptotically embedded
\end{thm}
\proof 
(i)
The maps of orbits are uniformly straight because Morse quasiisometric embeddings are.
Furthermore, since relatively Morse implies relatively asymptotically embedded, see Theorem~\ref{thm:RM->UR},
$\Ga$ has tied-up horospheres.

(ii) The first equivalence is the combination of Theorems~\ref{thm:str-ae} and~\ref{thm:rlsmbrlstr}.
The second equivalence is Yaman's theorem, see  Theorem \ref{thm:RRCA=RAE}. 

(iii) is Theorem \ref{thm:ae->str}. 

(iv) is Theorem~\ref{thm:BdZarden}. 
\qed

\medskip
In rank one,
all conditions become equivalent:

\begin{cor}\label{cor:rank1}
If the symmetric space $X$ has rank one,
then the following properties are equivalent for discrete subgroups $\Ga<G$:

(i) relatively Morse

(ii) relatively straight

(iii) relatively asymptotically embedded

(iv) relatively RCA

(v) relatively boundary embedded

(vi) geometrically finite
\end{cor}
\proof The implications (i)$\Ra$(ii)$\Leftrightarrow$(iii)$\Leftrightarrow$(iv)$\Ra$(v) hold in arbitrary rank.
In rank one, relative RCA amounts to the usual Beardon-Maskit
condition which is equivalent to geometric finiteness (see \cite{Bowditch}),
thus (iv)$\Leftrightarrow$(vi). 
Furthermore, (vi)$\Leftrightarrow$(i), see Theorem \ref{thm:RM=GF}. 

To get from (v) to the other conditions,
we observe that for non-elementary subgroups (v)$\Ra$(iii) holds 
because the limit set is the unique minimal nonempty $\Ga$-invariant closed subset of $\geo X$
and hence must equal the image of the boundary embedding. 
In the elementary case, we have (v)$\Ra$(vi) since in rank one 
all elementary discrete subgroups are geometrically finite. 
\qed

\section{Appendix 1: Geometry of nonpositively curved 
symmetric spaces and their ideal boundaries}

In this appendix we collect various definitions and facts about symmetric spaces of noncompact type, their ideal boundaries and their isometries. 

\medskip 
Let $X$ be a symmetric space of noncompact type. 
Throughout the paper we will be using two compactifications of $X$: 

1. The visual compactification (see \cite{BGS}), $\ol{X}= X\cup \geo X$, where $\geo X$ is defined as the set of {\em asymptotic} equivalence classes of geodesic rays in $X$: Two rays are equivalent if they are at finite Hausdorff distance from each other. We will also occasionally use the notion of {\em strongly asymptotic} geodesic rays : These are rays $\rho_i: \R_+\to X$, $i=1, 2$ such that
$$
\lim_{t\to\infty} d(\rho_1(t), \rho_2(t))=0. 
$$ 
The space of classes of strongly asymptotic rays is a topological space, obtained by taking the quotient space of the space of all geodesic rays in $X$. The latter is topologized by the topology of uniform convergence on compacts. Given a point $\xi\in \geo X$, one 
defines $X_{\xi}$, the {\em space of strong asymptote classes at $\xi$} as follows: Consider the set $Ray(\xi)$ consisting of all geodesic rays in $X$ asymptotic to $\xi$ (and, as before, equipped with the topology of uniform convergence on compacts). Then take the quotient of $Ray(\xi)$ by the equivalence relation, where two rays are equivalent if they are strongly asymptotic. One can identify the resulting space $X_\xi$ as follows. Pick a point $\hat\xi\in \geo X$ opposite to $\xi$ and consider the parallel set $P(\xi, \hat\xi)$ which is the  union of all geodesic rays in $X$ which are forward/backward asymptotic to $\xi, \hat\xi$. Each point $x\in P(\xi, \hat\xi)$ defines the ray $x\xi$ asymptotic to $\xi$. 
The projection of the subset of such rays to $X_\xi$ is a homeomorphism. Thus, we can identify $X_\xi$ with the parallel set $P(\xi, \hat\xi)$. 
One metrizes $X_\xi$ by
$$
d(\rho_1, \rho_2)= \lim_{t\to\infty} d(\rho_1(t), \rho_2(t)). 
$$
Then the projection $P(\xi, \hat\xi)\to X_\xi$ is an isometry. We refer to \cite[2.8]{anolec} for details.

Given a subset $A\subset X$ we let $\geo A$ denote the accumulation set of $A$ in $\geo X$.

2. The {\em Finsler bordification} $\ol{X}^{Fins}$ of $X$, also known as the {\em maximal Satake compactification} of $X$, obtained by equipping $X$ with a regular polyhedral Finsler metric $d^{\bar{\theta}}$, then compactifying $(X, d^{\bar{\theta}})$ by Finsler horofunctions on $X$.

{\bf Building notions.} 

The visual boundary $\geo X$ of a symmetric space $X$ admits a structure as a 
thick spherical building (the Tits building of $X$), which is a certain {\em spherical} simplicial complex, see \cite{Eberlein, BGS}. This complex is either connected (if $X$ has rank $\ge 2$) or discrete (if $X$ has rank $1$, equivalently, $X$ is negatively curved).  In the connected case,  this building is equipped with the path-metric induced by the spherical metrics on simplices. This metric space has diameter $1$. In the case when the building is discrete,  the distance between distinct 
points is $\pi$.  

The connected component $\Isom_0(X)$ of the isometry group of $X$ acts isometrically on this spherical building, and transitively on facets (top-dimensional simplices). The quotient 
$\geo X/\Isom_0(X)$ is then a single spherical simplex, denoted $\simod$, the {\em model spherical chamber} of the building $\geo X$. We will use the notation $\theta: \geo X\to \simod$ for the {\em type projection}, the quotient map $\geo X\to \geo X/\Isom_0(X)$. The full isometry group $\Isom(X)$ acts on both $\geo X$ and $\simod$ and the map $\theta$ is equivariant with respect to these actions. We let $G$ denote the kernel of the action of $\Isom(X)$ on $\simod$; 
i.e., $G$ is the subgroup of type preserving isometries.
It is a semisimple Lie group and has finite index in $\Isom(X)$.

For an algebraically inclined reader, the spherical building is defined as a simplicial complex whose vertices are conjugacy classes  of maximal parabolic subgroups of $G$. If $W$ is the Weyl group of $X$ (equivalently, the relative Weyl group of $G$) and $r$ is the rank of $X$ (equivalently, the split real rank of $G$), then $W$ is a reflection group acting 
 isometrically on the unit sphere $\amod$ of dimension  $r-1$ and $\simod$ is a fundamental chamber of this action. From the building viewpoint, $\amod$ is the {\em model spherical apartment} of the Tits building $\geo X$. 
 
 We return to the geometric discussion of $\geo X$. We let $\iota: \simod\to\simod$ denote the {\em opposition involution}: It is the projection to $\simod$ of {\em Cartan involutions} of $X$. We will use the notation $\bar\zeta=\theta(\zeta)$ for elements of $\simod$. Similarly (cf. \cite[\S 2.2.2]{anolec}), we let $\taumod\subseteq\simod$ denote  faces of $\simod$; these are the {\em model simplices} in the Tits building; we will use the notation $\tau$ for simplices in $\geo X$ 
of type $\taumod$, $\theta(\tau)=\taumod$. Given a face $\taumod$ of $\simod$, we let $\Dt \simod$ denote the union of faces of $\simod$ disjoint from $\taumod$; we use the notation $\ost(\taumod)$ for the complement $\simod \setminus \Dt \simod$, the {\em open star} of $\taumod$ in $\simod$. 

Let $W_{\taumod}$ denote the stabilizer of the face $\taumod$ in $W$. A (convex) subset $\Theta\subset \simod$ is said to be {\em Weyl-convex} if its $W_{\taumod}$-orbit is a convex subset in the sphere $\amod$ (see \cite[Def. 2.7]{anolec}). We will {\em always} use the notation $\Theta$ for compact $\iota$-invariant Weyl-convex subsets of the open star $\ost(\taumod)$.

 Note that, when $X$ has rank one, the data $\taumod, \Theta$ are obsolete since $\simod$ is a singleton. In this case, we also have that $\D\simod=\emptyset$ and $\Theta=\inte(\simod)=\simod$ is clopen. 
 
Two points in $\geo X$ are called {\em antipodal} if their distance in the Tits building $\geo X$ equals $\pi$. Equivalently, these points are connected by a geodesic in $X$. Equivalently, antipodal points are swapped by a Cartan involution of $X$. Similarly, two simplices in $\geo X$ are said to be antipodal if they are swapped by a Cartan involution. We will use the notation $\tau, \hat\tau$ for pairs of antipodal simplices. Their types in $\simod$ are swapped by the opposition involution $\iota$. Two simplices $\tau, \tau'$ in $\geo X$ are antipodal if and only if $\iota(\theta(\tau'))=\theta(\tau)$ and $\tau, \tau'$ contain antipodal points. 
 
Identify $\simod$ with a simplex in $\geo X$. Then $G$-stabilizers of faces $\taumod$ of $\simod$  
 are {\em standard parabolic subgroups} of $G$, they are denoted $P_{\taumod}$; these are closed subgroups of $G$. The set of simplices of type $\taumod$ in $\geo X$ is identified with the quotient 
 $G/P_{\taumod}$, which is a smooth compact manifold, called the {\em flag-manifold} $\Flagt$ of type $\taumod$, see \cite[\S 2.2.2, 2.2.3]{anolec}.  Given $\tau\in \Flagt$, one defines 
 the {\em open Schubert cell} $C(\tau)\subset \Flagt$, which is an open subset of  $\Flagt$ 
 consisting of elements opposite to $\tau$,  see \cite[\S 2.4]{anolec}. We will use the notation $\inte\tau$ for the {\em open simplex} obtained by removing from $\tau$ all its proper faces. The notation $\st(\tau)$ is used to denote the {\em star} of $\tau$ in $\geo X$, the union of all faces of 
 $\geo X$ containing $\tau$. Similarly, $\ost(\tau)$, the {\em open star} of $\tau$ in $\geo X$ is obtained from $\st(\tau)$ by removing all faces disjoint from $\tau$. Accordingly,
 $\Dt \simod$ is defined as $\simod \setminus \ost(\taumod)$. 
 
 Given a Weyl-convex compact subset $\Theta\subset \ost\taumod\subset \simod$, we will define the 
 {\em $\Theta$-star} of a simplex $\tau$ in $\geo X$ of type $\taumod$ as the preimage of 
 $\Theta$ under the restriction $\theta: \st(\tau)\subset \geo X\to \simod$.

\medskip
 {\bf Symmetric space notions.} An isometry of $X$ is called a {\em transvection} if it preserves a geodesic in $X$ and acts trivially on its normal bundle. Transvections in $X$ are precisely the compositions of pairs of Cartan involutions.

 A {\em flat} in $X$ is an isometrically embedded totally geodesic Euclidean subspace of $X$. The dimension of a {\em maximal flat} in $X$ is the {\em rank} of $X$. We will use the notation $F$ for maximal flats in $X$.% and $f$ for arbitrary flats. 
 The {\em parallel set} $P(l)$ of a geodesic $l$ in $X$ is the union of all maximal flats in $X$ containing $l$. We let $\tau_\pm$ denote the smallest simplices in $\geo X$ containing 
 the elements $\xi_\pm\in \geo l$; these simplices are antipodal. We will use the notation 
 $P(\tau_-,\tau_+)$ for the parallel set $P(l)$, since $P(l)$ can be described as the union of geodesics in $X$ forward/backward asymptotic to points in the open simplices $\inte\tau_\pm$.

Given a maximal flat $F\subset X$, the stabilizer $G_F$ of $F$ in $G$ acts transitively on $F$. 
 For each $x\in F$ the intersection $G_x\cap G_F$ acts on $F$ as a finite reflection group, isomorphic to the Weyl group $W$ of $X$. One frequently fixes a base-point $x=o\in X$ and the {\em model maximal flat} $\Fmod\subset X$ containing $o$; the stabilizer $G_x$ is then denoted $K$, it is a maximal compact subgroup of $G$. A fundamental domain for the $W$-action on $\Fmod$ is 
 the {\em model Euclidean Weyl chamber}, it is denoted $\Delta$. The ideal boundary $\geo \Delta$  is then identified with $\simod$, the model chamber in $\geo X$.

The Euclidean Weyl chamber $\De$ is the Euclidean cone over the simplex $\simod$. Thus, we can also ``cone-off" various objects from $\simod$. In particular, we define  the {\em $\taumod$-boundary} 
$\Dt\De$ of $\De$  as the union of rays $o\zeta$, $\zeta\in \Dt \simod$  (see \cite[\S 2.5.2]{anolec})
and the $\Theta$-cone $\Delta_\Theta$, as the union of rays $o\zeta$, $\zeta\in \Theta$, where 
$\Theta\subset \ost\taumod\subset \simod$.

%This subset of $\De$ will be used to define 

 The group of transvections along geodesics in $F$ is usually denoted $A$, then $G=KAK$ is the {\em Cartan decomposition} of $X$. The more refined form of this decomposition is $G=KA_+K$, where $A_+\subset A$ is the subsemigroup consisting of transvections mapping $\Delta$ into itself. 
 Geometrically speaking, the Cartan decomposition states that each $K$-orbit $Ky\subset X$ intersects  $\Delta$ in exactly one point. The projection $c: y\mapsto Ky\cap \Delta$ is 1-Lipschitz 
 since  each orbit $Ky$ meets $\Fmod$ orthogonally and transversally. This projection leads to the 
 notion of $\Delta$-distance on $X$ (see \cite[\S 2.6]{anolec} and \cite{KLM}): Given a pair of points $x, y\in X$, find $g\in G$ such that 
 $g(x)=o$ and $z=g(y)\in \Delta$. Then
 $$
 \ov{oz}:= d_\Delta(x,y). 
 $$  
 The vector $d_\Delta(x,y)$ is the complete $G$-congruence invariant of pairs $(x,y)\in X^2$ and 
 $$
 d(x,y)= ||d_\Delta(x,y)||
 $$ 
 where $||\cdot||$ is the Euclidean norm on $\Fmod$ (regarded as a Euclidean vector space with $o$ serving as zero).   Since $c$ is 1-Lipschitz, the $\Delta$-distance satisfies the triangle inequality
 $$
 ||d_\Delta(x,y) - d_\Delta(x,z)||\le ||d_\Delta(y,z)||= d(y,z). 
 $$
We refer the reader to \cite{KLM} for in-depth discussion of {\em generalized triangle inequalities} 
satisfied by the $\Delta$-valued distance function on $X$. 

We say that a nondegenerate segment $xy\subset X$ is $\taumod$-regular if 
$$
d_\Delta(x,y)\notin  \Dt\De
$$
 and is $\Theta$-regular if 
 $$
d_\Delta(x,y)\notin  \De_\Theta. 
$$
See \cite[\S 2.5.3]{anolec}. Suppose that $\taumod$ is $\iota$-invariant and $\Theta\subset \ost(\taumod)$ is an $\iota$-invariant Weyl-convex compact subset of $\ost(\taumod)\subset \simod$.  

A {\em Finsler geodesic} $\varphi$ in $X$ is said to be 
{\em $\taumod$-regular}, resp. $\Theta$-regular, or, simply, a $\Theta$-Finsler geodesic if it is $\taumod$-regular, respectively, $\Theta$-regular, as a subset of $X$.   

 \medskip
 {\bf Weyl cones and diamonds in $X$.} Fix a simplex $\tau$ in $\geo X$ and a point $x\in X$. 
 We define the {\em Weyl cone $V(x, \st(\tau))$} with the tip $x$ over the star $\st(\tau)\subset \geo X$ as the union of geodesic rays emanating from $x$ and asymptotic to points of $\st(\tau)$. Similarly, assuming that $\tau$ has the type 
  $\taumod$  and $\Theta\subset \ost(\taumod)\subset \simod$ is a Weyl-convex compact subset, we define the $\Theta$-cone  $V(x, \st_\Theta(\tau))$ with the tip $x$ over the $\Theta$-star 
 $\st_\Theta(\tau)\subset \geo X$ as the union of geodesic rays emanating from $x$ and asymptotic to points of $\st_\Theta(\tau)$. It was proven in  \cite[\S 2.5]{anolec} that Weyl cones  
 $V(x, \st(\tau))$ and $\Theta$-cones $V(x, \st_\Theta(\tau))$ are convex in $X$. 
 
 {In particular, they satisfy the {\em nested cones property:}}
 
{1. If $y\in  V(x, \st(\tau))$ then $V(y, \st(\tau))\subset V(x, \st(\tau))$.} 
 
 {2. If $y\in  V(x, \st_\Theta(\tau))$ then $V(y, \st_\Theta(\tau))\subset V(x, \st_\Theta(\tau))$.}

Intersecting cones, we define {\em diamonds} in $X$ (see \cite[\S 2.5]{anolec}). 
Take two antipodal simplices $\tau_+, \tau_-$ of the type 
$\taumod= \iota\taumod$ and points $x_\pm \in X$ such that 
$$
x_+\in V(x_-, \st(\tau_+)), x_-\in V(x_+, \st(\tau_-)). 
$$ 
 Then the intersection 
 $$
 \diamot(x_-, x_+)= V(x_-, \st(\tau_+)) \cap V(x_+, st_\Theta(\tau_-))
 $$
 is the $\taumod$-diamond with the tips $x_\pm$. Similarly, suppose that 
 $\Theta\subset \ost(\taumod)\subset \simod$ is a Weyl-convex $\iota$-invariant compact subset and 
 $$
x_+\in V(x_-, \st_\Theta(\tau_+)), x_-\in V(x_+, \st_\Theta(\tau_-)). 
$$ 
 Then  the $\Theta$-diamond is defined as the intersection 
  $$
 \diamoTh(x_-, x_+)= V(x_-, st_\Theta(\tau_+)) \cap V(x_+, st_\Theta(\tau_-)). 
 $$
 Thus, diamonds are also convex in $X$. 
 
 Diamonds have a nice interpretation in terms of Finsler geometry of $X$: It is proven in \cite{bordif} that  $\diamot(x_-, x_+)$ is the union of all Finsler geodesics in $X$ connecting 
 the points $x_\pm$. Similarly, $\diamoTh(x_-, x_+)$ 
 is the union of all $\Theta$-regular Finsler geodesics in $X$ connecting  the points $x_\pm$. 

\medskip 
 {\bf Regular sequences and groups.} 
 
The reader should think of {\em regularity} conditions for subgroups $\Ga< G$ as a way of strengthening the discreteness assumption: The discreteness condition means that 
the sequence of distances $d(x, \ga_n(x))$ diverges to infinite for every sequence of distinct elements $g_n\in \Ga$. Regularity conditions on $\Ga$ require certain  forms of divergence to infinity of vector-valued distances $d_\De(x, g_n x)$.  
 
We first consider sequences in the euclidean model Weyl chamber $\De$.
Recall that $\Dt\De=V(0,\Dt\simod)\subset\De$ 
is the union of faces of $\De$ which do not contain the sector $V(0,\taumod)$.
Note that $\Dt\De\cap V(0,\taumod)=\D V(0,\taumod)=V(0,\D\taumod)$. 
The following definitions are  taken from \cite[\S 4.2]{anolec}.  

\begin{dfn}
(i) A sequence $(\de_n)$ in $\De$ is called 
 {\em $\taumod$-regular} if it drifts away from $\Dt\De$, i.e. 
$$ d(\de_n,\Dt\De) \to+\infty .$$

(ii) A sequence $(x_n)$ in $X$ is {\em $\taumod$-regular} 
if for some (any) base point $o\in X$ the sequence of $\De$-distances $d_{\De}(o,x_n)$ in $\De$
has this property.

(iii) A sequence $(g_n)$ in $G$ is {\em $\taumod$-regular}, 
if for some (any) point $x\in X$ the orbit sequence $(g_nx)$ in $X$ has this property.

(iv) A subgroup $\Ga<G$ is {\em $\taumod$-regular}
if all sequences of distinct elements in $\Ga$ have this property.
\end{dfn}

Next, we describe a stronger form of regularity following \cite[\S 4.6]{anolec}.

\begin{dfn}
\label{def:unifreg}
(i) A sequence $\de_n\to\infty$ in $\De$ is 
{\em uniformly $\taumod$-regular} if it drifts away from $\Dt\De$
at a linear rate with respect to 
its norm,
$$ \liminf_{n\to+\infty} \frac{d(\de_n,\Dt\De)}{\|\de_n\|} > 0.$$

(ii) A sequence $(x_n)$ in $X$ is {\em uniformly $\taumod$-regular}
if for some (any) base point $o\in X$ the sequence of $\De$-distances $d_{\De}(o,x_n)$ in $\De$
has this property.

(iii) A sequence $(g_n)$ in $G$ is {\em uniformly $\taumod$-regular}
if for some (any) point $x\in X$ the orbit sequence $(g_nx)$ in $X$ has this property.

(iv) A subgroup $\Ga<G$ is {\em uniformly $\taumod$-regular}
if all sequences of distinct elements in $\Ga$ have this property.
\end{dfn}

Note that (uniform) regularity of a sequence in $X$ 
is independent of the base point and stable under bounded perturbation of the sequence 
(due to the triangle inequality for $\De$-distances). A sequence $(x_n)$ is uniformly $\taumod$-regular if and only if there exists a compact $\Theta\subset \ost(\taumod)$ such that for each 
$x\in X$ all but finitely many vectors $d_\De(x, x_n)$ belong to $\Delta_\Theta$. 

\begin{rem}
{The definition of regularity of sequences in $G$ has the following dynamical interpretation, in terms of dynamics on the flag-manifold $\Flagt$, which generalizes the familiar convergence property for sequences of isometries of a Gromov-hyperbolic space $Y$ acting on the visual boundary of $Y$.} 

{A sequence $(g_n)$ in $G$ is said to be {\em $\taumod$-contracting} if there exists a pair of elements $\tau_\pm \in \Flagt$, such that the sequence $(g_n)$ converges to $\tau_+$ uniformly on compacts in $C(\tau_-)\subset \Flagt$. Here, as elsewhere, $C(\tau_-)$ is the open Schubert cell in $\Flagt$ consisting of simplices antipodal to $\tau_-$. In this situation, the simplex $\tau_+$ is the $\taumod$-limit of $(g_n)$. 
A sequence $(g_n)$ is $\taumod$-regular if and only if there exists a pair of bounded sequences $a_n, b_n\in G$ such that the sequence of compositions $c_n g_n b_n$ is $\taumod$-contracting. Equivalently, a sequence $(g_n)$ is $\taumod$-regular if and only if every subsequence in $(g_n)$ contains a further subsequence which is $\taumod$-contracting. We refer the reader to \cite{anolec} for details.}
\end{rem}

For a subgroup $\Ga<G$, uniform $\taumod$-regularity is equivalent to 
the {\em visual limit set} $\La\subset\geo X$ being contained 
in the union of the open $\taumod$-stars, where $\La$ is the accumulation set in $\geo X$ 
of some (any) $\Ga$-orbit in $X$. 

A subgroup $\Ga<G$ is (uniformly) $\taumod$-regular iff it is 
(uniformly) $\iota\taumod$-regular.

\section{Appendix 2 (by Grisha Soifer): Auslander's Theorem}\label{sec:Auslander}

\begin{thm} \label{thm:auslander}
Let $G$ be a Lie group which splits as a semidirect product 
$G= N\rtimes K$ where $N$ is connected nilpotent and $K$ is compact. 
Then each discrete subgroup $\Ga< G$ is finitely generated and virtually nilpotent. 
\end{thm}
\proof This theorem first appeared in Auslander's paper \cite{Au1}, but its proof was flawed. The theorem can be derived  from a more general result \cite{BK} (in the torsion-free case). 
Notice that in the paper we are only interested in subgroups of finitely generated subgroups linear groups, which are virtually torsion-free by the Selberg Lemma. Hence, in this setting, Auslander's Theorem can be viewed as a corollary of \cite{BK}.  If $\Ga$ is assumed to be finitely generated then 
Auslander's Theorem can be viewed as a corollary of Gromov's Polynomial Growth Theorem (which is much easier in this setting since $\Ga$ 
is already assumed to be a subgroup of a connected Lie group).

Nevertheless, for the sake of completeness, we present a direct and 
self-contained proof, which is well-known to  experts, but which we could not find in the literature. 

\medskip 
{\textit{Step 1.}} Let us show that, after passing to a finite index subgroup in $\Gamma$, 
we can assume that $\Gamma$ is a discrete subgroup of a closed connected solvable subgroup of $G$. 
Indeed, since $K$ is compact, the quotient $K/K^0$ is compact, where $K^0$ is the identity component of $K$;  hence $\Gamma \cap N K^0$ is a finite index subgroup in $\Gamma$. 
Therefore we can assume that $K$ is a connected compact Lie group. 
Let $R$ be the solvable radical of $G$. Obviously,
$R= N L$, where $L$ is an abelian normal compact subgroup of $K$. 
Let  $\pi : G \rightarrow G/R$ be the quotient homomorphism.  
By another Auslander theorem, \cite[Theorem 8.24]{R},  the connected component of the closure  
$\overline{\pi(\Gamma)}$ in  $G /R$ is an abelian group. 
Since $G/N$ is compact,  $G /R$ is compact as well.  
Let $\overline{\Gamma} = \pi^{-1} (\overline{\pi(\Gamma)}^0) \cap \Ga$. Clearly,  $\overline{\Gamma} $ is a finite index subgroup of $\Gamma$ and   $\overline{\Gamma} $ is a subgroup of the solvable group
 $\pi^{-1} (\overline{\pi(\Gamma)}^0)$. \\ 
 
 {\textit{Step 2.}} From now on, we will assume that $G$ is a solvable connected group, $G =N\rtimes K$, where $N$ is a connected nilpotent   group and $K$ is a compact abelian group.  
\par  Let us show that we can assume that $N$ is simply connected. 

\begin{lem}\label{lem:step1}
For every  connected nilpotent group Lie $N$ there exists a maximal compact subgroup $T< N$ which is characteristic in $N$, such that the quotient $N/T$ is a simply connected nilpotent Lie group.  
\end{lem} 
\proof Consider a  compact subgroup $C < N$ and  the adjoint representation $Ad_N: N \rightarrow GL({\mathfrak n})$, where ${\mathfrak n}$ is the Lie algebra of $N$. Since $C$ is a  compact group,  
for every $c\in C$ we have that $Ad_N (c)$ is a semisimple linear transformation. 
On the other hand,  $Ad_N (c)$ is nilpotent since the group $N$ is nilpotent. Therefore  $Ad_N (c) ={\mathbf 1}$, 
hence, $C \subseteq Z(N)$. Thus, each compact subgroup of $N$ is contained in the center $Z(N)$ of $N$ 
and, by commutativity, the union of all compact subgroups of $Z(N)$ forms a compact subgroup $T< Z(N)$. 
Since the center is a characteristic subgroup and each automorphism of $N$ sends compact subgroups to compact subgroups, $T$ is a characteristic subgroup of $N$.   
\qed  
 
 \medskip
 In our setting, the maximal compact subgroup $T< N$ will be a normal subgroup of $G$. 
 By the compactness of $T$ and discreteness of $\Ga$, the intersection $\Gamma \cap T$ is a finite subgroup of $\Gamma$. Consider the quotient homomorphism $\pi : G \rightarrow G_1=G/T$. The kernel of $\pi|_\Gamma$ is a finite normal subgroup of $\Gamma$. The quotient $\Gamma_2=\pi(\Gamma)$ is a discrete subgroup of the group $G_1= N_1 K_1$, where $K_1= \pi(K)$ and $N_1=\pi(N)$ is a simply connected nilpotent group. We will work with the subgroup $\Ga_1:=\pi(\Gamma)$ of $G_1$. Once we know that $\Gamma_1$ is virtually a uniform lattice in a connected nilpotent Lie group, it is finitely generated and, hence, polycyclic. It then will follow that  the  group $\Gamma$ itself is  residually finite according to \cite{Hirsch}; cf. Lemma 11.77 in \cite{Drutu-Kapovich}. Thus, $\Gamma$ contains a finite index subgroup $\tilde\Gamma_1$ such that $\pi: \tilde\Gamma_1\to \Gamma_1$ is injective. This reduces the problem to the case when $N$ is simply connected and $K$ is compact abelian. 
 
Since $G$ is a connected solvable group with simply connected  nilpotent radical $N$,  there exists a faithful linear representation $\rho: G \rightarrow  GL(d,\R) $ such that $\rho(n)$ is a unipotent matrix for every $n \in N$; see e.g. \cite{Ho} and also \cite{S1}. 
We will identity $G$ with $\rho(G)$. Then $N$ is an algebraic subgroup of $GL(d, \R)$ (since $N$ is unipotent, both the exponential and logarithmic maps of $N$ are polynomial). 
 
 \medskip 
 {\textit{Step 3.}} This is the key step in the proof. Let $N_2 < N$ be the Zariski closure of $\Gamma \cap N$. Since $K$ is abelian, we have  
 $[\Gamma, \Gamma] < N$; since $N$ is normal in $G$,  $\Gamma$ normalizes 
 $N_2$.\\
 
 Recall that a discrete subgroup of a simply connected (algebraic) nilpotent group $H$ is Zariski dense in 
 $H$  if and only if it is a cocompact lattice in $H$, see \cite[Theorem 2.3, page 30]{R}.  
 Therefore, in our case,  $\Gamma \cap N= \Gamma \cap N_2$  is a cocompact lattice in $N_2$. In particular, this intersection is finitely generated. 
As we noted above, the subgroup $\Ga$ normalizes $N_2$.  Our next goal is to prove that $N_2\Ga= N_2\rtimes \Ga$ is a closed subgroup of $G$.  This will be a corollary of a more general lemma about Lie groups:

\begin{lemma}
Let $H, N$ be closed subgroups of a Lie group $G$, such that 
the intersection $N\cap H$ is a cocompact subgroup of $N$ and $H$ normalizes $N$. 
Then the subgroup $NH$ is closed in $G$. 
\end{lemma}
 \proof Let $K\subset H$ be a compact such that $HK=N$. Consider a convergent sequence $(g_i)$ 
 in $NH$. Then every $g_i$ can be written as $g_i= k_i h_i$, $k_i\in K\subset N, h_i\in H$.  After extraction, $k_i\to k\in K\subset N$, hence, $h_i\to h$ and, since $H$ is closed, $h$ is in $H$. Therefore, 
 $$
 g_i\to kh\in NH. \qed 
 $$
 
 Specializing to our situation, where the role of the closed subgroup $H$ is played by $\Ga$, and taking into account that $N_2$ is normalized by $\Ga$,  we obtain

\begin{cor}
$\Gamma N_2< G$ is a closed Lie subgroup in $G$ (in the classical topology). 
\end{cor}

 \begin{cor}
 The identity component of $\Gamma N_2$ coincides with that of $N_2$.
 \end{cor}
 \proof This follows from countability of $\Ga$. \qed  
 
 \medskip 
 Let ${\mathfrak n}_2$ be the Lie algebra of $N_2$. By the above observation about the identity 
 component, the Lie algebra of $\Gamma N_2$ coincides with  ${\mathfrak n}_2$.  Let 
 $Ad : \Gamma N_2 \rightarrow GL ({\mathfrak n}_2)$ be the adjoint representation. We will use the following idea due to Margulis: There exists a full-rank  lattice $\Delta$ in the vector space ${\mathfrak n}_2$ such that $\exp (\Delta)$ is a finite index subgroup of  $\Gamma \cap N_2$;  see e.g. \cite[sect. 3.1]{Margulis}.

 The number of subgroups of the given index in a finitely generated group, such as $\Gamma\cap N$, is finite.  
 Therefore, taking into account the fact that $\Ga\cap N$ is normal in $\Ga$, 
 there exists a subgroup  of finite index  $\overline{\Gamma}$ of $\Gamma$ such that $\exp \Delta$ is $\overline{\Gamma}$-invariant. 
 Since  $\mid\Gamma : \overline{\Gamma} \mid < \infty $ we can and will assume that $\exp \Delta$ is $\Gamma$--invariant, i.e. is invariant under the action of $\Ga$ by conjugation on $N$. 
 Thus $\Delta$ is $Ad(\Gamma)$-invariant, and, by identifying $\Delta$ with $\Z^m$, $m = \dim {\mathfrak n}_2$,  we have $Ad \gamma \in  GL (m, \mathbb{Z})$ for 
 each $ \gamma \in \Gamma$.    After passing to a further finite index subgroup of 
 $\Gamma$, we can assume that $Ad (\Ga) <  SL (m, \mathbb{Z})$.

 Consider the Jordan decomposition $Ad \gamma = \gamma_s \gamma_u$  of  
 $Ad \gamma$,  $\gamma \in \Gamma$, where  $ \gamma_s$ is the semisimple and $ \gamma_u$ 
 is the unipotent part of the decomposition. Since the maps $Ad \gamma \mapsto \gamma_u$ and   
 $Ad \gamma \mapsto \gamma_s$  are restrictions of  ${\mathbb Q}$--rational maps 
 $GL(m, {\mathbb R})\to GL(m, {\mathbb R})$ we have  $ \gamma_u \in  SL (n , \mathbb{Q})$ and $ \gamma_s \in  SL (m , \mathbb{Q})$.  (See \cite[p. 158]{AM}.)

\begin{lem} For each unipotent element $u\in GL(m,\Q)$, there exists $q=q_u\in \N$ such that  
$u^q \in GL(m,\Z)$. 
 \end{lem} 
 \proof There exists $q_1 =q_1(m)$ such that for  each unipotent element $u\in GL(m, \R)$,
 $$
 ({\mathbf 1}-u)^{q_1 +1} =0. 
 $$

Assume now that  $u \in GL(m, \mathbb{Q})$; let $M$ denote 
the product of  denominators of all matrix entries of $u_1, \ldots, u_1^{n_1}$ and 
set  $q_2 := M\cdot q_1!$, $u_1 := {\mathbf 1} -u$. Then  
$$
u^{q_2}= \sum_{k=0}^{q_2} {q_2 \choose k} u^{k}_1=  \sum_{k=0}^{q_1} {q_2 \choose k} 
u^{k}_1. 
$$
Clearly,
$$
{q_2 \choose k} u^k_1 \in Mat(m, \mathbb{Z})
$$
for all $k \leq n_1$. Hence, for $q=q_2$, $u^q\in GL(m,\Z)$. \qed

 \medskip 
 In our case, given that  $ \gamma_u \in  SL (m , \mathbb{Q})$, we conclude that there exists a positive integer $q_{\gamma}$ such that $\gamma_u ^{q_{\gamma}} \in  SL (m, \mathbb{Z})$.   Since $\gamma_s \gamma_u =\gamma_u \gamma_s $ and $Ad \gamma \in  SL (m, \mathbb{Z})$, we have $\gamma_s ^{q_{\gamma}} \in  SL (m, \mathbb{Z})$.   
 On the other hand,  every $ \gamma \in \Gamma$ is the product $ \gamma = n k$ where 
 $ n \in N, k \in K$. As we noted above,  $Ad (n)$ is unipotent; 
 by compactness of $K$,  for every eigenvalue $\lambda$ of   
 $\gamma_s $ we have $|\lambda| =1$.   Hence  $\gamma_s ^{q_{\gamma}} $ is a finite order element of the discrete group $ SL (m , \mathbb{Z})$. Thus there exists $p_{\gamma}$ such that  $\gamma_s ^{p_{\gamma}}  = {\mathbf 1}$.  From this it follows that  $\gamma^{p_{\gamma}} \in N$.  
  
 \begin{lem}
$\Ga$ is finitely generated. \end{lem} 
\proof The subgroup $\Gamma \cap N_2$ is finitely generated because it is a cocompact  lattice. For the same reason, the projection of $\Ga$ to the connected abelian group $G/N_2$ is discrete. Therefore, this projection is 
finitely generated as well. It follows that $\Ga$ itself is finitely generated. \qed

\medskip  
Therefore, since the projection of $\Ga$ to $G/N_2$ is a finitely generated torsion group,  this projection has to be finite. 
In particular, $\Ga\cap N_2$ has finite index in $\Ga$. 
Since $\Ga\cap N_2$ is a cocompact lattice,  Theorem \ref{thm:auslander} follows.   \qed

Addresses:

\noindent M.K.: Department of Mathematics, \\
University of California, Davis\\
CA 95616, USA\\
email: kapovich@math.ucdavis.edu

\noindent
KIAS, 85 Hoegiro, Dongdaemun-gu,\\ 
Seoul 130-722, South Korea\\

\noindent B.L.: Mathematisches Institut\\
Universit\"at M\"unchen \\
Theresienstr. 39\\ 
D-80333, M\"unchen, Germany\\ 
email: b.l@lmu.de


\begin{thebibliography}{BLP05}

\bibitem[AM]{AM} 
H. Abbaspur, M. Moskowitz,  ``Basic Lie theory,'' World Scientific, 2007.  

\bibitem[Ah]{Ahlfors}
L. Ahlfors,  {\em Fundamental polyhedrons and limit point sets of Kleinian groups}, 
Proc. Nat. Acad. Sci. U.S.A. {\bf 55} (1966) pp. 251--254.


\bibitem[Au1]{Au1} 
L. Auslander, {\em Bieberbach's theorems on space groups and discrete uniform subgroups of solvable Lie groups II}, Amer. J. Math. {\bf 83} (1961) pp. 276-- 280.

\bibitem[Au2]{Au2} 
L. Auslander, {\em On radicals of discrete subgroups of Lie groups,} Amer. J. Math., {\bf 85} (1963) no. 2,  pp. 145--150.

\bibitem[Ba]{Ballmann}
W. Ballmann,
``Lectures on spaces of nonpositive curvature". With an appendix by Misha Brin. 
DMV Seminar, vol.\ 25, Birkh\"auser Verlag, Basel, 1995.

\bibitem[BGS]{BGS}
W. Ballmann, M. Gromov and V. Schroeder, ``Manifolds of nonpositive curvature'', Birkh\"auser Verlag, 1985. 



\bibitem[BM]{BM}
A. Beardon, B. Maskit, 
{\em Limit points of Kleinian groups and finite sided fundamental polyhedra},  
Acta Math. {\bf 132} (1974) pp. 1--12. 


\bibitem[BK]{BK}
I. Belegradek, V. Kapovitch, 
{\em Classification of negatively pinched manifolds with amenable fundamental groups},  
Acta Math. {\bf 196} (2006), no. 2, pp. 229--260. 

\bibitem[Be]{Benoist}
Y.\ Benoist,
{\em Propri\'et\'es asymptotiques des groupes lin\'eaires}, Geom. Funct. Anal.  {\bf 7} (1997), no. 1, pp. 1--47.



\bibitem[Bou]{Bourdon}
M.\ Bourdon,
{\em Structure conforme au bord et flot g\'ed\'esique d'un $CAT(-1)$-espace},
Enseign. Math. (2) {\bf 41} (1995), no. 1-2, pp. 63--102. 

\bibitem[Bo1]{Bowditch93}
B.\ H. Bowditch,  {\em Geometrical finiteness for hyperbolic groups}, J. Funct. Anal. {\bf 113} (1993), no. 2, pp. 245--317.


\bibitem[Bo2]{Bowditch}
B.\ H. Bowditch, {\em Geometrical finiteness with variable negative curvature},
Duke Math. J.  {\bf 77} (1995) pp. 229--274.

\bibitem[Bo3]{Bowditch1998}
B. H. Bowditch, {\em  Spaces of geometrically finite representations}, 
Ann. Acad. Sci. Fenn. Math. {\bf 23} (1998), no. 2, pp. 389--414.
 
 \bibitem[Bo4]{Bowditch_config}
B.\ H. Bowditch, 
{\em Convergence groups and configuration spaces}, 
in ``Geometric group theory down under'' (Canberra, 1996), 
pp.\ 23--54, de Gruyter, Berlin, 1999. 



\bibitem[Bo5]{Bowditch2012}
B. H. Bowditch, {\em Relatively hyperbolic groups},  
Internat. J. Algebra Comput. {\bf 22} (2012), no. 3, 1250016, 66 pp.

\bibitem[BH]{BH}
M. Bridson and A. Haefliger, ``Metric spaces of non-positive curvature,'' Springer- Verlag, Berlin, 1999.

\bibitem[Bu]{Bumagin}
I. Bumagin, {\em On definitions of relatively hyperbolic groups}, in 
``Geometric Methods in Group Theory," Contemporary Mathematics, Vol. 372 (American Mathematical Society, 2005), pp. 189--196.  

\bibitem[BS]{BS}
S. Buyalo, V. Schroeder,
{\em Elements of asymptotic geometry}, 
EMS 2007.

\bibitem[CLT]{CLT}
D. Cooper, D. Long, S. Tillmann, {\em On convex projective manifolds and cusps}, 
Advances in Mathematics, {\bf 277} (2015) pp. 181--251. 

\bibitem[CDP]{CDP}
M. Coornaert, T. Delzant, A. Papadopoulos, ``Geom{\'e}trie et th{\'e}orie des groupes: 
Les groupes hyperboliques de Gromov'', 
Lecture Notes in Math. 1441, Springer, Berlin (1990). 

\bibitem[DY]{DY}
F. Dahmani,  A. Yaman, 
{\em Bounded geometry in relatively hyperbolic groups},  
New York J. Math. {\bf 11} (2005) pp. 89--95. 


      
\bibitem[DK]{Drutu-Kapovich}
C. Drutu, M. Kapovich, ``Geometric Group Theory''.  Colloquium Publications, AMS, Vol. 63, 2018.  
       
\bibitem[DS]{Drutu-Sapir}
C. Drutu, M. Sapir, {\em Tree-graded spaces and asymptotic cones of groups}, Topology {\bf 44} (2005), no. 5, pp. 959--1058.       
       
\bibitem[E]{Eberlein}
P. Eberlein, ``Geometry of nonpositively curved manifolds'', University of Chicago Press, 1997. 

\bibitem[F]{Farb}
B. Farb, {\em Relatively hyperbolic groups}, Geom. Funct. Anal. {\bf 8} (1998), no. 5, pp. 810--840. 
       
\bibitem[FG]{FG}   
V.  Fock, A. Goncharov, {\em Moduli spaces of local systems and higher Teichm\"uller theory}, 
Publ. Math. IHES, {\bf 103} (2006) pp. 1--211.    
       
\bibitem[Ge]{Gerasimov}       
V. Gerasimov, 
{\em Expansive convergence groups are relatively hyperbolic}, 
Geom. Funct. Anal. {\bf 19} (2009), no. 1, pp. 137--169.        
       
       
\bibitem[GP1]{Gerasimov-Potyagailo}       
V. Gerasimov, L. Potyagailo, 
{\em Non-finitely generated relatively hyperbolic groups and Floyd quasiconvexity}, 
Groups, Geometry and Dynamics, {\bf 9} (2015), no. 2, pp. 369--434.    
             
\bibitem[GP2]{GerasimovPotyagailo2}       
V. Gerasimov, L. Potyagailo, 
{\em  Similar relatively hyperbolic actions of a group}, 
Int. Math. Res. Not. IMRN 2016, no. 7, pp. 2068--2103. 

\bibitem[GdlH]{GhysdelaHarpe}
E. Ghys, P. de la Harpe, eds. ``Sur les groupes hyperboliques d'apr{\`e}s Mikhael Gromov.'' 
Progress in Mathematics. Birkh{\"a}user Boston, Inc., 1990.
              
\bibitem[G]{Gromov}
M. Gromov, {\em Hyperbolic groups}, In: ``Essays in group theory,'' Math. Sci. Res. Inst.
Publ. 8, Springer, New York (1987) pp. 75--263. 


\bibitem[GW]{GW}
O.\ Guichard, A.\ Wienhard, 
{\em Anosov representations: Domains of discontinuity and applications}, 
Invent. Math. {\bf 190} (2012) no. 2, pp. 357--438. 

\bibitem[Gui]{Guivarch}
Y. Guivarc'h, {\em Produits de matrices al\'eatoires et applications}, Erg. Th. Dyn. Sys. {\bf 10} (1990) pp. 483--512.

\bibitem[He]{Helgason}
S. Helgason, ``Differential geometry, Lie groups and symmetric spaces'',  AMS series Graduate Studies in Mathematics, 2001. 


\bibitem[Hi]{Hirsch}
K.~A. Hirsch, \emph{On infinite soluble groups, {I}{I}}, Proc. Lond. Math. Soc.
  {\bf 44} (1938) pp. 336--344.



\bibitem[Ho]{Ho}
G. Hochschild, {\em On representing analytic groups with their automorphisms}, Pacific Journal of Math. {\bf 21} 
(1978)  pp. 333--336. 


\bibitem[Hr]{Hruska}
C. Hruska, {\em Relative hyperbolicity and relative quasiconvexity for countable groups}, Algebr. Geom.
Topol. {\bf 10} (2010) pp. 1807--1856.

\bibitem[J]{Jacobson}
N. Jacobson, 
{\em Completely reducible Lie algebras of linear transformations},
Proc. Amer. Math. Soc. {\bf 2} (1951) pp. 105--113.

\bibitem[KL1]{bordif}
M.\ Kapovich, B.\ Leeb, {\em Finsler bordifications of symmetric and certain locally symmetric spaces}, Geometry and Topology, {\bf 22} (2018) 2533--2646. 



\bibitem[KL2]{manicures}
M. Kapovich, B. Leeb, 
{\em Discrete isometry groups of symmetric spaces}, 
MSRI Lecture Notes. ArXiv e-print, arXiv:1703.02160, 2017. 
To appear in ``Handbook of Group Actions, IV". 



\bibitem[KL3]{relmorse-2}
M. Kapovich, B. Leeb, {\em Relativizing characterizations of Anosov subgroups, II,} in preparation.

\bibitem[KLM]{KLM} M. Kapovich, B. Leeb and J. J. Millson, {\em Convex functions on symmetric spaces, side lengths of polygons and the stability inequalities for weighted configurations at infinity},  
Journal of Differential Geometry,  vol. 81, 2009, p. 297-- 354. 


\bibitem[KLP1]{morse}
 M. Kapovich, B. Leeb and J. Porti, Morse actions of discrete groups on symmetric spaces, arXiv  e-print, 2014.  arXiv:1403.7671. 

\bibitem[KLP2]{mlem}
M. Kapovich, B. Leeb and J. Porti, {\em A Morse lemma for quasigeodesics in symmetric spaces and euclidean buildings}, Geometry and Topology, {\bf 22} (2018) pp. 3827--3923.


\bibitem[KLP3]{anosov}
M.\ Kapovich, B.\ Leeb, J.\ Porti, 
{\em Some recent results on Anosov representations}, Transformation Groups, {\bf 21} (2016) no. 4, pp. 1105--1121. 

\bibitem[KLP4]{coco15}
M.\ Kapovich, B.\ Leeb, J.\ Porti,
{\em Dynamics on flag manifolds: domains of proper discontinuity and cocompactness},
Geometry and Topology, {\bf 22} (2017) pp. 157--234. 


\bibitem[KLP5]{anolec} M. Kapovich, B. Leeb and J. Porti,  
{\em Anosov subgroups: dynamical and geometric characterizations}, European Journal of Mathematics, {\bf 3} (2017) pp. 808--898.


\bibitem[KiLe]{Kim-Lee}
S. Kim, J. Lee, {\em Three-holed sphere groups in $PGL(3, \R)$}, in preparation. 

\bibitem[La]{Labourie}
F. Labourie, {\em Anosov flows, surface groups and curves in projective space}, 
Invent. Math. {\bf 165} (2006) no. 1, pp. 51--114.

\bibitem[Le]{habil}
B.\ Leeb, 
{\em A characterization of irreducible symmetric spaces 
and Euclidean buildings of higher rank by their asymptotic geometry}, 
Bonner Mathematische Schriften 326 (2000), 
see also  arXiv:0903.0584 (2009).


\bibitem[Mar]{Marden}
A. Marden, {\em The geometry of finitely generated Kleinian groups}, Ann. of Math. (2) {\bf 99} (1974) pp. 383--462. 

\bibitem[M]{Margulis} 
G. Margulis, {\em Arithmetic properties of discrete subgroups}, 
Russian Math. Surveys {\bf 29:1} (1974)  pp. 107--156. 	


\bibitem[O]{Osin}
D. Osin, {\em Relatively hyperbolic groups: intrinsic geometry, algebraic properties and algorithmic problems},
Mem. AMS {\bf 179} (2006) no. 843 vi+100pp. 

\bibitem[Pr]{Prasad}
G. Prasad, 
{\em R-regular elements in Zariski-dense subgroups}, 
Quart. J. Math. Oxford Ser. (2) {\bf 45} (1994) pp. 541--545. 


\bibitem[R]{R} 
H. Raghunathan,  ``Discrete subgroups of Lie groups,'' 
Berlin-Heidelberg, Springer,  (1972).

\bibitem[Ra]{Rat} 
J. Ratcliffe, ``Foundations of hyperbolic manifolds''. Second edition. Graduate Texts in Mathematics, 149. Springer, New York, 2006.

\bibitem[Ro]{Ronan}
M. Ronan, {\em Buildings: main ideas and applications. II. Arithmetic groups, buildings and symmetric spaces,} Bull. London Math. Soc., 24 (1992) pp. 97--126. 

\bibitem[Sch]{Sch}
R. Schwartz, {\em Pappus' theorem and the modular group}, Inst. Hautes \'Etudes Sci. Publ. Math. No. {\bf 78} (1993) pp. 187--206 (1994). 
 
\bibitem[S]{S1} G. Soifer, 
{\em Matrix representation of solvable Lie groups and their lattice automorphisms}, Func. Anal. Appl. {\bf 11} (1977) pp. 91-- 92.
 

\bibitem[Th]{Thurston}
W. Thurston, {\em The geometry and topology of three-manifolds}, Princeton lecture notes, 1981. 

\bibitem[Ti]{Tits}
J. Tits, {\em Free subgroups in linear groups}, 
J. Algebra,  {\bf 20} (1972) pp. 250--270.

\bibitem[Tu1]{Tukia1985}
P. Tukia, {\em On isomorphisms of geometrically finite Moebius groups}, Inst. Hautes \'Etudes 
Sci. Publ. Math. No. {\bf 61} (1985) pp. 171--214.

\bibitem[Tu2]{Tukia1994}
P. Tukia, {\em Convergence groups and Gromov's metric hyperbolic spaces}, 
New Zealand J. Math. {\bf 23} (1994), no. 2, pp. 157--187.

\bibitem[Tu3]{Tukia_conical}
P. Tukia, 
{\em Conical limit points and uniform convergence groups}, 
J. Reine Angew. Math. {\bf 501} (1998) pp. 71--98.

\bibitem[V]{Vaisala}
J. V\"{a}is\"{a}l\"{a}, {\em Gromov hyperbolic spaces},  
Expo. Math. {\bf 23} (2005), no. 3, pp. 187--231. 

\bibitem[Y]{Yaman}
A. Yaman, {\em A topological characterisation of relatively hyperbolic groups}, 
J. Reine Angew. Math. {\bf 566} (2004) pp. 41--89. 

\bibitem[Z]{Zhu}
F. Zhu, {\em Relatively dominated representations}, Annales de l'Institut Fourier, Volume {\bf 71} (2021) no. 5, pp. 2169--2235.


\end{thebibliography}
\end{document}